\documentclass[12pt]{amsart}
\usepackage{wrapfig,epsfig,epsf,graphicx,amsmath,amssymb}
\usepackage[normal,footnotesize]{caption2}

\newcommand{\Ass}{A_{\rm ss}}

\newcommand{\B}{\mathcal B}
\newcommand{\e}{\varepsilon}
\newcommand{\h}{h_{\rm ss}}

\newcommand{\hhat}{\hat{h}_{\rm ss}}

\newcommand{\E}{\mathcal E}

\newcommand{\LL}{{\mathcal L}}
\newcommand{\N}{{\text{\scriptsize N}}}
\newcommand{\Pss}{P_{\! \rm ss}}

\newcommand{\R}{\mathbb R}
\newcommand{\T}{\mathbb T}

\renewcommand{\d}{\prime} 
\newcommand{\dd}{{\prime \prime}} 
\newcommand{\ddd}{{\prime \prime \prime}}

\newcommand{\half}{\frac{1}{2}}
\newcommand{\dt}{\Delta t}
\newcommand{\dx}{\Delta x}
\renewcommand{\hbar}{\overline{h}}
\newcommand{\full}[2]{\frac{d #1}{d #2}}

\renewcommand{\part}[2]{\frac{\partial #1}{\partial #2}}

\begin{document}
\title[]
{Heteroclinic orbits, mobility parameters and stability for thin film type equations}
\author[]
{R. S. Laugesen\thanks{laugesen@math.uiuc.edu} and
M. C. Pugh\thanks{mpugh@math.upenn.edu}}
\address{Department of Mathematics, University of Illinois, Urbana, IL
61801}
\address{Department of Mathematics, University of Pennsylvania,
 Philadelphia, PA 19104}
\date \today

\begin{abstract}
We study numerically the phase space of the evolution equation
\[
h_t = -(h^n h_{xxx})_x - \B (h^m h_x)_x .
\]
Here $h(x,t) \geq 0$, $n>0$ and $m \in \R$, and the Bond number $\B$
is positive.

We pursue three goals: to investigate the nonlinear stability of the
positive periodic and constant steady states; to locate heteroclinic
connecting orbits between these steady states and the compactly
supported `droplet' steady states; and to determine how these orbits
change when the `mobility' exponents $n$ and $m$ are changed.

For example, when $n + 1 \leq m <n + 2$ we know from the companion
article that there can be three fundamentally different steady states
with the same period and volume.  The first is a constant steady state
that is a local minimum of the energy.  The second is a positive
periodic steady state that is a saddle for the energy and has higher
energy than the constant steady state.  The third is a periodic
collection of droplet steady states having lower energy than either
the positive or constant steady states.  Here, we find numerically
that the constant steady state appears to be asymptotically stable and
that perturbing the positive periodic steady state in one direction
yields a solution that tends to the constant steady state while
perturbing in the other direction yields a solution that appears to
touch down in finite time.

Also, we consider the effect of changing the mobility coefficients,
$h^n$ and $h^m$.  We change them in such a way that the steady states
are unchanged and find evidence that heteroclinic orbits between steady states are perturbed
but not broken.  We also find that when there appear to be touch--down
singularities, the exponent $n$ affects whether they occur in finite
or infinite time.  It also can affect whether there is one touch--down
or two touch--downs per period.

\end{abstract}

\maketitle

\baselineskip = 12pt

\vspace*{-.2in}
\begin{center}{\sc Contents}\end{center}
%
%
\[
\begin{array}{ccl}
\text{\S\ref{introduction},\S\ref{set_up}} & \text{---} & \text{Introduction and review of steady states} \\
\text{\S\ref{bifurcation}} & \text{---} & \text{Bifurcation diagrams
and weakly nonlinear analysis}\\
\text{\S\ref{simulations}} & \text{---} & \text{Heteroclinic connections}\\
\text{\S\ref{mobility}} & \text{---} & \text{The mobility exponents}\\
\text{\S\ref{numerics}} & \text{---} &
\text{Numerical methods}\\
\text{\S\ref{conclusions}} & \text{---} & \text{Conclusions and future
directions}
\end{array}
\]
\baselineskip = 18pt
\section{Introduction}
\label{introduction}

We study the evolution equation
\begin{equation}  \label{evolve}
        h_t = - (h^n h_{xxx})_x - \B (h^m h_x)_x ,
\end{equation}
where $n>0$ and $m \in \R$, and the Bond number $\B>0$.  This is a
special case in one space dimension of the equation $h_t = - \nabla
\cdot( f(h) \nabla \Delta h) - \nabla \cdot (g(h) \nabla h)$, which
has been used to model the dynamics of a thin film of viscous liquid.
The air/liquid interface is at height $z = h(x,y,t) \geq 0$ and the
liquid/solid interface is at $z=0$. The one dimensional equation
applies if the liquid film is uniform in the $y$ direction.

The fourth order term in equation (\ref{evolve}) reflects
surface-tension-type effects and the second order term can reflect
gravity, van der Waals interactions, thermocapillary effects or the
geometry of the solid substrate, for example.  Typically
\[
f(h) \sim h^n \quad \text{and} \quad g(h) \sim \pm \B h^m
\]
as $h \rightarrow 0$, where $1 \leq n \leq 3$ and $m \in \R, \B > 0$, hence our choice of $f$ and $g$ in (\ref{evolve}) as power
laws. See the introduction of our companion paper \cite{LP3} for
references to the modeling and mathematical literature.

\vskip 6pt

In \cite{LP2} we proved linear stability and instability results for
the positive periodic steady states of (\ref{evolve}). These results
were for zero-mean perturbations that have the same period as the
steady state, or shorter period --- longer perturbations always lead
to linear instability \cite{LP2}. We summarize the linear stability
and instability results in the bifurcation diagrams in
\S\ref{bifurcation}; the main points are roughly that positive
periodic steady states are linearly unstable if $m-n <0$ or $m-n \geq
.80$, are linearly stable if $0 < m-n \leq .75$ (albeit with a zero
eigenvalue arising from the translation invariance of
(\ref{evolve})$\,$), and could be either stable or unstable if $.75 <
m-n < .8$. The case $m-n=0$ is degenerate.

For the linearly {\em stable} steady states, it is natural to ask: do
numerical simulations show asymptotic stability of the steady state
(modulo translation in space), meaning that a solution that starts
from a perturbed steady state will relax back either to the steady
state or to one of its translates? In \S\ref{q1.5} we show the answer
is ``Yes'', for $m$ and $n$ with $m-n =.5$ and other values.

For the steady states $\h$ that are linearly {\em unstable}, two
natural questions occur. First, is the steady state {\it nonlinearly}
unstable? We find this seems always to be the case, as evidenced by
the many numerical simulations in \S\ref{simulations}. Second, do
perturbations of the steady state subsequently converge to some {\it
other} steady state, and can we predict what this long-time limiting
state will be?  This is a difficult question. The three most likely
{\it candidates} for a long-time limit are: a positive periodic steady
state, the constant steady state $\overline{\h}$, or one or more
compactly supported `droplet' steady states. (Droplet steady states
might have zero or nonzero contact angles, but in \cite{LP2,LP3} we
have considered only the zero contact angle case.) How can one predict
which, if any, of these candidates will be the long-time limit?

One way to make predictions is by comparing the {\em energy levels} of
these different kinds of steady state, using the energy
\[
\E(h(\cdot,t)) = \int_0^X \left[ \frac{1}{2} \, h_x(x,t)^2 - 
\frac{\B}{(m-n+2)(m-n+1)} \,
h(x,t)^{m-n+2}
\right] dx .
\]
This energy is known to be dissipated by the evolution, $\frac{d\
}{dt} \E(h(\cdot,t)) \leq 0$, and in our companion article \cite{LP3}
we proved a number of results establishing relative energy levels of
positive periodic, constant and droplet steady states. For example, we
proved in \cite[Theorem~6]{LP3} that if $m-n<0$ or $m-n \geq 1$ then
the positive periodic steady state$\h$ has higher energy than the
constant steady state $\overline{\h}$, suggesting that perturbations
of $\h$ might converge to $\overline{\h}$. Indeed, we find numerically
in \S\ref{simulations} for a variety of $m$ and $n$ values that many
perturbations of $\h$ do subsequently evolve towards the mean. In
another unstable case, if $-2 < m-n < 0$ or $1 \leq m-n < 2$ then
\cite[Theorem~7]{LP3} shows there is a zero angle droplet steady state
with length less than the period of $\h$, with the same area as $\h$,
and with less energy than $\h$. We find numerically in
\S\ref{simulations} that many perturbations of $\h$ evolve towards
touch-down, where the solution goes to zero at a point.  Our code
stops at that time, so we cannot show relaxation to the droplet steady
state, but this certainly seems a likely eventual outcome.

Our theorems and simulations thus allow us to predict at least a
couple of likely long-time limits. But when there is more than one
steady state with lower energy than $\h$, we have not been able to
predict to where a given perturbation of $\h$ will relax. And we are
so far unable to predict the amount of translation that might occur in
the long-time limit (for example, where will the maximum of the
limiting steady state occur?). Impressive results on such a
translation problem have recently been obtained by \cite{bricmont} for
the Cahn--Hilliard equation $h_t = - h_{xxxx} - ((1 - 3 h^2) h_x)_x$
on the whole real line.

\vskip 6pt 

In this paper we also ask how the evolution is affected by changes in
the coefficient functions $h^n$ and $h^m$. In particular, we consider
the evolution of a fixed initial condition under the equation
(\ref{evolve}) for two different pairs $(n_1,m_1)$ and $(n_2,m_2)$ of
exponents. One can show that the steady states of the two evolution
equations are the same if $m_1 - n_1 = m_2 - n_2$.  In that case we
ask: given the same initial data, do the two solutions have the same
long-time limit?  Or, can changing the mobility exponents $n$ and $m$
can change the long-time limit?  We also ask whether changing the
mobility parameters can affect the number and type of finite--time
singularities.

Before outlining the paper, we note that throughout the paper we
consider only {\em zero--mean} perturbations. This seems reasonable
from a physical standpoint, because such perturbations correspond to a
disturbance of the fluid that alters the profile without adding
additional fluid. Mathematically it is reasonable because the
evolution equation (\ref{evolve}) preserves volume for spatially
periodic solutions: $\int h(x,t) \, dx = \int h(x,0) \, dx$ for all
time $t$.  Thus zero--mean perturbations allow the possibility of
relaxation back to the original steady state, while nonzero--mean
perturbations do not.

\subsection*{Outline of the paper}
This paper reports on numerical simulations that are inspired by, and
extend, the theoretical results in our linear stability paper
\cite{LP2} and our companion paper \cite{LP3} on energy stability and
relative energy levels of steady states.  The stability results of
\cite{LP2,LP3} are summarized by the bifurcation diagrams in
\S\ref{bifurcation} below, and we will remind the reader of the
relevant results and their implications when describing our numerical
simulations. So the reader need not digest the earlier papers
\cite{LP2,LP3} before reading this one, although it would help to have
those papers at hand.

\vskip 6pt 

The bifurcation diagrams in Section~\ref{bifurcation} summarize the
known results from \cite{LP2,LP3} on stability of steady states. The
section also contains a weakly nonlinear stability analysis.

Section~\ref{simulations} presents a detailed numerical study of the
evolution equation (\ref{evolve}), for a number of different exponents
$n$ and $m$. We focus especially on initial data close to a steady
state. Our stability and energy level results from \cite{LP2,LP3} lead
to many predictions for the behavior of the solution, both short and
long time, and these predictions are generally borne out by our
simulations.  Strikingly, the period and area of the initial data, the
Bond number $\B$ and the value $m-n$ seem to be sufficient information
to reduce the possible long-time behavior to just a couple of options.

We find evidence for heteroclinic connections between different types
of steady state (such as: periodic to constant, periodic to droplet,
and constant to droplet).  One would like to know how robust the
behaviors observed in our numerical simulations are under changes in
the coefficients $h^n$ and $\B h^m$.  In particular, what happens when
the mobility exponents $n$ and $m$ are changed in a way that $n-m$ is
unchanged, leaving the steady states of the evolution unchanged?  In
Section~\ref{pert_hets} we give numerical evidence that such changes
perturb, but do not break, heteroclinic orbits.

Another question investigated in Section~\ref{mobility} is how
changing $m$ and $n$ affects singularity formation.  For the equation
$h_t = - (h^n h_{xxx})_x$, there has been extensive computational work
studying how the choice of $n$ affects the spatial structure of
singularities and whether they occur in finite or in infinite time
\cite{similar}.  Specifically, simulations suggest there is a critical
exponent $1 < n_c < 2$ such that if $n > n_c$ then solutions are
positive for all time while solutions can touch down in finite time if
$n < n_c$.  For our equation (\ref{evolve}), we numerically estimate
$n_c$ in \S\ref{time_pinch}.  We find that its value depends on the
difference $m-n$.  Then in \S\ref{one_two_sec}, we demonstrate that
the mobility can also affect the {\em number} of apparent
singularities per period.

Section~\ref{numerics} discusses our numerical methods.

Section~\ref{conclusions} sums up the paper, and recalls the `gradient
flow for the energy' interpretation of the evolution equation
(\ref{evolve}).

\section{Terminology, definition of the energy, and review
of steady states}
\label{set_up}

\subsection*{Terminology} 
We write $\T_X$ for a circle of circumference $X>0$. As usual, one
identifies functions on $\T_X$ with functions on $\R$ that are
$X$-periodic and calls them {\it even} or {\it odd} according to
whether they are even or odd on $\R$.

\begin{figure}
  \begin{center}
    \includegraphics[width=.4\textwidth]{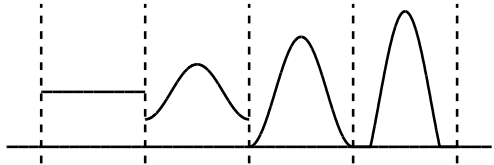}
  \end{center}
  \vspace{-.4cm}  
  \caption{\label{steady} Four types of steady state.}
\vspace{-0.4cm}
\end{figure}
Positive periodic steady states are assumed to satisfy the steady
state equation classically.  A droplet steady state $\h(x)$ (see
Figure \ref{steady}) is by definition positive on some interval $(a,b)$ and
zero elsewhere, with $\h \in C^1[a,b]$; we require $\h$ to satisfy the
steady state equation classically on the open interval $(a,b)$, and to
have equal acute contact angles: $0 \leq \h^\d(a) = - \h^\d(b) <
\infty$. (Throughout the paper, if a function has only one independent
variable then we use ${}^\d$ to denote differentiation with respect to
that variable: $\h^\d = (\h)_x$.)

We say a droplet steady state $\h$ has `zero contact angle' if
$0=\h^\d(a)=-\h^\d(b)$, and `nonzero contact angle' otherwise.  A
`configuration' of droplet steady states is a collection of
steady droplets whose supports are disjoint. 

\subsection*{Definition of the energy} 
There is a well-known dissipated energy (or Liapunov function) for the
evolution equation (\ref{evolve}).  It is defined for $\ell \in
H^1(\T_X)$ to be
\begin{equation} \label{energy_is}
\E(\ell) = \int_0^X \left[ \frac{1}{2} (\ell^\d)^2 - H(\ell) \right] dx ,
\end{equation}
where $H$ is a function with $H^\dd(y) = \B y^{m-n}$. This energy is
strictly dissipated: if $h(x,t)$ is a smooth solution of
(\ref{evolve}) then $(d/dt) \E(h(\cdot,t)) \leq 0$, with equality if
and only if $h$ is a steady state (cf.\ \cite[\S2.1]{LP3}). The
energy, like the evolution equation, is invariant under translation:
hence it cannot distinguish between a steady state and its translates.

\vskip 6pt

A precise definition of linear stability can be found in
\cite[Appendix~A]{LP2}, and of energy stability in \cite[\S2]{LP3},
but those definitions are not needed to understand this paper.

\subsection*{Brief review of steady states}
\label{power_law_review}
Here we quickly review the basic facts about steady states needed to
appreciate this paper. For more on the steady states and their
properties, and for justifications of the following remarks, see
\cite[\S2.3]{LP3} and \cite{LP1} and the references therein.

\vskip 6pt 
We start with a non-constant positive periodic steady state $\h \in
C^4(\T_X)$ of the evolution (\ref{evolve}). We translate $\h$ so that
its minimum occurs at $x=0$. The steady state equation for
(\ref{evolve}) integrates to give $\h^n \h^\ddd + \B \h^m \h^\d = C$
for some constant $C$.

The constant $C$ (the flux) equals zero, by integrating $\h^\ddd + \B
\h^{m-n} \h^\d = C \h^{-n}$ over a period. Hence the steady state
satisfies
\begin{equation} \label{ss1}
\h^\ddd + r(\h) \h^\d = 0 ;
\end{equation}
here $r(y) := \B y^{q-1}$ and 
\[
\fbox{$q := m-n+1$.}
\]
This exponent $q$ determines many properties of the
steady state, including (usually) its linear stability.

Integrating, we find the steady states have a nonlinear oscillator formulation:
\begin{equation} \label{ss3}
  \h^\dd + \frac{\B \h^q - D}{q} = 0
\end{equation}
for some constant $D$, when $q \neq 0$. For $q=0$ the analogous
equation is $\h^\dd + \B \log \h - D = 0$.  This oscillator equation
contains three constants: $q$, $\B$, and $D$, but $\B$ and $D$ can be
removed by rescaling $\h$ and $x$ (see \cite[eq.\ (7)]{LP3}), leading
to the `canonical' steady state equation
\begin{eqnarray} \label{ss4}
k^\dd + \frac{k^q - 1}{q} & = & 0, \qquad q \neq 0, \\ 
\label{ss5}
        k^\dd + \log k            & = & 0, \qquad q   =  0 .
\end{eqnarray}

Every positive periodic steady state $\h$ can be rescaled to such a
function $k=k_\alpha$ with $k_\alpha^\d(0)=0$, where we write $\alpha=k(0) \in (0,1)$ for the initial value. Conversely, for each $q
\in \R$ and $\alpha \in (0,1)$ there exists a unique smooth positive
periodic $k_\alpha$ satisfying equations (\ref{ss4}--\ref{ss5}) and
with $k_\alpha(0)=\alpha, k_\alpha^\d(0)=0$. The same holds for
$\alpha=0$ when $q>-1$, although $k_0$ may be only $C^1$-smooth at
$x=0$, where $k_0(0)=0$ (see \cite[Theorem~3.2]{LP1}).

We write
\[
        P=P_\alpha=P(\alpha) \qquad \text{and} \qquad 				A=A_\alpha=A(\alpha)
\]
respectively for the least period of $k_\alpha$ and for the area under
its graph, $A=\int_0^P k_\alpha(x) \, dx$. Note that $P_\alpha$ and
$A_\alpha$ approach $2\pi$ as $\alpha \rightarrow 1$. As seen in
\cite{LP2,LP3}, the function
\[
E(\alpha) := P(\alpha)^{3-q} A(\alpha)^{q-1} = P(\alpha)^2 [A(\alpha)/P(\alpha)]^{q-1}
\]
plays a large role in the stability theory of the steady states, in
part due to its rescaling invariance:
\begin{equation} \label{invariant}
\B \Pss^{3-q} \Ass^{q-1} = P(\alpha)^{3-q} A(\alpha)^{q-1} (=E(\alpha)),
\end{equation}
where $\Pss$ and $\Ass$ are the period and area of the steady state
$\h$.
\section{Bifurcation diagrams and weakly nonlinear analysis}
\label{bifurcation}

In this section we present bifurcation diagrams that encode the linear
stability information in \cite{LP2,LP3}. Then we relate our
bifurcation diagrams to a weakly nonlinear analysis near the constant
steady state.

\subsection{Bifurcation diagrams}
Figure~\ref{bifurcs} gives bifurcation diagrams for various values of
$q$.  These $q$-values are representative members of the intervals
$(-\infty,1)$, $(1,1.75)$, $(1.75,1.794)$, and $(1.795,\infty)$.  (The
value $1.794$ approximates a critical exponent; see \cite[\S5.1]{LP1}
for further details.)

In all these diagrams we are considering the $2 \pi$-periodic problem
with $\B = 1$.  We construct the diagrams as follows.  Given $q$, we
first compute the rescaled steady states $k_\alpha$ for a range of
$\alpha \in (0,1)$.  For each $k_\alpha$, we use its period $P_\alpha$
to determine a constant $D$ in the rescaling \cite[eq.\ (7)]{LP3} that
yields a positive steady state $\h$ of period $\Pss=2 \pi$.  We then
plot the amplitude of $\h$ versus the scale invariant quantity $E = \B
\Pss^{3-q} \Ass^{q-1}$, determining the linear stability of the steady
state by the results in \cite[\S3.2]{LP2}, particularly
\cite[Theorem~9]{LP2}. The conclusions on linear stability summarized
in the figure are rigorously proved in \cite{LP2} except for $1 < q <
2$, in which case we are relying on a combination of analytical and
numerical results.

\begin{figure}[h] 
 \vspace{-.5cm}
    \begin{center}
    \includegraphics[width=.2\textwidth]{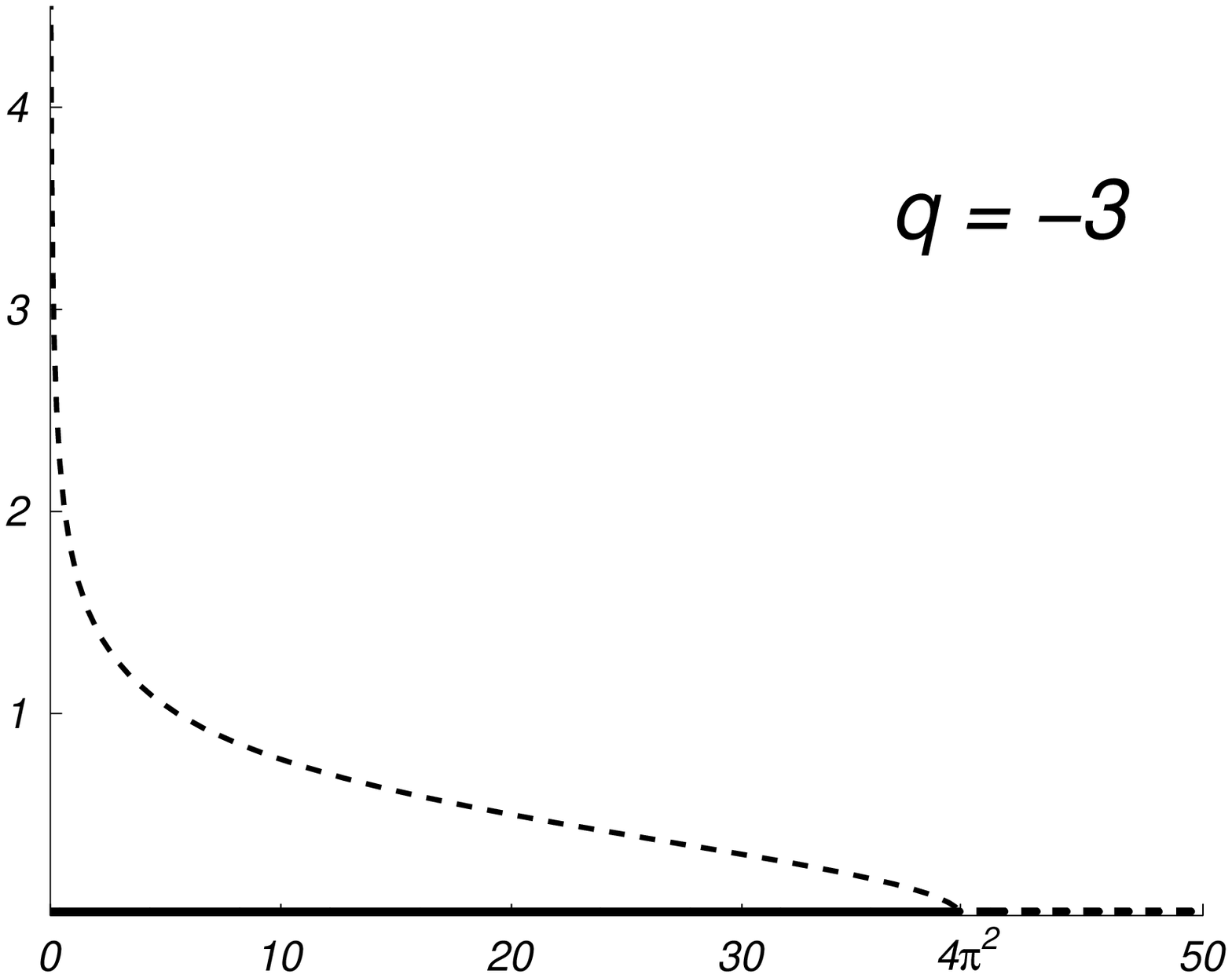}
    \hspace{.5cm} 
    \includegraphics[width=.2\textwidth]{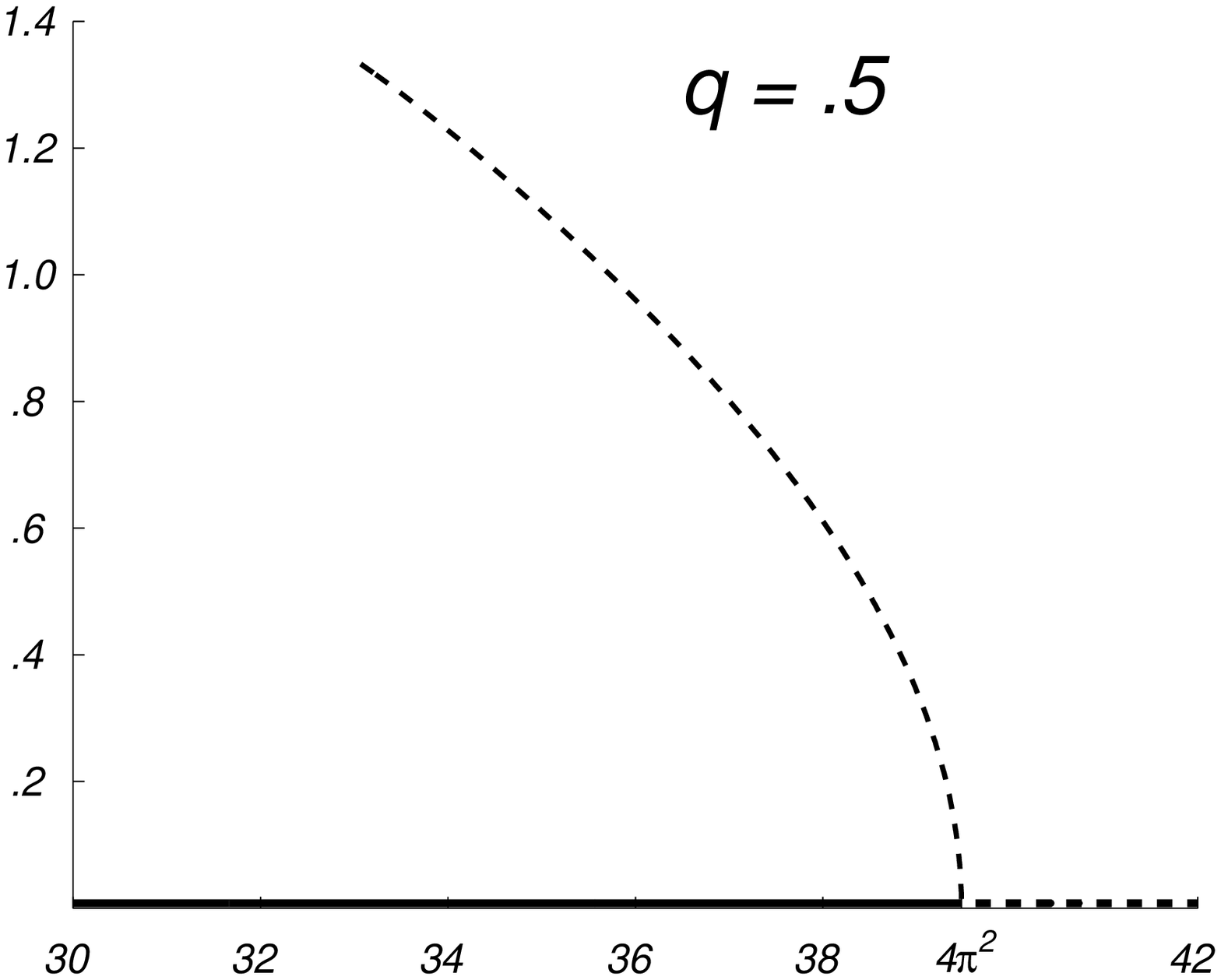}
    \hspace{.5cm} 
    \includegraphics[width=.2\textwidth]{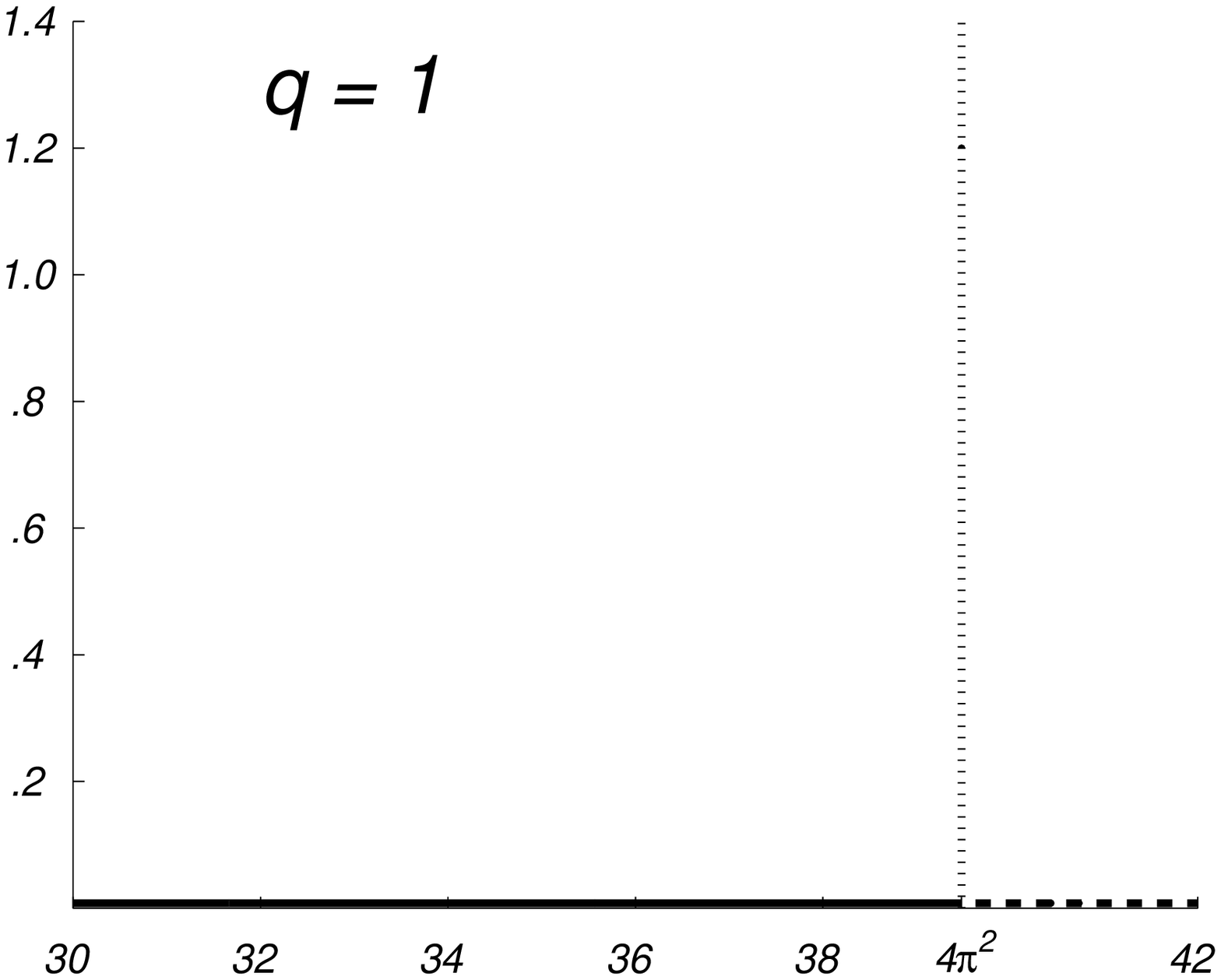}
    \end{center}
 \vspace{-.5cm}
\end{figure}
\begin{figure}[h] 
 \vspace{-.4cm}
    \begin{center}
    \includegraphics[width=.2\textwidth]{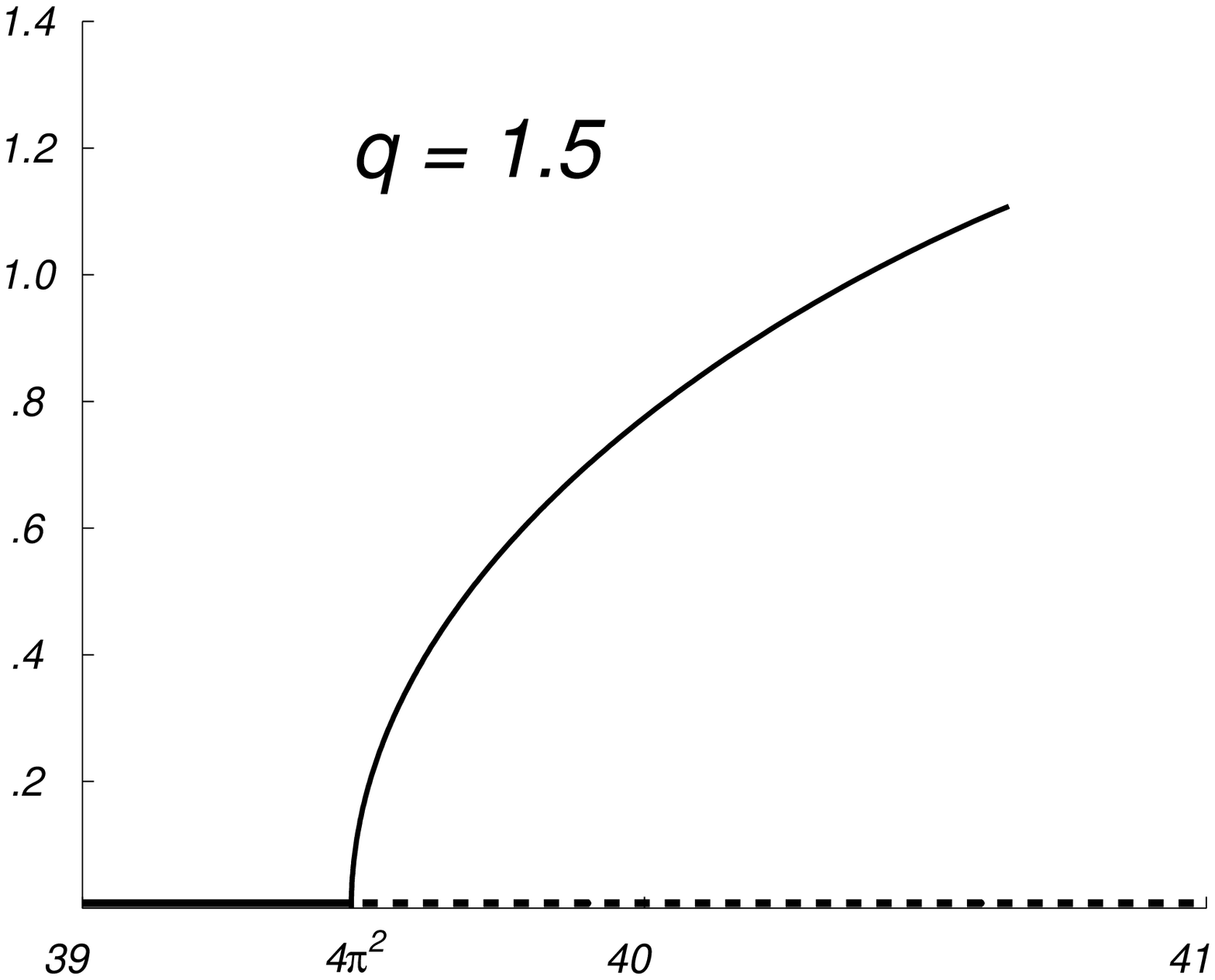}
    \hspace{.5cm} 
    \includegraphics[width=.2\textwidth]{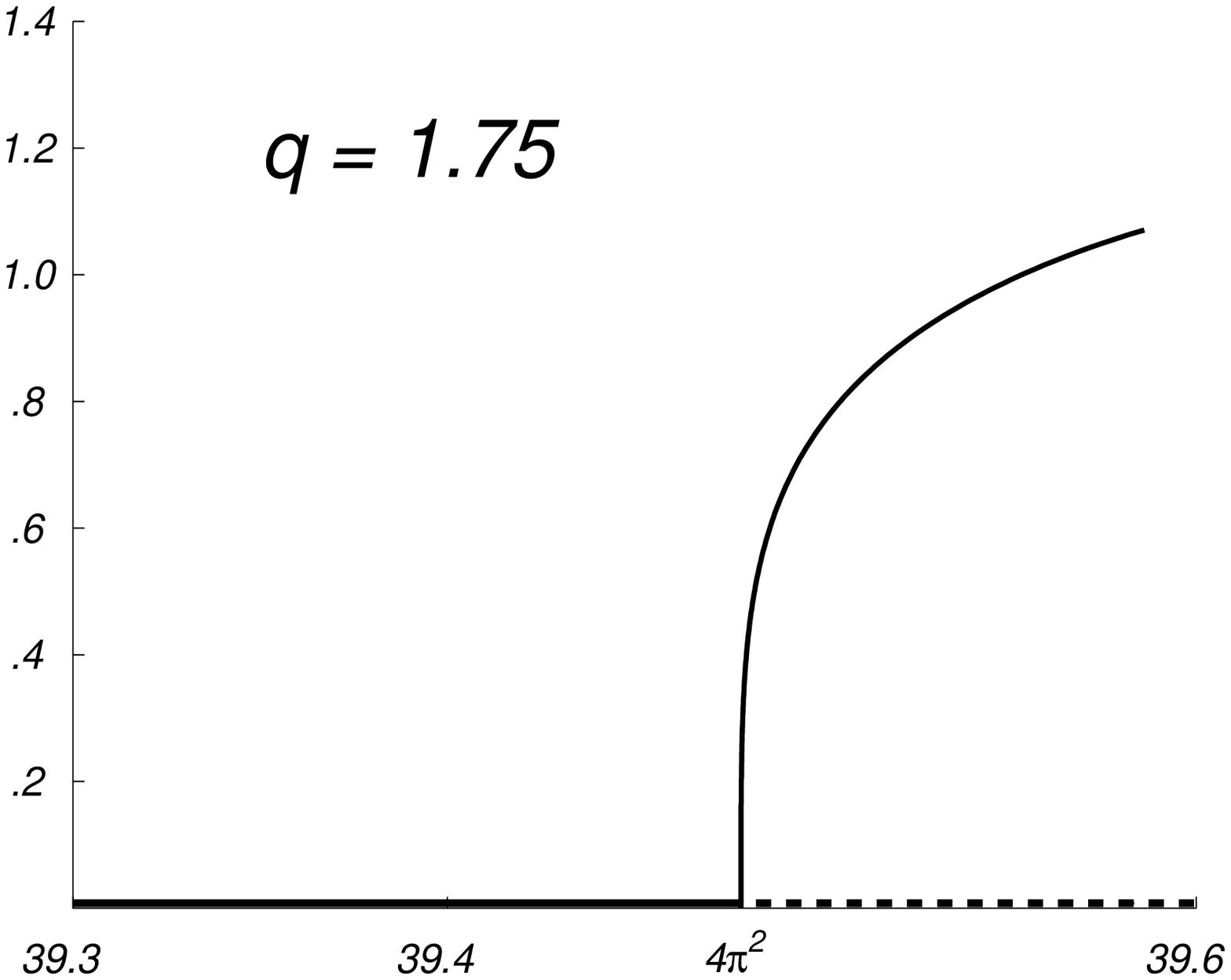}
    \hspace{.5cm} 
    \includegraphics[width=.2\textwidth]{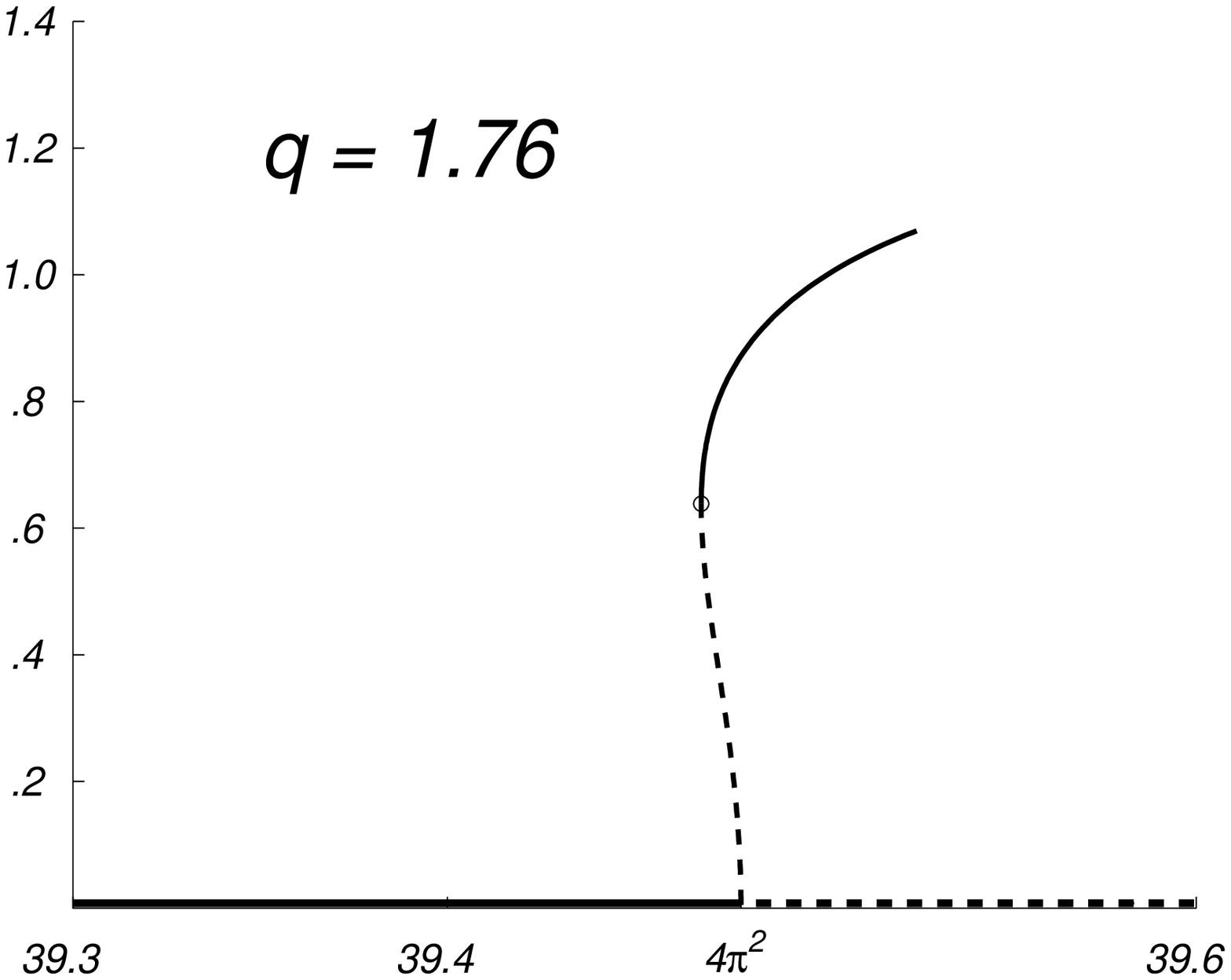}
    \end{center}
 \vspace{-.5cm}
\end{figure}
\begin{figure}[h] 
 \vspace{-.4cm}
    \begin{center}
    \includegraphics[width=.2\textwidth]{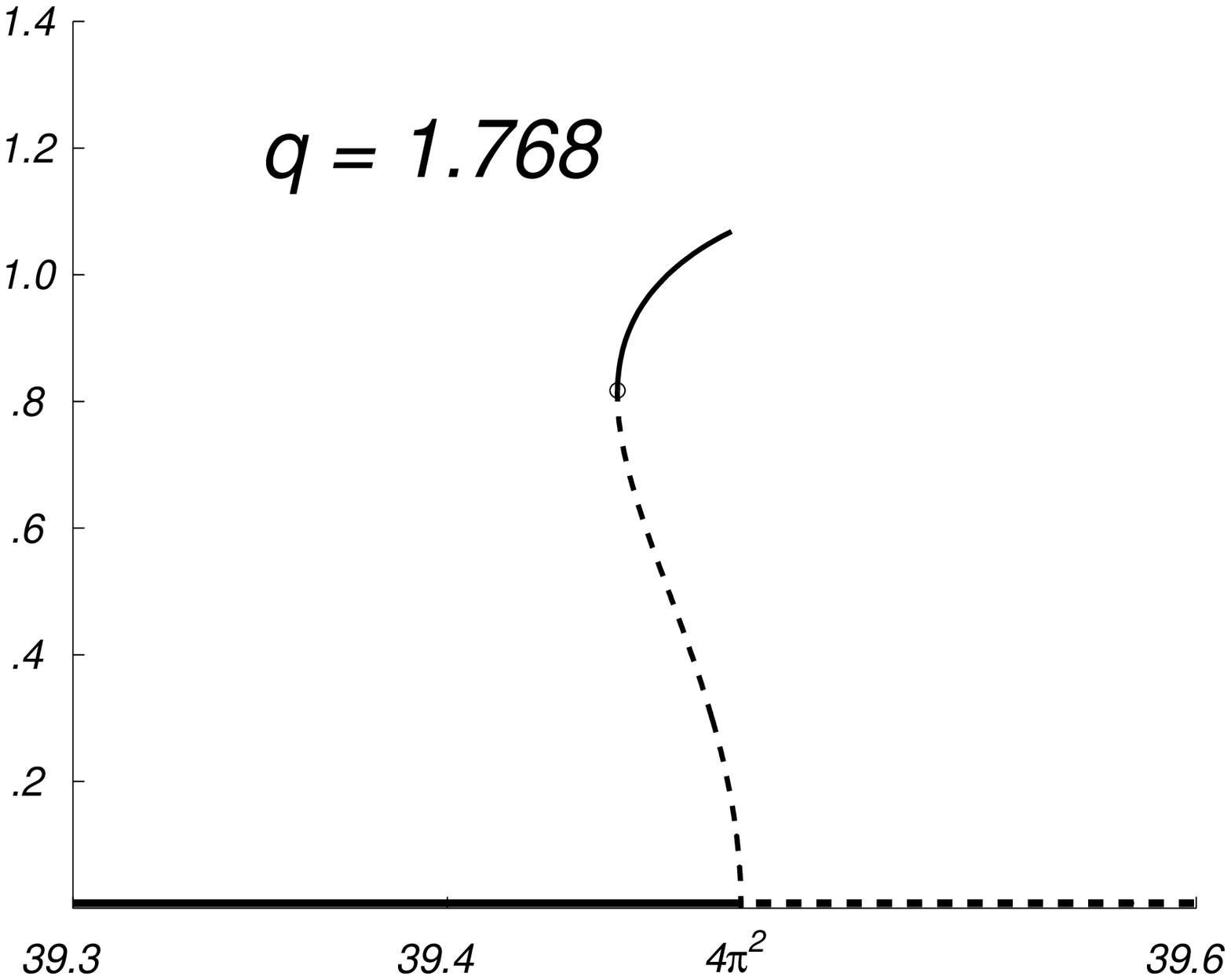}
    \hspace{.5cm} 
    \includegraphics[width=.2\textwidth]{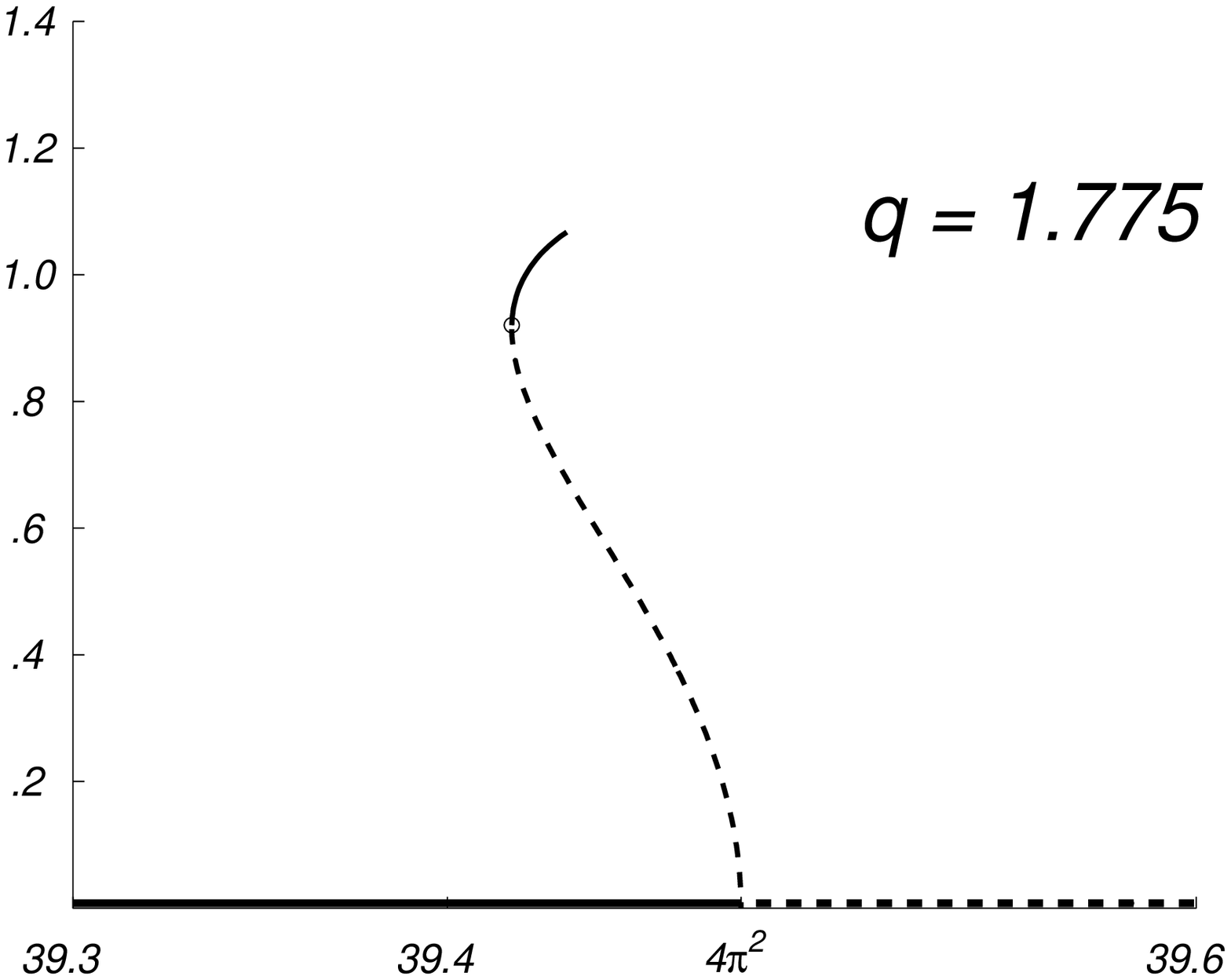}
    \hspace{.5cm} 
    \includegraphics[width=.2\textwidth]{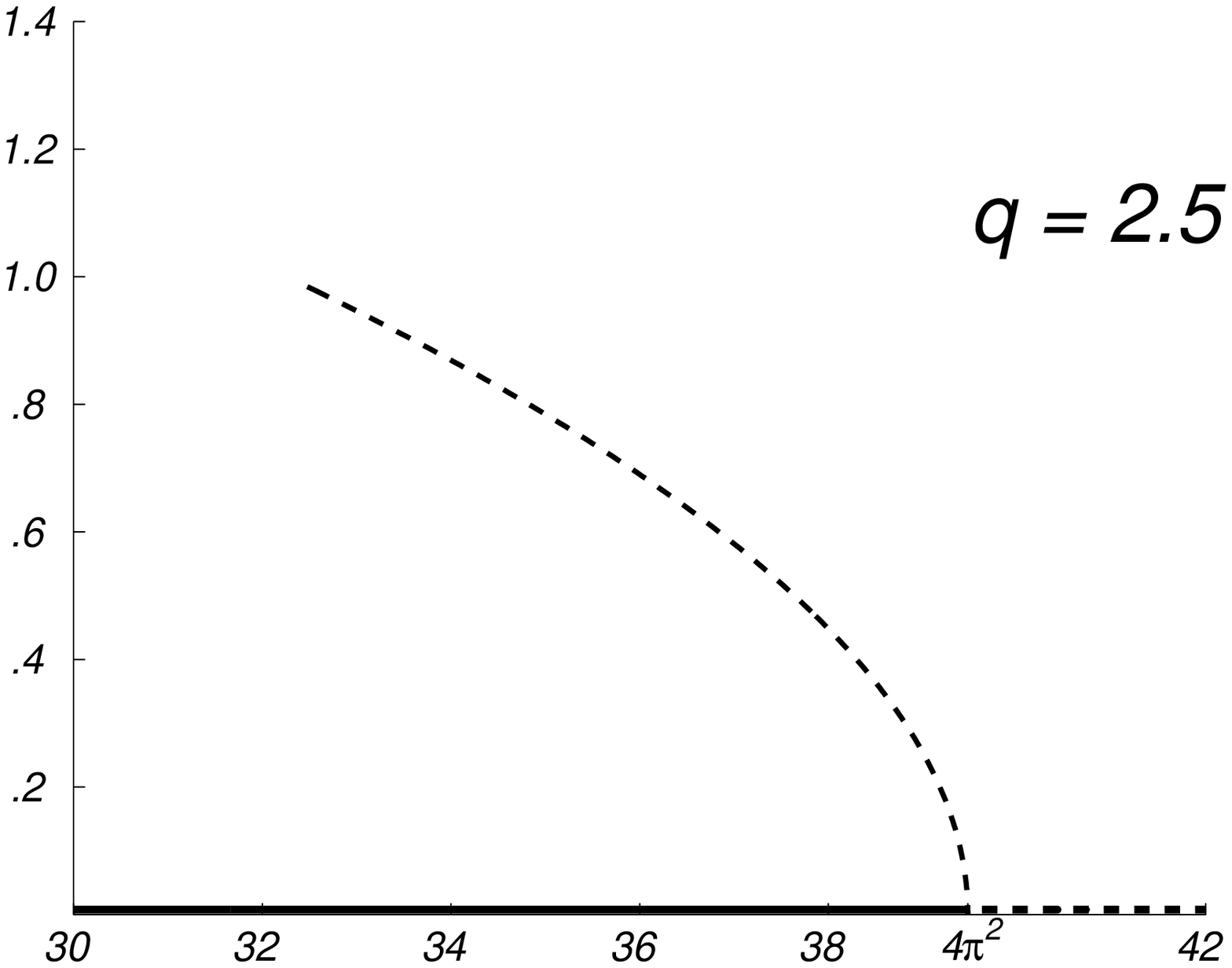}
    \end{center}
 \vspace{-.25cm}
\renewcommand{\figurename}{Fig.$\!$$\!$}
\setcaptionwidth{5in}
\caption{\label{bifurcs} The $x$-axis is $\B \Pss^{3-q} \Ass^{q-1}$
and $y$-axis is $(h_{max}-h_{min})/2$.  Dashed: linearly unstable;
dotted: linearly neutrally stable; solid: linearly stable.}
\end{figure}

The horizontal axes of these diagrams also show the linear stability
of the {\em constant} steady state, with respect to zero-mean
perturbations, where the constant steady state is considered as having
period $\Pss$ and area $\Ass$. This linear stability information is
taken from \cite[Theorem~10]{LP3}.

Qualitatively, the diagrams change continuously as $q$
increases. Choosing a period other than $2 \pi$ or a Bond number other
than $\B = 1$ simply dilates the $y$-axis of the diagrams.
%
%
We discuss the diagrams further when we present numerical simulations in
\S\ref{simulations}.

\subsection{Weakly nonlinear analysis}
In the following, we sketch the weakly nonlinear analysis for the
evolution equation (\ref{evolve}). This was done for $q=-3$ in \cite[\S2.3]{WB00}. 

In short, we consider physical parameters such that the constant
steady state has one mode which is barely linearly stable or is barely
linearly unstable, while all other modes are strongly damped.  This
yields a separation of timescales which allows one to find a reduced
representation of the PDE in terms of an ODE governing the amplitude
of the unstable mode.  We refer the reader to \S5.1 of
\cite{Manneville}.

Let $h$ be an $X$-periodic solution of $ h_t = - (h^n h_{xxx})_x - \B
(h^m h_x)_x $ with mean value $\hbar$.  We rescale the solution to
have period $2 \pi$ and mean value $1$, and we also rescale time:
$$
\zeta = \frac{2 \pi x}{X},
\qquad
\hbar \eta = h,
\qquad \mbox{and} \qquad
t' = \left( \frac{2 \pi}{X} \right)^{\! \! 2} t.
$$
The rescaled evolution equation is $ \eta_{t'} = - \sigma( \eta^n
\eta_{\zeta \zeta \zeta})_\zeta - \mu(\eta^m \eta_\zeta)_\zeta $ where
$$
\sigma = \left(\frac{2 \pi}{X}\right)^{\! \! 2} \hbar^{\,n} 
\qquad \text{and} \qquad
\mu = \B \hbar^{\,m}.
$$
Linearizing about $\hbar = 1$, we find that linearly unstable modes
exist if $0 < k < k_c$.  Proceeding in the usual manner, we introduce
a small parameter $\delta$ that corresponds to moving through the
critical wave number $k_c = 1$:
$$
\frac{\mu}{\sigma} = 1 + Q \delta^2
$$
where $Q = \pm 1$.  We then introduce a slow time-scale $\tau =
\delta^2 t'$ and expand the solution in orders of $\delta$: $
\eta(\zeta,\tau) = 1 + \delta \eta_1(\zeta,\tau) + \delta^2
\eta_2(\zeta,\tau) + \delta^3 \eta_3(\zeta,\tau) + O(\delta^4).$ For
simplicity, we assume the solution is even.  By the usual arguments,
$\eta_1(\zeta,\tau) = A(\tau) \cos(\zeta)$, $\eta_2(\zeta,\tau) =
B(\tau) \cos(2 \zeta)$, and $\eta_3(\zeta,\tau) = C(\tau) \cos(3
\zeta)$.  Putting this ansatz into the evolution equation and
expanding in orders of $\delta$, we find there are no $O(1)$ or
$O(\delta)$ terms.  At $O(\delta^2)$ and $O(\delta^3)$, one determines
the amplitudes $B(\tau)$ and $C(\tau)$ in terms of $A(\tau)$ which, in
turn, satisfies
$$
\full{A}{\tau} = Q \sigma A(\tau) - \kappa A(\tau)^3
\qquad \mbox{where} \qquad
\kappa = \frac{\sigma}{6} (q-1) ( 1.75 - q ).
$$
The dynamics of the amplitude $A(\tau)$ depend on the signs of the
Landau constant $\kappa$ and of the linear term.  If $\kappa > 0$ then
for $Q = 1$ the constant $A(\tau) \equiv 0$ is linearly unstable and
$A(\tau)$ saturates to the linearly stable amplitude $A_c(\tau) \equiv
\sqrt{\sigma/\kappa}$.  This corresponds to a supercritical
bifurcation.  If $\kappa < 0$ then for $Q = -1$ the constant $A(\tau) \equiv
0$ is linearly stable and the steady amplitude $A_c(\tau) \equiv
\sqrt{-\sigma/\kappa}$ is linearly unstable.  This corresponds to a
subcritical bifurcation.  Since $\sigma > 0$,
$$
1 < q < 1.75 \Longrightarrow \mbox{{\it supercritical
bifurcation}}
\quad \mbox{and} \quad
q < 1, 1.75 < q \Longrightarrow \mbox{{\it subcritical
bifurcation}}.
$$

The Landau constant $\kappa$ is determined by the mean value $\hbar$,
the period $X$, and the exponents $m$ and $n$, but not the Bond number
$\B$.  Subcritical bifurcations are often seen in systems that can
have finite-time pinching (rupture) singularities {\it e.g.} \cite[\S3.2]{BBW}, \cite[\S{IV}]{OB}.

The above weakly nonlinear analysis is consistent with our analytical and
computational results.  Specifically, recalling the bifurcation
diagrams in Figure~\ref{bifurcs} we note that positive periodic steady
states are linearly unstable for $q < 1$ and $q \geq 2$. For $1.75 < q
< 2$, when there is only one positive periodic steady state it is
linearly unstable.  When there are two, the one of smaller amplitude
$(h_{max}-h_{min})/2$ is linearly unstable.  For $1 < q \leq 1.75$ the
positive periodic steady state is linearly stable.

\section{Simulations, and heteroclinic orbits connecting steady states}
\label{simulations}
We now numerically simulate solutions of the evolution equation
(\ref{evolve}), for a wide range of initial data near steady
states.  This has not been done before. The solutions obtained display
a great variety of stability and long--time behaviors.  Our stability
theorems \cite{LP2,LP3} often allow us to {\em predict} the
numerically observed short--time behavior, and our theorems on the
energy levels of steady states \cite{LP3} often allow us to guess the
long--time limit of the evolution. As part of this we predict (and
find strong evidence for) heteroclinic connections between certain
steady states.

We expect our numerical investigations of the power law evolution
(\ref{evolve}) will provide resources, ideas and motivation for
researchers studying $h_t = - (f(h)h_{xxx})_x - (g(h) h_x)_x$ with
non-power law coefficient functions $f$ and $g$. There are some such numerical studies already. For example, the papers \cite{Gruen00,OB} consider an $f$ that is degenerate ($f(0)=0$) and $g$'s that are not power laws, and there is a large literature on the Cahn--Hilliard equation (for which $f \equiv 1$ and $g$ is a quadratic).

\vskip 6pt
Recall the steady states depend on the parameter
\[
q = m-n+1.
\]
We take seven values of $q$: 
\[
q = -3, \; 0.5, \; 1, \; 1.5, \; 1.768, \; 2.5, \; 4,
\]
representatives of the intervals
$\{(-\infty,-1],(-1,1),(1,1.75],(1.75,1.79),(1.8,3),[3,\infty) \}$ in
which our theorems suggest the solutions will display distinct
behaviors.  For the $q=-3$ case, we take $n=3$ and $m=-1$, making
(\ref{evolve}) a `van der Waals' equation previously studied by other
authors.  Otherwise, we take $n=1$ and $m=q$.  In any event, we find
in \S\ref{mobility} that our numerical simulations are not greatly
affected if we change $n>0$ and $m$ in a manner that keeps $q$ fixed
({\it i.e.} that keeps $m-n$ fixed).
%


\vskip 12pt
\subsection{$\mathbf{q=-3}$: the van der Waals case}\ 
\label{qm3}

{\it Characteristic features for $q \in (-\infty,-1]$: 
positive periodic steady states are linearly unstable, and there are no
droplet steady states with acute contact angles.} (See bifurcation
diagram~\ref{bifurcs}a and \cite[\S2.2]{LP1}.)

\subsubsection{$q=-3$. Perturbing the positive periodic 
steady state} 
\label{qm3pp}

First we explain how to find a steady state having period $\Pss = 2
\pi$ and having some specified value for the area (or volume)
$\Ass=\int_0^{2 \pi} \h(x) \, dx$.

Given the Bond number
$\B$ (a physical parameter), a positive steady state satisfies (\ref{ss3}):
\[
   \h^\dd + \frac{\B \h^q - D}{q} = 0
\]
for some constant $D$.  The problem is to find the value of $D$ for
which there is a steady state of period $2\pi$ and area $\Ass$.  For
$q \leq -1$, if $\Ass > 2 \pi \, \B^{1/(1-q)}$ then there exists a
non-constant positive periodic steady state $\h$ with period $2\pi$
and area $\Ass$ (see \cite[\S5.1]{LP1}).  The rescaling \cite[eq.\
(7)]{LP3} then implies there is an admissible value of $D$.

For simplicity, instead of starting with the Bond number $\B$ and
finding this admissible values of $D$, we instead fix $D = 1$ and
determine the interval of admissible $\B$ values. That is, we choose
$\alpha \in (0,1)$ and consider the steady state $k_\alpha$ of the
rescaled equation, as in \S\ref{power_law_review}.  Its period
$P(\alpha)$ together with $D=1$ and $\Pss = 2 \pi$ then uniquely
determine $\B$ by \cite[eq.\ (24)]{LP3}, and hence the steady state
$\h$ by \cite[eq.\ (7)]{LP3}.  (Choosing a different $\alpha$ would
yield a different $\B$ and $\h$.)  Finally, we locate a nearby
`finite-difference steady state' on $\N$ meshpoints (see \S
\ref{compute_id}) and study numerically its stability under
perturbation.

For $q = -3$ we carry out this construction with $\alpha =0.2145$ and
period $P_\alpha = 16.32$,
resulting in Bond number $\B = 0.003259$,
and a non-constant positive periodic finite-difference steady state
$\h$ of least period $2 \pi$ and with area $A = 6.884$.
Note the positive periodic steady state is linearly unstable, by
bifurcation diagram~\ref{bifurcs}a with $\B \Pss^{3-q} \Ass^{q-1} =
.08930$.

\vskip 3pt
\noindent {\bf Remark:} For the rest of \S\ref{simulations}, whenever
we refer to a `positive periodic steady state', we implicitly mean a
$2\pi$-periodic {\it finite-difference steady state} that has its
minimum at $x=0$ (see \S\ref{compute_id}).  Also, we always find our
periodic steady states as above, by fixing $D=1$ then choosing $\alpha
\in (0,1)$ and then determining $\B$ and $\Ass$.  The only exception
is for the $q = 1.768$ simulations discussed in \S\ref{q1.768}.

\vskip 3pt 
Now that we have a steady state for $q=-3$, we use a
perturbation of it as initial data in the equation
\begin{equation} \label{vanderwaal}
h_t = - (h^3 h_{xxx})_x - \B (h^{-1}h_x)_x.
\end{equation}
Here $n=3,m=-1$ and $q=m-n+1=-3$. We study this equation for the rest
of \S\ref{qm3}. It was proposed by Williams and Davis \cite{WD82} to
model a thin liquid film with net repulsive van der Waals
interactions, and more recently it has been studied by Zhang and
Lister \cite{ZL99} and by Witelski and Bernoff \cite{WB99,WB00}.

\noindent $\bullet$ {\it Even perturbations.}

First we perturb $\h$ with the even zero-mean perturbation $\pm \e
\h^\dd$ and numerically study the resulting solution.

The steady state $\h$ is linearly unstable and since the perturbation
$+ \e \h^\dd$ lowers the maximum of $\h$ and raises the minimum, one
might hope the resulting solution would converge to the constant
steady state $h \equiv \overline{\h}$. If this happens for all small
$\e$, then this would be strong evidence for existence of a
heteroclinic orbit connecting $\h$ to the constant steady state.
There are a number of theoretical reasons to suspect such heteroclinic
orbits exist: (i) $\h$ is energy unstable in the directions $\pm
\h^\dd$ by \cite[Theorem~2]{LP3}, (ii) the energy of $\h$ is higher
than that of the constant steady state $\overline{\h}$ by
\cite[Theorem~6]{LP3} (also observed numerically by Witelski and
Bernoff \cite[\S3]{WB00}), and (iii) the constant steady state is a
local minimum of the energy $\E$ by \cite[Theorem~10]{LP3}. [To deduce
(iii) from \cite[Theorem~10]{LP3} requires $\B \overline{\h}^{\, q-1}
X^2 < 4 \pi^2$, which we now establish: $\B \overline{\h}^{\, q-1} X^2
= \B X^{3-q} \Ass^{q-1} = E(\alpha)$ for some $\alpha \in (0,1)$, by
rescaling as in (\ref{invariant}), and $E(1)=4\pi^2$. So we want to
show $E(\alpha)<E(1)$; this holds because $E^\d>0$ by
\cite[Theorem~11]{LP2}, since $q<1$.]

To seek evidence for a heteroclinic connection, we start with initial
data $\h + 10^{-4} \h^\dd$.  Here, we've normalized $\h^\dd$ to have
$L^\infty$ norm $1$; all perturbations are similarly normalized.  We
have also considered $\e$ smaller than $10^{-4}$, in most of the
simulations below; we found that these smaller perturbations resulted
in qualitatively the same behavior, until one hits the level of
roundoff error.

\begin{figure}[h]
\begin{center}   
    \epsfig{file=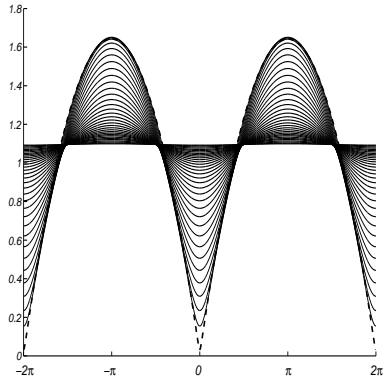, width=2in, height=2in}
\end{center}
\renewcommand{\figurename}{Fig.$\!$$\!$}   
\setcaptionwidth{5.7in}
\vspace*{-.2cm}
\caption{\label{qm3_relax} $q=-3,n=3$. Dashed: initial data $\h
+10^{-4} \h^\dd$.  Local extrema stay fixed in space and solution
relaxes to the mean.}
\end{figure}
The steady state $\h$ has very large curvature at its local minima,
and so we need a large number of meshpoints to resolve the initial
data $\h + 10^{-4} \h^\dd$ with spectral accuracy.  We find that for a
solution on $[0,2\pi)$ we need $2048$ meshpoints: $\h(x) =
\sum_{k=-1023}^{1023} a_k \exp(i k x)$ has amplitudes that decay to
the level of round--off error ($a_k \sim 10^{-13}$ for $k \sim 1024$).
Figure~\ref{qm3_relax} shows the evolution of the solution; the
solution relaxes to the constant steady state. (This was shown
previously in \cite[Figure 4b]{WB00}.) We use the adaptive
timestepping described in \S\ref{acc}.  At first sight, this would
seem unnecessary since the solution is becoming more regular as it
evolves.  However, there is a short transient during which
$h_{min}(t)$ decreases.  The timestep initially decreases to
accurately track this fast behavior and then increases as the solution
relaxes to the mean.  We use adaptive time--stepping throughout our
work since we almost always observe such a short transient.  Also, in
a number of cases, the solutions become less smooth (curvatures
increase) as time passes, requiring refinement later in time.

\vskip 3pt

In the opposite direction, the perturbation $- \e \h^\dd$ raises the
maximum and lowers the minimum of the steady state $\h$.  Since $\h$
is energy unstable, we might expect the solution to subsequently
converge to a droplet steady state or to a configuration of droplet
steady states.  From \cite[\S2.2]{LP1}, if such a droplet exists it
must have $90^\circ$ contact angles, though we have not discussed such
steady states here in this paper or in \cite{LP1,LP2}.

Our numerical simulation of the solution with initial data $\h -
10^{-4} \h^\dd$ shows that, after a short transient, the minimum
height of the solution decreases in time, appearing to decrease to
zero in finite time.
As the minimum height decreases, the curvature increases, requiring
that after some time the number of meshpoints be increased to keep the
solution spectrally resolved.  We do this as follows.  We compute the
solution with $2048$ meshpoints until the computation stops
($h_{min}(t) = 0$, see \S\ref{stop}.) We look at the power spectrum of
the solution and choose a time right before the active part of the
power-spectrum is reaching the Nyquist frequency (see right plot of
Figure~\ref{n8192_ev}). That is, we find the last time at which the
$1023$rd Fourier amplitude of the solution is at the level of
round--off.  We take the solution at this time and compute its Fourier
coefficients, defined for wave numbers $-\N/2+1 \leq k \leq \N/2-1$
where $\N = 2048$.  We pad by zeros, extending the Fourier
coefficients to be defined for wave numbers $-\N+1 \leq k \leq \N-1$,
and then compute the inverse Fourier transform.  This yields a
function on $2\N$ meshpoints that is indistinguishable from the
solution at that time, to the level of round--off.  Using this
function as initial data, we continue the computation on $2 \N$
meshpoints, repeating this point--doubling process whenever the
solution becomes unresolved.

\begin{figure}[h] 
  \vspace{-.2cm}
    \begin{center}
      \includegraphics[width=.35\textwidth]{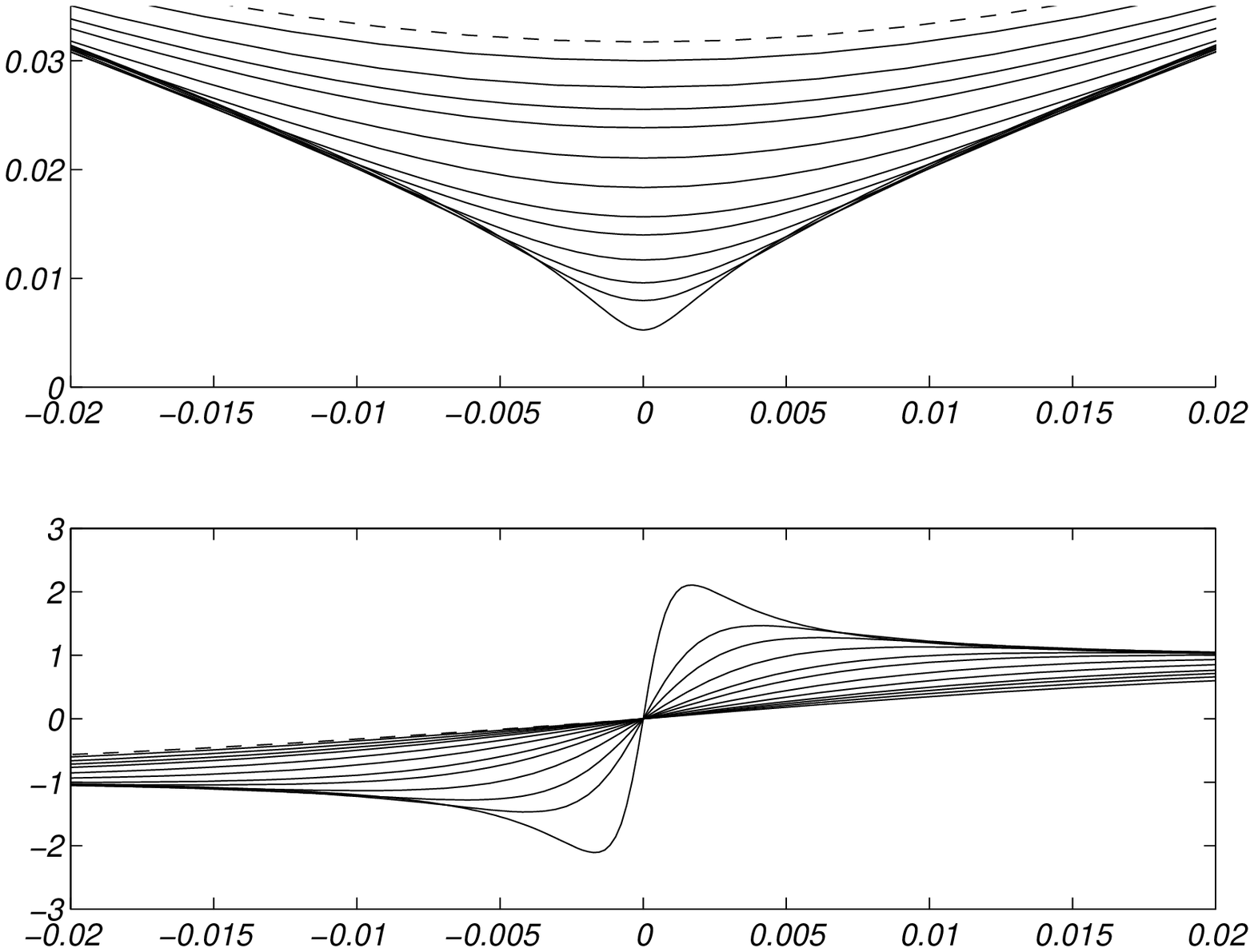}
      \hspace{1.5cm}      
      \includegraphics[width=.35\textwidth]{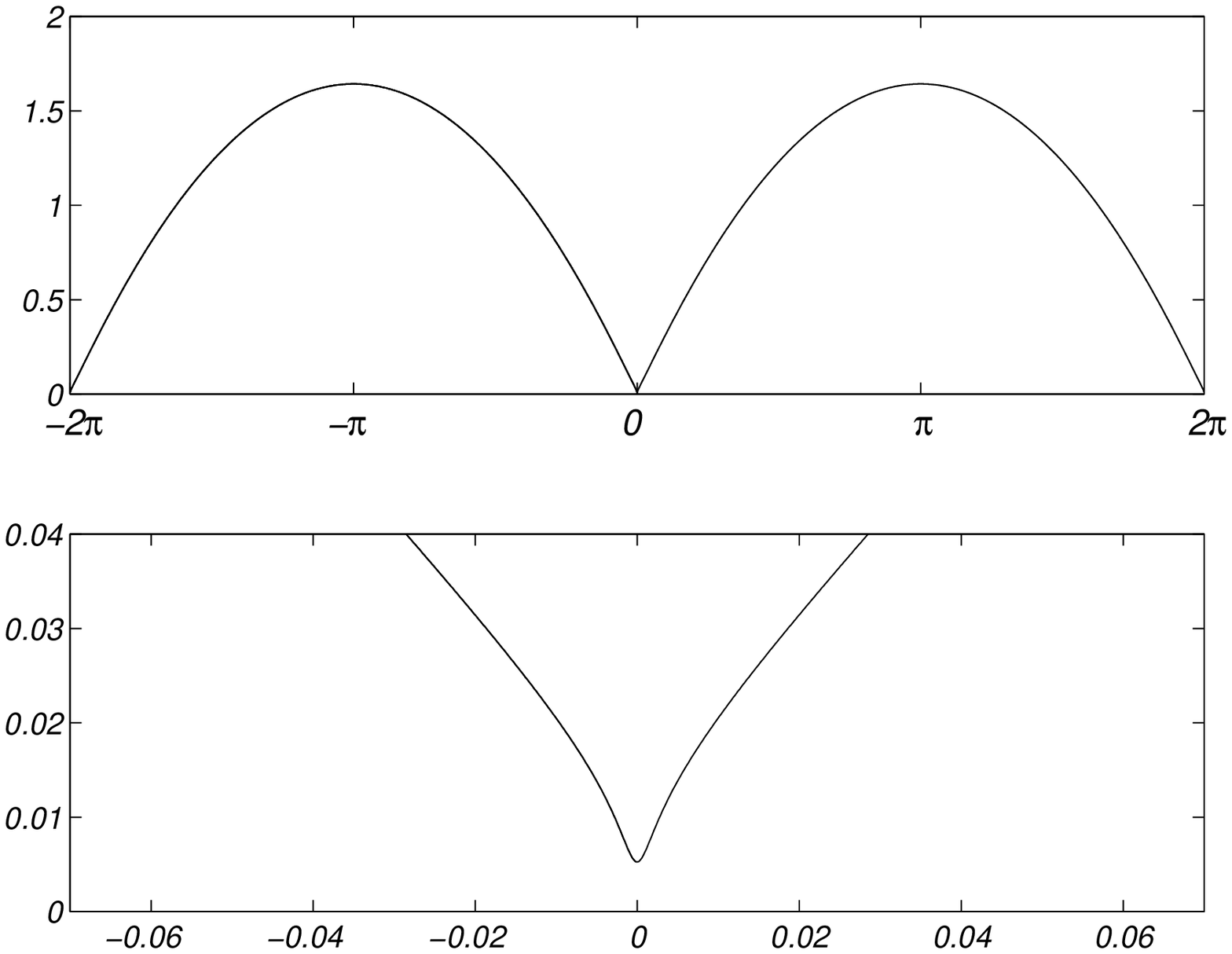}
    \end{center}
  \vspace{-.2cm}	  
\renewcommand{\figurename}{Fig.$\!$$\!$}
\setcaptionwidth{5.7in}
\caption{\label{n32768_evolve}Left top: dashed line, initial data $\h -
10^{-4} \h^\dd$; solid lines, solution $h$ at later times. Local minimum is
fixed in space and decreases monotonically after short transient.
Left bottom: dashed line, initial slope; solid lines, $h_x$ at same
times as above ($\|h_x\|_\infty$ increases in time).  Right top:
solution at final resolved time, $t = .0472496406249$. Right bottom:
close--up near $x=0$.}
\vspace{-.5cm}
\end{figure}

\begin{figure}[h] 
  \vspace{-.2cm}
    \begin{center}
      \includegraphics[width=.35\textwidth]{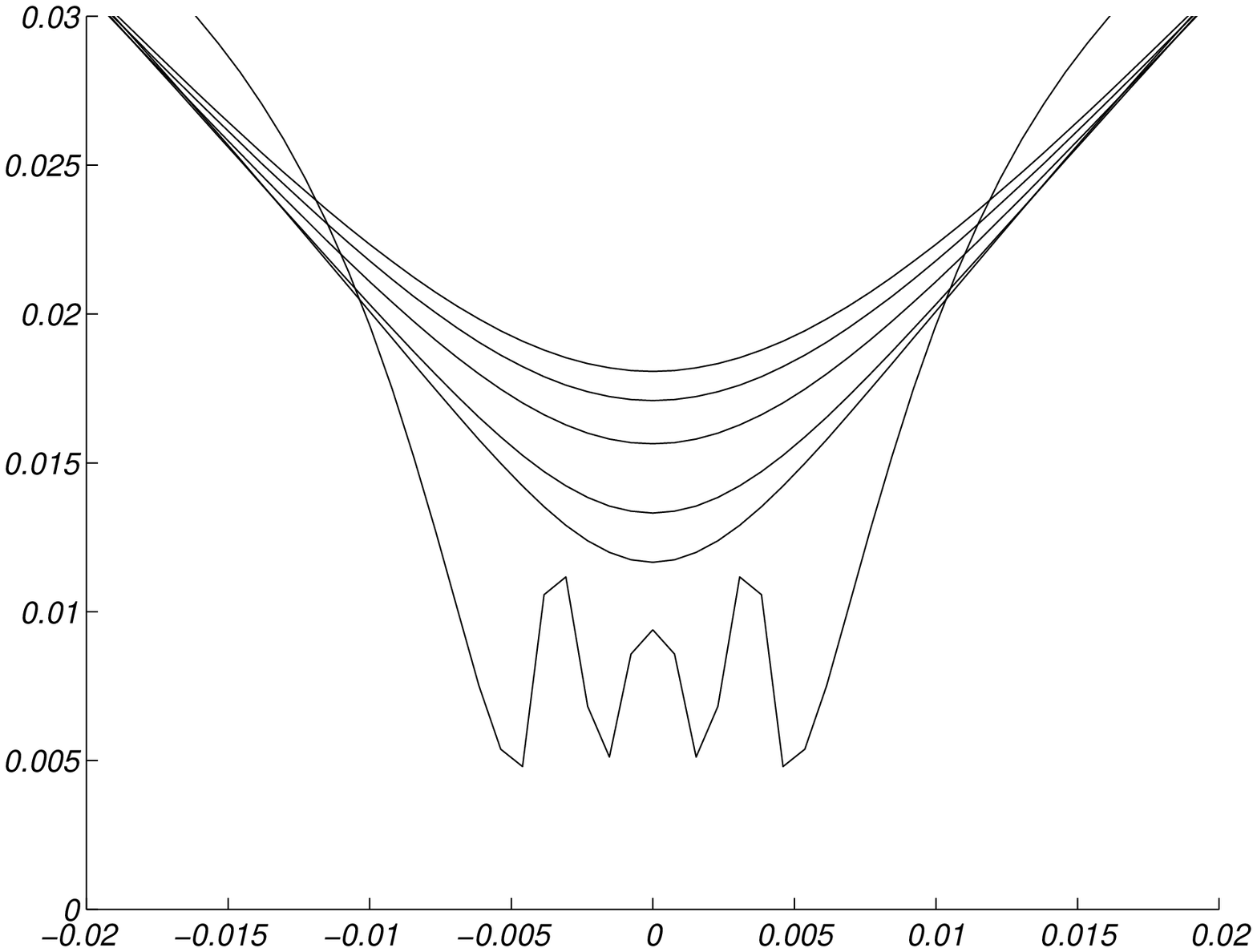}
      \hspace{1.5cm} 
      \includegraphics[width=.35\textwidth]{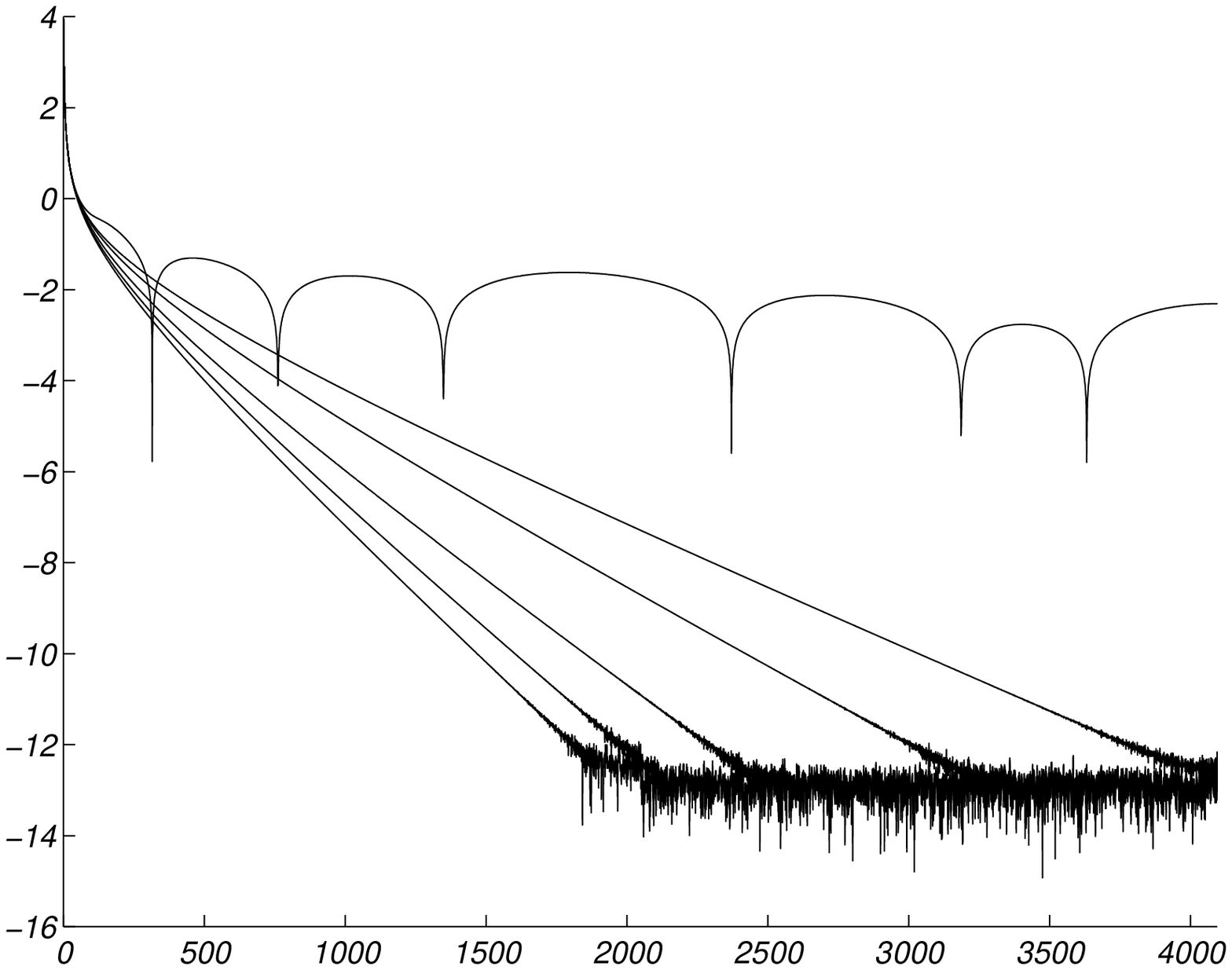}
    \end{center}
  \vspace{-.2cm}	  
\renewcommand{\figurename}{Fig.$\!$$\!$}
\setcaptionwidth{5.5in}
\caption{\label{n8192_ev}Left: evolution near $x=0$ of the solution
with $8192$ meshpoints. Right: the corresponding power spectra. The
power spectrum reaches the Nyquist frequency soon after $t = .04716$,
and the solution loses resolution.  The profile on the left with
multiple local minima is after $t = .04716$. The profiles with a
single local minimum are spectrally resolved.}
\vspace{-.75cm}	  
\end{figure}

In this way, we computed a resolved solution to time 
$t = .0472496406249$, the time when the $32,\!768$ meshpoint 
solution became unresolved.  The top left plot in
Figure~\ref{n32768_evolve} 
presents the evolution of the solution near $x=0$.  As before, the
local extrema are fixed in space, with the solution appearing to touch
down at one point per period. (This was shown previously in \cite[Figure 4c]{WB00}.) We did not design the code to study the
formation of finite--time singularities; the solution has decreased by
only 
a factor of $6.04$ 
due to limitations of our uniform--mesh code (see \S \ref{stop}).
Computing the derivative $h_x$ of the solution, we find that its
maximum and minimum values grow as time passes, as in the bottom left
plot of Figure~\ref{n32768_evolve}.  
These extremum points of $h_x$ move in time, appearing to converge to
$x=0$ as the singular time approaches.  This is consistent with a
solution that touches down with $90^\circ$ contact angles in finite
time. The right plots in Figure~\ref{n32768_evolve} show the final
resolved solution, which suggests $90^\circ$ contact angles are
developing.

The work of Zhang and Lister \cite[\S5]{ZL99} on similarity solutions
suggests that $h(x,t) \sim \B^{1/4} (t_c - t)^{1/5} H(x/(t_c -
t)^{2/5})$ as touchdown approaches; here $t_c$ is the time of
touchdown and $H$ is a particular positive function with $H(\eta) \sim
(0.807) |\eta|^{1/2}$ for large $\eta$.
%
%
Our computations are consistent with the above ansatz. Further, if we
make the ansatz then we can estimate $t_c$, since taking the ratio of
the computed values of $h(0,t)$ at two late times $t_1$ and $t_2$
gives the value of $(t_c - t_1)^{1/5}/(t_c - t_2)^{1/5}$, from which
$t_c$ can be determined.  We find $t_c$ is slightly larger than the
final resolved time.  See \cite{WB99,WB00} for more on the similarity
solutions to (\ref{vanderwaal}).
%

The most plausible future behavior for the solution in
Figure~\ref{n32768_evolve} is that the solution might touch down in
finite time and become a nonnegative weak solution. Then it might
relax, as a weak solution, to a droplet steady state. Alternatively,
since the energy of the original periodic steady state is higher than
the energy of the constant steady state, the solution might perhaps
touch down in finite time and then at some later time become positive
and smooth again, ultimately relaxing to the constant steady state.
This is certainly possible since the solution shown in the right of
Figure \ref{n32768_evolve} has higher energy than the constant steady
state.  However, we did not write our code to study finite--time
singularities or weak solutions, and consequently we cannot
distinguish which (if any) of the above options might be happening.

We close with a graphic demonstration of the kind of spurious effects
that can occur if the solution is {\em not} spectrally resolved.  In
the left plot of Figure~\ref{n8192_ev} we present the $8192$ meshpoint
solution that starts from the same initial data as in
Figure~\ref{n32768_evolve}. The right plot of Figure~\ref{n8192_ev}
shows the corresponding power spectra: as the solution evolves, higher
and higher frequencies are needed to resolve the ever-sharpening local
minimum, until at the final time resolution has been lost. At this
time the solution has multiple oscillations. We found those
oscillations then grew and the solution seemed to touch down in finite
time with two droplets per period, one large and one small. But the
small droplet is a numerical artifact, in view of the absence of small
droplets in the resolved solution shown in Figure~\ref{n32768_evolve}.

\noindent $\bullet$ {\it Odd perturbations.}

We now turn to odd perturbations.  An obvious choice is the
perturbation $\e \h^\d$; it arises naturally from a leftwards
translation of $\h$, since $\h(x+\varepsilon) \approx \h(x)+
\varepsilon\h^\d(x)$. This perturbation shifts the local minimum of
the initial data from $x=0$ to a point slightly to the left of $x =
0$, and it also lowers the minimum value since for small $x < 0$,
\[
\h(x) + \e \h^\d(x) = \h(0) + \e \h^\dd(0) x + O(x^2) < \h(0) .
\]

Given initial data $\h + 10^{-4} \h^\d$, we expect the solution either to
approach zero at some point or to relax to the constant steady state, since $\h$ is energy unstable in the
direction $\pm \h^\d$ by \cite[Theorem~2]{LP3}. The simulation confirms this
expectation, yielding a solution that touches down in a fashion very similar to the solution shown in
Figure~\ref{n32768_evolve}.  One difference is that the location of
the local minimum is no longer fixed in space; it moves as $h_{min}(t)
\to 0$.  The similarities are that the minimum value of the initial
data is very close to zero, relative to the bulk of the initial data,
and that the local minimum of the solution appears to decrease to zero
in finite time while the bulk of the solution is essentially unchanged
from the initial data.


\noindent $\bullet$ {\it Random perturbations.}  

The equation (\ref{vanderwaal}) is
translation invariant so there is nothing to distinguish any
particular point in space.  Up until now, we used either even or odd
perturbations, making $x=0$ appear somewhat distinguished.  To check
the degree to which the behaviors described above depend on the choice
of perturbation, we performed a number of runs with random
perturbations.

The random perturbations are constructed as follows.  The $2048$
meshpoints can resolve $1023$ frequencies, so we chose $\{a_k \}$ and
$\{\phi_k \}$ to be two sets of $1023$ uniformly distributed random
numbers in $[0,1]$ and defined the perturbation
\[
\phi(x) = \sum_{k=1}^{1023} a_k \exp(-.036 k) \cos(k x + 2\pi \phi_k ). 
\]
This perturbation has zero mean, and has random amplitudes and phases
at each wave--number.  The decay rate $0.036$ is chosen so that the
amplitudes at wave numbers $k=900$ and higher are at the level of
round--off error, so that the initial data is spectrally resolved.  We
then divided $\phi$ by its $L^\infty$ norm.

We found that all the solutions resulting from applying such a random
perturbation to $\h$ either relaxed to the constant steady state or
else appeared to touch down in finite time.  An intuitive rule might
be that if the minimum value of the initial data is greater than that
of $\h$ then the solution should relax to the constant steady state
and if the minimum value is less than that of $\h$ then the solution
should touch down in finite time.  In practice, we found that most
solutions respected this intuition, and the solutions that appeared to
touch down in finite time had their gross dynamics as described
earlier (their finer dynamics concern the position of the local
minimum as a function of time).  However, there were exceptions ---
hardly surprising since the evolution equation is a fourth order PDE,
not second order, and the intuitive rule has the flavor of a
comparison principle.

\subsubsection{$q=-3$. Perturbing the constant steady state}
Consider the constant steady state $h \equiv \overline{h}$, for some
Bond number $\B$. If $\B \overline{h}^{q-1}<1$ then $\overline{h}$ is
a strict local minimum of the energy by \cite[Theorem~10]{LP3}, and is
nonlinearly stable. In this case a non-constant positive periodic
steady state $\h$ exists with period $2 \pi$ and mean value
$\overline{h}$ by \cite[Theorem~12]{LP2}, since $0 < \B
\overline{h}^{q-1} (2\pi)^2 < (2\pi)^2$; it is linearly unstable (see
bifurcation diagram~\ref{bifurcs}a).  If $\B \overline{h}^{q-1}>1$
then the constant steady state $\overline{h}$ is a saddle point for
the energy and we expect it to be unstable. In this case there is no
stable positive periodic steady state to which a perturbation could
converge (see bifurcation diagram \ref{bifurcs}a), suggesting that the
solution will touch down in finite or infinite time.

\noindent $\bullet$ {\it Cosine perturbations.}

We take $\B = 1$.  Consider the constant steady state $\overline{h} =
2$, a local minimum of the energy since $\B \overline{h}^{q-1}<1$.
Here the constant steady state is linearly stable and there exists a
linearly unstable positive periodic steady state, by above.  Given
initial data $2 - 10^{-4} \cos x$, the solution appears numerically to
relax to the constant steady state.

\begin{figure}[h]
    \vspace{-.3cm}	
\begin{center}   
    \epsfig{file=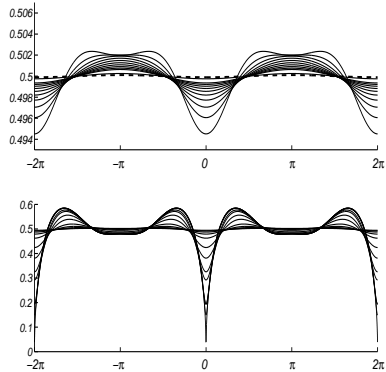, width=2in, height=2in}
\renewcommand{\figurename}{Fig.$\!$$\!$}
\setcaptionwidth{5.7in}
    \vspace{-.2cm}	
\caption{\label{qm3_meanpert}$q=-3,n=3$.  Dashed: initial data $1/2 -
10^{-4}\cos(x)$.  The local minimum is at $x=0$ at all times and
decreases monotonically to zero.}
    \vspace{-.5cm}	
\end{center}   
\end{figure}
On the other hand, the constant steady state $\overline{h} \equiv 1/2$
is linearly unstable and is a saddle point of the energy. Indeed, for
initial data $1/2 - 10^{-4} \cos x$ we find a solution that appears to
touch down in finite time.  Figure \ref{qm3_meanpert} shows this
evolution over two periods.  The top plot shows the short--time
dynamics: the local minimum decreases, while the local maximum
increases for a while.  The top then ``flattens'' and two local maxima
form one to each side of the flat region.  The bottom plot shows the
later-time dynamics: the solution appears to touch down at one point
per period and continues to have two local maxima per period.  The
final profile presented does not look like any known steady state, and
so we expect the solution to continue evolving after touching down.
%

\noindent $\bullet$ {\it Random perturbations.}

Using random perturbations, we verified that the above results are
robust: the steady state $\overline{h} = 2$ is asymptotically stable and
$\overline{h} = 1/2$ is unstable.

\vskip 6pt
Incidentally, we verified that the evolution is exponential in time
near the periodic and constant steady states.  This is consistent with
the nonlinear behavior being dominated by the linear theory when the
solution is sufficiently near a steady state. We found there was a
short transient before the exponential behavior began, suggesting that
the direction $\h^\dd$ is near but not equal to the first
eigendirection.  During the transient time, the solution is locating
this eigendirection.

\subsection{$\mathbf{q=0.5}$}\ 
\label{q.5}

{\it Characteristic features of $q \in (-1,1)$: positive periodic
steady states are linearly unstable. A `Mountain pass' scenario can
occur --- the energy of the non-constant positive periodic steady
state is higher than the energies of the constant steady state and a
zero-contact angle droplet steady state.}  (See bifurcation
diagram~\ref{bifurcs}b and remarks after \cite[Theorem~11]{LP3}.)

\vskip 3pt
We take $n=1, m=.5, q=0.5$ and compute solutions of 
\[
h_t = - (h^1 h_{xxx})_x - \B (h^{.5} h_x)_x.
\]

\subsubsection{$q=0.5$, perturbing the positive periodic steady state}
We rescale a steady state $k_\alpha$ with minimum height
$\alpha = 0.03000$ and period $P(\alpha) = 6.049$.
This yields a Bond number $\B = 0.9628$
and positive periodic steady state $\h$ with period $2\pi$ and area 
$\Ass = 7.980$.
Note $\h$ is linearly unstable, by bifurcation diagram~\ref{bifurcs}b
with $\B \Pss^{3-q} \Ass^{q-1} = 35.59$.

\noindent $\bullet$ {\it Even perturbations.}

As in \S\ref{qm3pp}, we perturb $\h$ with $\pm 10^{-4} \h^\dd$.  

For the initial data $\h + 10^{-4} \h^\dd$, we expect to see evidence
of a heteroclinic connection to the constant steady state
$\overline{\h}$, for the same reasons given for the $q=-3$ case, and
indeed our simulation turned out to be very similar to that shown in
Figure~\ref{qm3_relax}.

\begin{figure}[h]
 \vspace{-.2cm} 
    \begin{center}
       \includegraphics[width=.35\textwidth]{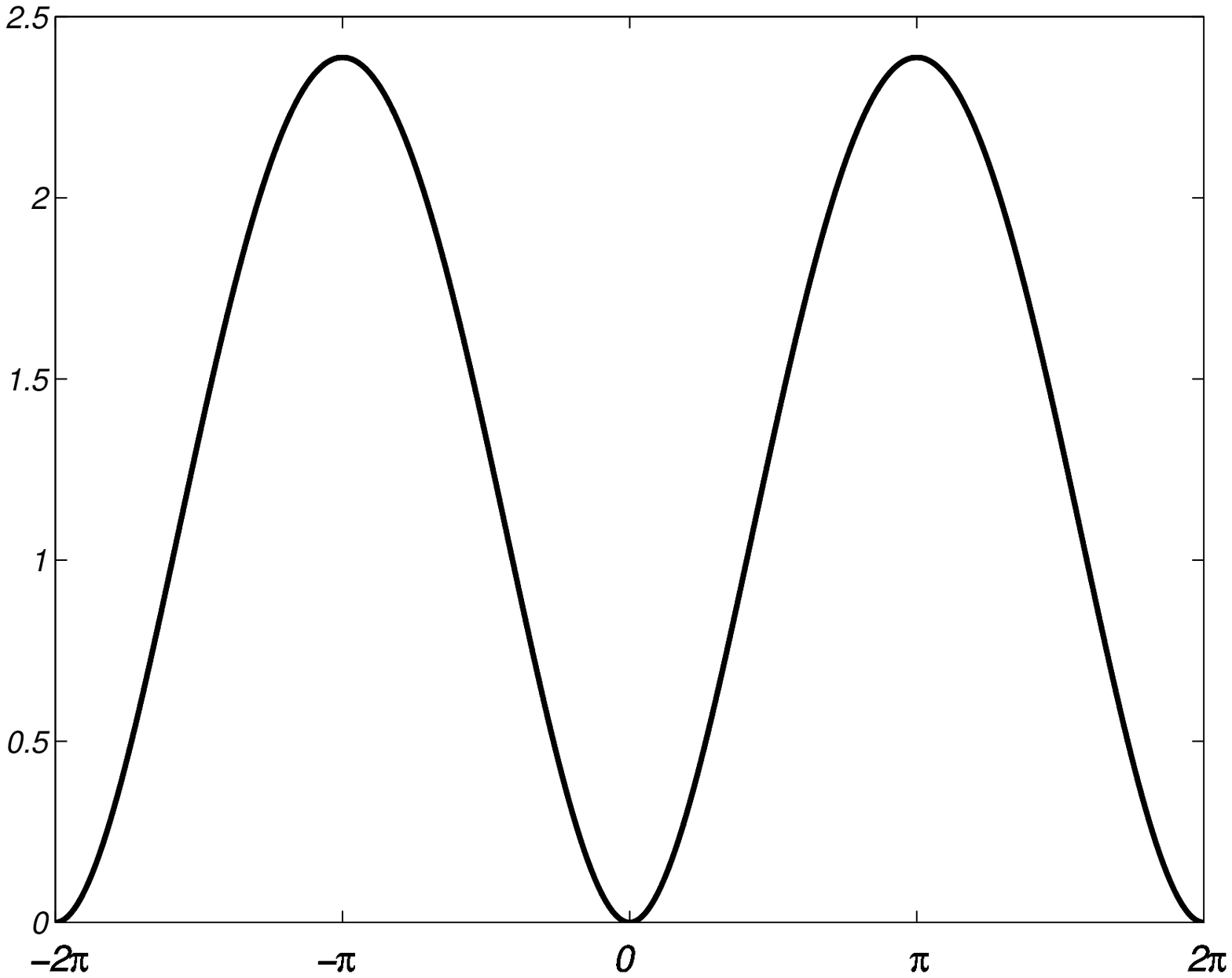} 
       \hspace{1.5cm}
       \includegraphics[width=.35\textwidth]{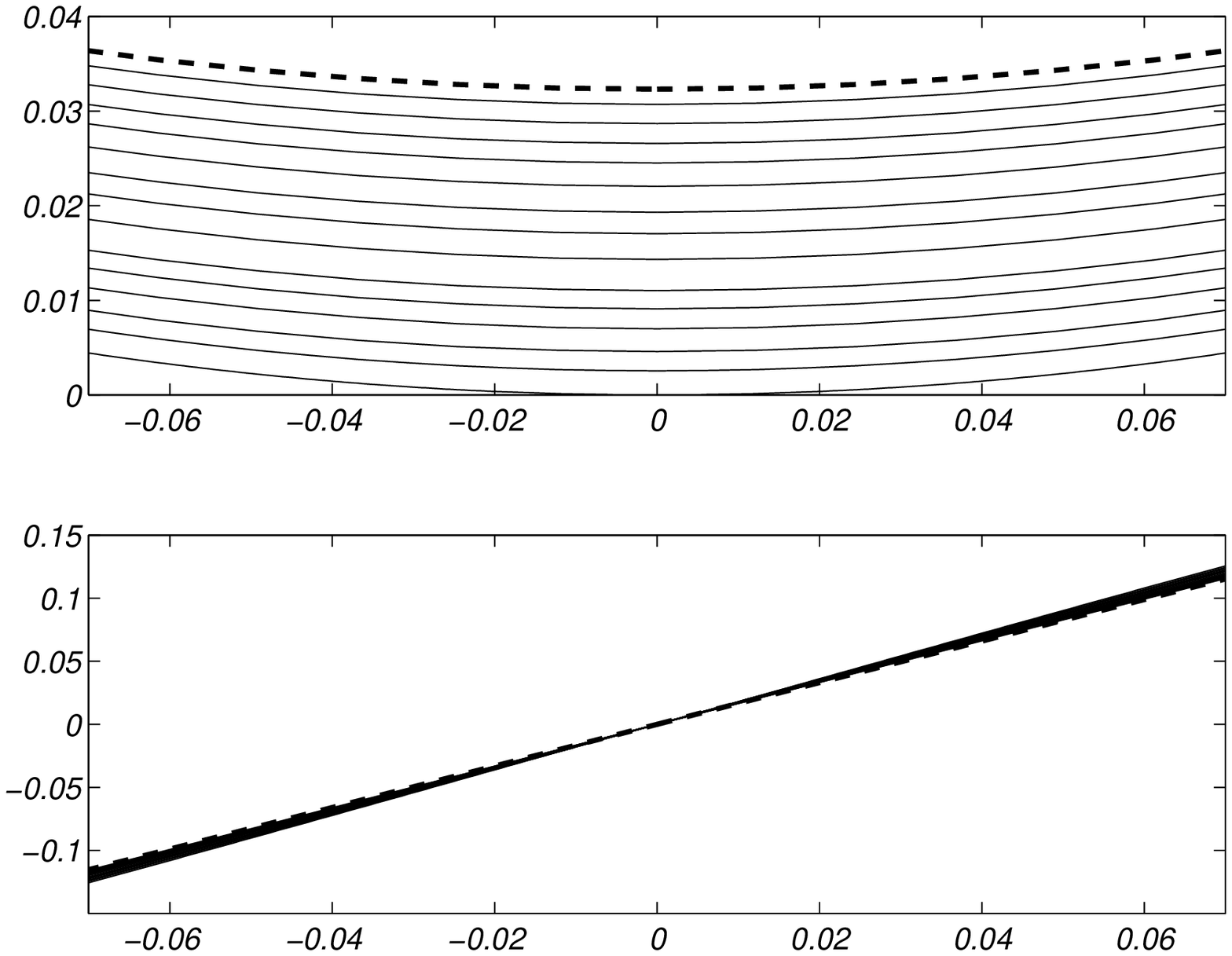}
    \end{center} 
    \vspace{-.25cm} 
\renewcommand{\figurename}{Fig.$\!$$\!$}
\setcaptionwidth{5.5in}
\caption{\label{q.5_pinch} $q=.5$, $n=1$. Left: $h$ at final resolved
time, $t = 7.762540625$, shown over two periods with $8192$ meshpoints per period. Right top: close--up of
$h$ near pinch--off; minimum decreases monotonically after a short
transient; dashed line is initial data. Right bottom: close--up of
$h_x$ at same times.}
\vspace{-.75cm}	  
\end{figure}

Next, for the initial data $\h - 10^{-4} \h^\dd$ we find the solution
appears to touch down in finite time.  Like the $q=-3$ simulation,
$h_{min}(t)$ is located at $x=0$ and, except for a short transient,
decreases monotonically in time.  Figure~\ref{q.5_pinch}a presents $h$
at the final resolved time.
%

Since $-1<q<1$, \cite[Theorem~7]{LP3} tells us there exists a
zero-angle droplet steady state $\hhat$ that has the same area as
$\h$, has length less than $2\pi$, and has lower energy than $\h$. And
indeed Figure~\ref{q.5_pinch}a shows a profile that appears to be
close to a zero-contact-angle solution.  As further evidence of this,
in the top plot of Figure~\ref{q.5_pinch}b we present a close--up of
the evolution near the touch--down point.  In the bottom plot we
present $h_x$ at those times: the slope does not appear to be forming
a jump discontinuity, which it would have to be doing were the
solution converging to a solution with nonzero contact angle.

\noindent $\bullet$ {\it Odd perturbations.}

The odd perturbation $10^{-4} \h^\d$ yields an evolution qualitatively
like that of $q=-3$ above. The solution appears to touch down in
finite time, with the late-time behavior much like that shown in the
right top plot of Figure~\ref{q.5_pinch}.

\noindent $\bullet$ {\it Random perturbations.}

For random perturbations, we observed the same type of dynamics as
seen in the $q=-3$ case.  Some perturbations led to solutions that
relaxed to the constant steady state, and others yielded solutions
that appeared to touch down in finite time with a late--time evolution
like in the right top plot of Figure~\ref{q.5_pinch}.

\vspace{.2cm}

\noindent{\it Mountain pass:} Our simulations above numerically
confirm for $q = .5$ the following `mountain pass' scenario, which is
possible whenever $-1<q<1$.

Choose a positive periodic steady state $\h$ such that $\B
(2\pi)^{3-q} \Ass^{q-1}$ lies between $E_0(q)$ and $4\pi^2$ (see
\cite[Figure~6]{LP3} and refer to the existence result
\cite[Theorem~12]{LP2}). This positive periodic steady state is
linearly unstable and has higher energy than the constant steady state
$h \equiv \overline{\h}$, by \cite[Theorem~7]{LP2} and
\cite[Theorem~6]{LP3}.  It also has higher energy than the zero-angle
droplet steady state $\hhat$ by \cite[Theorem~7]{LP3}. Further, the
constant steady state $\overline{\h}$ is linearly stable and is a local minimum of the
energy, by \cite[Theorem~10]{LP3}.
%
%

Thus the positive periodic steady state $\h$ appears to sit at a
`mountain pass' between the constant steady state $\overline{\h}$
(which lies at the bottom of an energy well) and the droplet state
$\hhat$. If in addition $\B (2\pi)^{3-q} \Ass^{q-1} > L(q)$ (see
\cite[Figure~6]{LP3} for $L(q)\:$) then the constant steady state has
higher energy than the droplet, by \cite[Theorem~11]{LP3}.  But
regardless of that, one would expect that perturbing $\h$ in one
direction would lead to a solution that converges to the constant
steady state while perturbing in the other direction would lead to a
single droplet.  Our numerics for $q=.5$ are all consistent with this
expectation.

\subsubsection{$q=0.5$, perturbing the constant steady state}

\begin{figure}[h]
    \vspace{-.7cm}	
\begin{center}   
    \epsfig{file=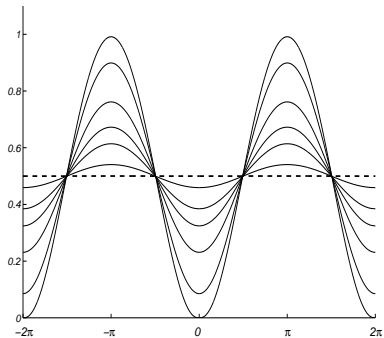, width=2in, height=1.8in}
\renewcommand{\figurename}{Fig.$\!$$\!$}
\setcaptionwidth{5.7in}
    \vspace{-.4cm}	
\caption{\label{q.5_meanpert}$q=.5,n=1$. Dashed: initial data
$1/2-.0001\cos x$. The local minimum and maximum are fixed in space
and, after a short transient, decrease (increase) monotonically.}
\end{center}   
\vspace{-.5cm}
\end{figure}
As in the $q=-3$ case, we take $\B = 1$ and consider perturbing the
constant steady states $h \equiv 2$ ($\B \overline{h}^{q-1}<1$) and $h \equiv 1/2$ ($\B \overline{h}^{q-1}>1$). The first steady state is linearly stable while the second is
linearly unstable.  We perturbed the two states with a range of
perturbations, and found that all perturbations of $h \equiv 2$
relaxed to $h \equiv 2$, while all perturbations of $h \equiv
1/2$ led to apparent finite--time touch--down.  

In Figure~\ref{q.5_meanpert}, for example, we plot the evolution of
the solution with initial data $1/2 - .0001 \cos(x)$.  The extrema are
located at $x=0$ and $2 \pi$ and, after a short transient, increase
(decrease) monotonically.  The evolution shown here is very standard;
we observed this type of behavior more often than the type shown in
Figure~\ref{qm3_meanpert}

\subsection{$\mathbf{q=1}$}\ 
\label{q1}

{\it Characteristic features for $q =1$: all positive periodic steady
states are linearly neutrally stable.} (See bifurcation
diagram~\ref{bifurcs}c and \cite[Lemma~4]{LP3}.)

\vskip 3pt 

Here we consider the $q=1$ case, for which the non-constant positive
periodic steady states are neutrally stable.  We take $n=m=1$ and
compute solutions of
$$
h_t = -(h^1 h_{xxx})_x - \B (h^1 h_x)_x.
$$ 
Goldstein et al.\ \cite[Fig. 3a]{GPS97} found that fairly large
multi-modal perturbations of such steady states relax to steady
states.  They found the solution may relax to a different steady state
than the one of which it was initially a perturbation.

\begin{figure}[h] 
    \vspace{-.2cm}      
    \begin{center}
      \includegraphics[width=.35\textwidth]{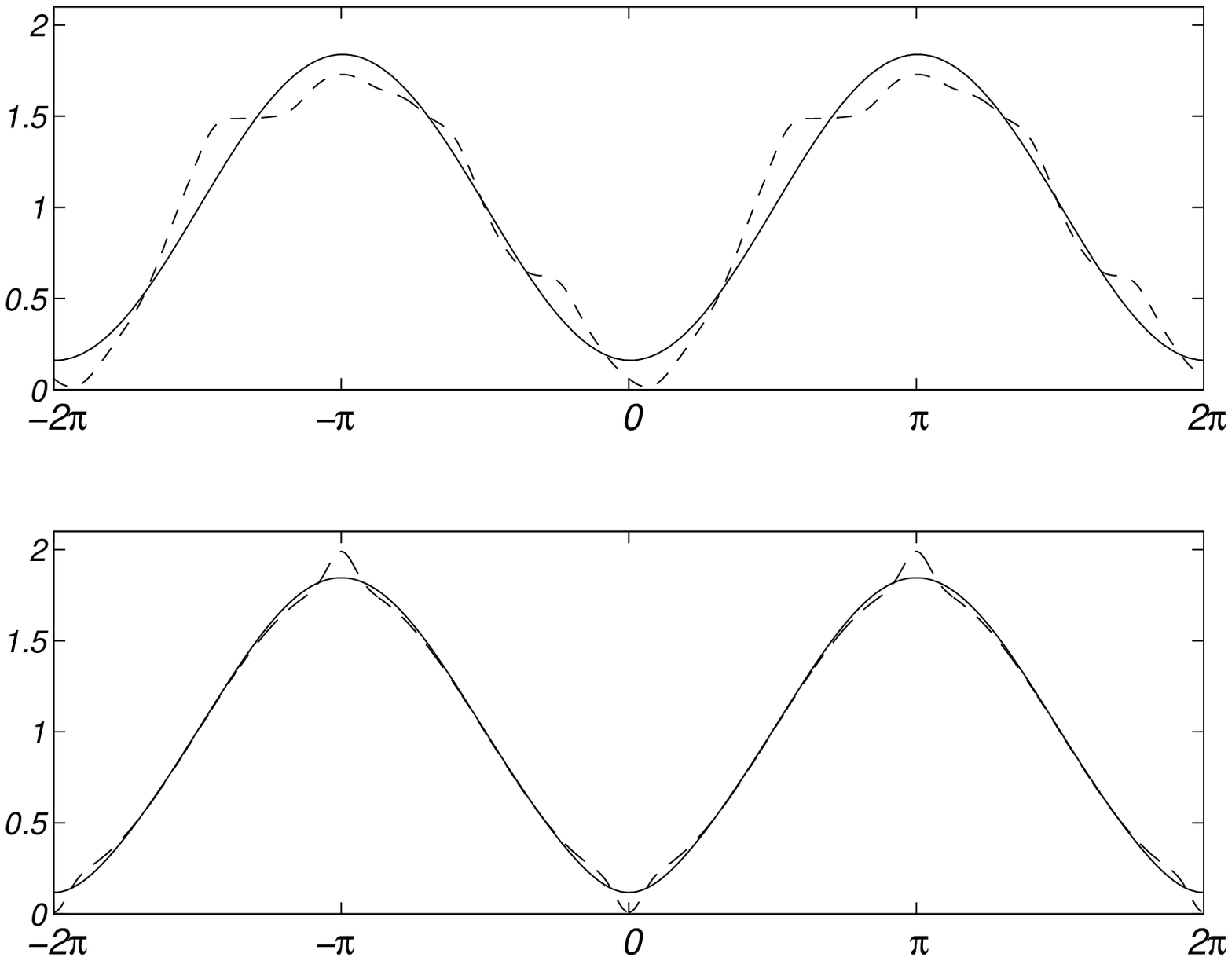}
      \hspace{1.5cm}
      \includegraphics[width=.35\textwidth]{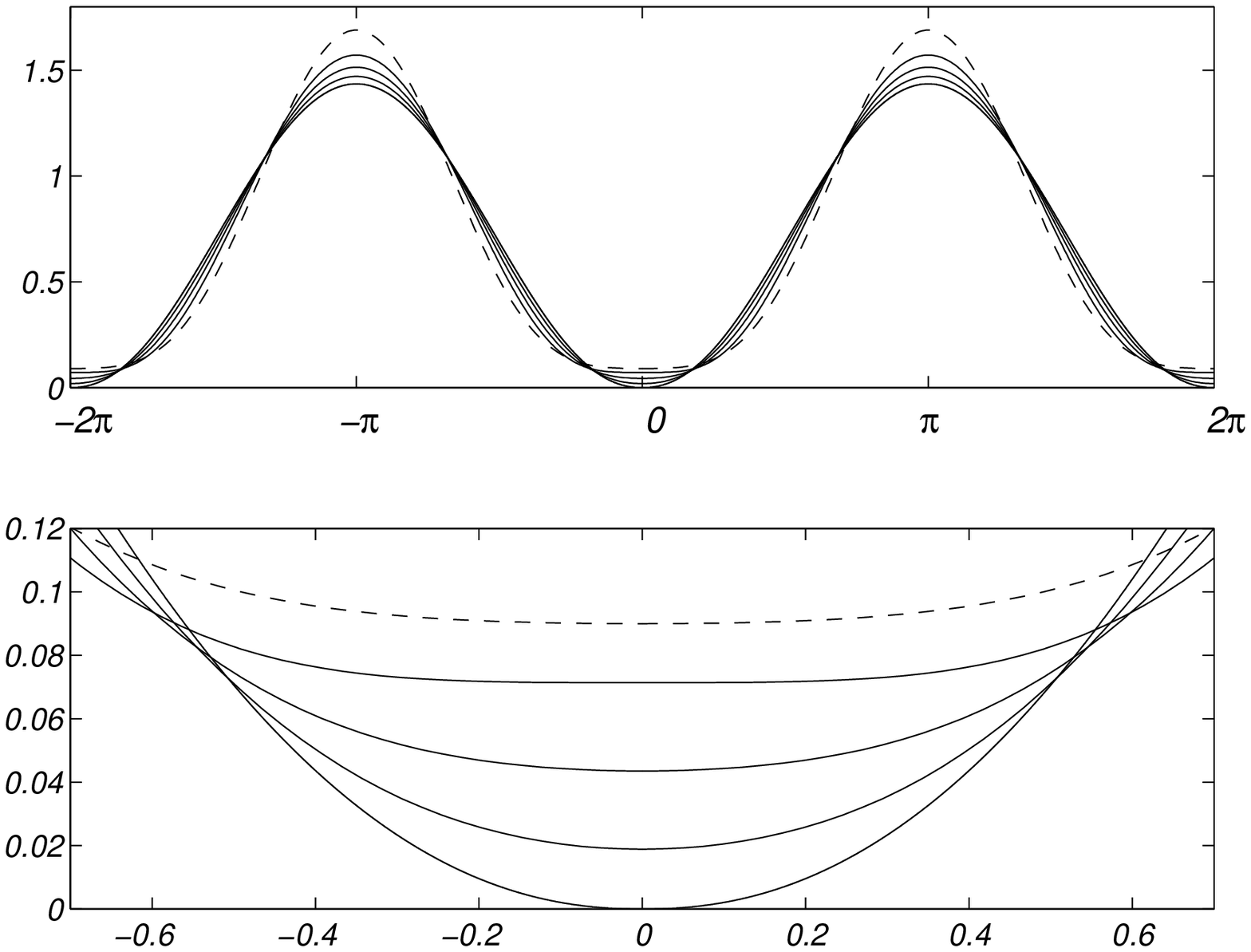}
    \end{center}
    \vspace{-0.2cm}
\renewcommand{\figurename}{Fig.$\!$$\!$}
\setcaptionwidth{5.5in}

\caption{\label{q1_fig}$q = 1$, $n=1$; dashed line: initial data.  Top
left: $h_0(x) = 1 - .8 \cos(x) + .3 v(x)$ where $v$ is a random
perturbation. Bottom left: $h_0(x) = 1 - .8 \cos(x) - .19 \exp(-100
\sin^2(x/2)) + .19 \exp(-100 \sin^2((x-\pi)/2))$.  In both cases,
solution relaxes to a steady state.  Right: $h_0(x) = .7-.8 \cos(x) +
.19 \cos(2 x)$ with minimum value $.09$.  Extrema are fixed in space
and change monotonically (after a short transient), with
$h_{min}(t)$ touching down in finite time.  The final resolved
solution is on $32,\!768$ meshpoints with $h_{min}=1.02 \; 10^{-7}$.}

\vspace{-0.75cm}
\end{figure}

For $\B = 2(1-\cos(\Delta x))/(\Delta x)^2 = 1-O(\Delta x^2)$, one
finds a finite--difference steady state by sampling $a + b \cos(x) + c
\sin(x)$ on a $\Delta x$--uniform mesh (see \S\ref{compute_id}).  In
the left plot of Figure~\ref{q1_fig}, we present two simulations
confirming that each positive periodic steady state is nonlinearly
stable and that a small perturbation of a positive periodic steady
state converges to a (potentially different) positive periodic steady
state.  In the top left plot, we present the evolution from initial
data $1 - .8 \cos(x) + .3 v(x)$ where $v$ is a zero--mean random
perturbation.  In the bottom left plot, we present the evolution from
initial data $ 1 - .8 \cos(x) - .19 \exp(-100 \sin^2(x/2)) + .19
\exp(-100 \sin^2((x-\pi)/2))$.  In both cases, the solution relaxes to
a positive periodic steady state with an amplitude close to $.8$ and a
local minimum close to $x=0$.  We find that the smaller the
perturbation, the closer the long--time limit is to the original
steady state, numerically demonstrating nonlinear stability.  We have
no rule for predicting the amplitude of the long--time limit and,
unless the perturbation is even, we have no way of predicting the
position of the local minimum.

Since these simulations suggest that the non-constant positive
periodic steady states are nonlinearly stable, one might guess that
one cannot find a solution that touches down in finite time.  This is
certainly what we observed for $q=1.5$ and $q=1.768$ in
\S\ref{q1.5}--\ref{q1.768}. And as the bottom left plot of
Figure~\ref{q1_fig} suggests, initial data that has a sharp local
minimum will likely not evolve towards touch--down; the local minimum
will retract in time, as expected for a solution of a surface--tension
driven flow. But initial data that is very flat near its local
minimum, such as $h_0(x) = .7 - .8 \cos(x) + .19 \cos(2 x)$, does
appear to lead to touch--down in finite time, as shown in the top right
plot of Figure~\ref{q1_fig}.  The bottom right plot shows the local
evolution near the touch--down point.  

The $q=1$ case remains mysterious in many ways, because there are
infinitely many $2 \pi$-periodic steady states all having the same
mean value; in the $q \neq 1$ case there are at most two.

{\em Note:} numerical simulations for $q=1$ were earlier presented
with $m=n=1$ in \cite{GPS97}, with $m=n=2$ in \cite[\S8]{Gruen00}, and
with $m=n=3$ in \cite{BPLW}. The latter two articles do not consider
Bond numbers for which periodic steady states might be observed.

\pagebreak

\subsection{$\mathbf{q=1.5}$}\ 
\label{q1.5}

{\it Characteristic features for $q \in (1,1.75]$: positive periodic
steady states are linearly stable.} (See bifurcation
diagrams \ref{bifurcs}d--e.)

\vskip 3pt
We take $n=1$, $m=1.5$, $q=1.5$, and compute solutions of
\[
h_t = - (h^1 h_{xxx})_x - \B (h^{1.5} h_x)_x.
\]

\subsubsection{$q=1.5$. Perturbing the positive periodic
steady state}
\label{q1.5_pert}
We rescale a steady state $k_\alpha$ with minimum height
$\alpha = 0.2145$, period $P_\alpha = 6.453$.
With $D=1$ this yields a Bond number $\B = 1.083$
and positive periodic steady state $\h$ with period $\Pss = 2\pi$ and
area $\Ass = 5.516$, and with minimum at $x=0$.
This positive periodic steady state is linearly stable, unlike for the
other $q$-values considered so far; see bifurcation
diagram~\ref{bifurcs}d with $\B \Pss^{3-q} \Ass^{q-1} = 40.07$.

\noindent $\bullet$ {\it Even perturbations.}  

We perturb $\h$ with
$\pm 10^{-4} \h^\dd$, and expect to see the solution relax back to
$\h$, since $\h$ is linearly stable and since translation in space is
ruled out by the evenness of $\h$ and of our
perturbation. Numerically, we indeed observe relaxation back to $\h$:
the extrema are at $x=0$ and $\pi$ and, after a short transient, they move  monotonically to the maximum and minimum values of $\h$.
%
%

%
\begin{figure}[h]
    \vspace{-0.5cm}	
\begin{center}   
    \epsfig{file=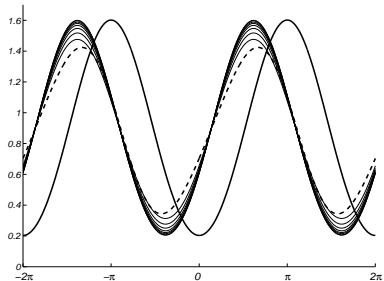, width=2in, height=1.5in}
\renewcommand{\figurename}{Fig.$\!$$\!$} 
\setcaptionwidth{5.7in}
    \vspace{-0.2cm}	
\caption{\label{q1.5_oddpert}$q=1.5$, $n=1$. Dashed line: initial data
$\h + .5 \h^\d$. Isolated solid line: $\h$.  Solution relaxes to a translate
of $\h$; extrema move in time.}
\end{center}   
\end{figure}

\noindent $\bullet$ {\it Odd perturbations.}  

We now perturb $\h$ with
$.5 \h^\d$, which is a large perturbation relative to those discussed
above.  In view of the linear stability and the fact that every
perturbation increases the energy, by \cite[Theorem~5]{LP3}, we expect
to observe relaxation back to a translate of $\h$. (Note the energy and our linear stability results are insensitive to translation.) Convergence back to $\h$ itself (with no translation) seems unlikely since the odd perturbation breaks the
evenness of the initial data.

Figure~\ref{q1.5_oddpert} displays the evolution of the solution with
initial data $\h + .5 \h^\d$; the solution certainly seems to converge
to a translate of $\h$.  We verified this type of behavior for a range
of $q \in (1,1.75]$.

\noindent $\bullet$ {\it Random perturbations.}  

Using random
perturbations, we verified the robustness of the above results: the
periodic steady state $\h$ is asymptotically stable, up to
translation.


\subsubsection{$q=1.5$. Perturbing the positive periodic
steady state with longer perturbations}

\begin{figure}[h]
\vspace{-0.5cm}	
\begin{center}   
    \epsfig{file=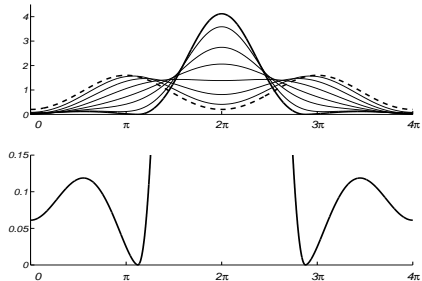, width=2.2in, height=1.5in}
\vspace{-.1cm}
\renewcommand{\figurename}{Fig.$\!$$\!$}
\setcaptionwidth{5.7in}
\vspace{-0.25cm}	
\caption{\label{long_pert} $q=1.5$, $n=1$. Top plot: dashed line is
initial data $\h - .0001 \cos(x/2)$; heavy solid line is final resolved
solution at $t=61.7305$, with $4096$ meshpoints.  Bottom plot: close--up
of final resolved solution.}
\end{center}   
\end{figure}
We have demonstrated above that the positive periodic steady state
$\h$ is asymptotically stable to small perturbations of the same
period, $2\pi$. However, it is linearly {\em unstable} to zero-mean
perturbations of {\em longer} period --- $4\pi, 6\pi$ and so on --- by
\cite[Theorem~1]{LP2}.

For initial data $\h - .0001 \cos(x/2)$, {\it i.e.} a perturbation with period
$4\pi$, the top plot of Figure~\ref{long_pert} presents the evolution
of the solution. The steady state $\h$ has minimum heights at $x = 0,
2\pi$, while the perturbation $- .0001 \cos(x/2)$ is even about $x=2 \pi$ and decreases
the initial value at $x=0$ and increases the value at $x=2\pi$. The solution appears to touch down in finite time,
though it does not do so at $x=0$; also the solution does not appear
to be converging to a single droplet.  The bottom plot of
Figure~\ref{long_pert} shows a close--up of the final resolved
solution.  The smaller droplet is not close to a steady droplet, since
it contains a local minimum within itself --- an impossibility for a
steady droplet.  We expect the solution would continue to evolve as a
nonnegative weak solution, relaxing either to a single steady droplet
or to some (unknown) configuration of steady droplets.

We also considered a number of {\em random} $4\pi$-periodic
perturbations, and always found that the solution appears to touch
down in finite time with one large droplet flanked as in
Figure~\ref{long_pert} by a smaller profile which contains a local
minimum within itself.

\subsubsection{$q=1.5$. Perturbing the constant steady state}
We take $\B = 1$. As always, the constant steady state $\h \equiv
\overline{h}$ is a strict local minimum of the energy if $\B
\overline{h}^{q-1}<1$ (or $\overline{h}<1$), and is a saddle point if
$\B \overline{h}^{q-1}>1$ (or $\overline{h}>1$).

\noindent $\bullet$ {\it Cosine perturbations.}

We consider $\overline{h}= .5$, $1.05$,
and $2$.  We expect $\overline{h} = .5$ will be stable to
all perturbations, and that $\overline{h} = 1.05$ and
$\overline{h} = 2$ will be unstable to some perturbations.  

For $\overline{h} = 1.05$, there exists a linearly stable positive
$2\pi$-periodic steady state with mean value $1.05$, by bifurcation
diagram~\ref{bifurcs}d, using that $E:=\B \overline{h}^{q-1} (2\pi)^2
\in (4 \pi^2,E_0(q))$ when $\B = 1, q = 1.5$ and $\overline{h} =
1.05$, since $E_0(q) = 40.67$.
(See \cite[\S3.1.2]{LP1} for the formula for $E_0(q)$.) A perturbation of $\overline{h} = 1.05$ might converge to this
positive periodic steady state, especially since: $\h$ is linearly
stable; $\h$ should have lower energy than $\overline{h}$, as
discussed after \cite[Theorem~6]{LP3}; and there is no zero-angle
droplet steady state with length less than $2\pi$ and the same area as
$\overline{h}$, by \cite[Theorem~8]{LP3}.  (Note: we cannot exclude
that the solution might converge to a droplet steady state with
nonzero contact angle.)

On the other hand, for $\overline{h} = 2$ there is no positive steady
state with least period $2\pi$ and mean value $2$, since $\B
\overline{h}^{q-1} (2\pi)^2 \notin (4 \pi^2, E_0(q))$.  (Or see bifurcation
diagram~\ref{bifurcs}d.)  For this reason, we expect that
perturbations of $\overline{h} = 2$ should lead to solutions that
converge to a configuration of steady droplets.

\begin{figure}[h]
\vspace{-0.6cm}	
\begin{center}   
    \epsfig{file=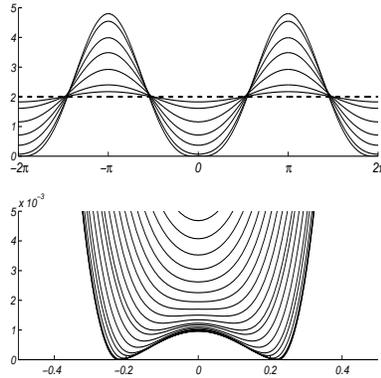, width=2in, height=2in}
\renewcommand{\figurename}{Fig.$\!$$\!$}
\setcaptionwidth{5.7in}
\vspace{-0.3cm}	
\caption{\label{q1.5_meanpert}$q=1.5$, $n=1$. Dashed line: $h_0(x) = 
2 - .0001 \cos(x)$.  The maximum points are fixed; the minimum is initially
at $x=0$ but then splits into two minima which move apart in time as they
decrease.}
\vspace*{-.5cm}
\end{center}   
\end{figure}
We take initial data $\overline{h} - 10^{-4} \cos x$.  We find that if
$\overline{h} = .5$ then the solution relaxes to the constant steady
state and if $\overline{h} = 1.05$ then the solution relaxes to the
positive periodic steady state.  In both cases, the dynamics are very
natural, somewhat like Figure~\ref{qm3_relax} (but run in reverse when
$\overline{h}=1.05$), with the local extremum points fixed in space
and the corresponding extremal values evolving monotonically in time,
after a short transient.  For $\overline{h} = 2$, the solution appears
to be converging to a configuration with two droplets per period (see
Figure \ref{q1.5_meanpert}).  All the solutions shown in this figure
are numerically resolved; the small droplet is real rather than a
numerical artifact, although it may later vanish as the solution
evolves.

\noindent $\bullet$ {\it Random perturbations.}  

Random zero-mean perturbations led to solutions with the same
dynamics: perturbations of $\overline{h}=.5$ relaxed to the mean,
perturbations of $\overline{h} = 1.05$ relaxed to a translate of the
positive periodic steady state $\h$, and perturbations of
$\overline{h} = 2$ appeared to converge to a configuration of two
droplets.

\subsection{$\mathbf{q=1.768}$}\ 
\label{q1.768}

{\it Characteristic features for $q \in (1.75,1.79)$: some positive
periodic steady states are linearly stable, while others
are linearly unstable; and there can be more than one positive
periodic steady state with the same period and area.} (See bifurcation
diagrams \ref{bifurcs}f--h, and \cite[\S5.1]{LP1}.)

\vskip 3pt
We take $n=1$, $m=1.768$, $q = 1.768$, and compute solutions of
$$
h_t = -(h^1 h_{xxx})_x - \B (h^{1.768} h_x)_x.
$$
\subsubsection{$q=1.768$. Perturbing the positive periodic steady
states}

For $q$-values in the interval $(1.75,1.794]$ (approx.) a new
possibility arises: a heteroclinic connection between two
fundamentally different positive periodic steady states.  We
investigate this possibility in what follows.

For Bond number $\B = 1.257$
we consider two distinct positive periodic steady states, $h_{{\rm
ss}1}$ and $h_{{\rm ss}2}$, that have least period $2\pi$, area $A =
4.661$, and have their local minima at $x=0$.  We denote the steady
state that has lower minimum value by $h_{{\rm ss}1}$, and the other
by $h_{{\rm ss}2}$.  Then we expect $h_{{\rm ss}1}$ to be linearly
stable and $h_{{\rm ss}2}$ to be unstable, by \cite[Theorem~9]{LP3}
and its accompanying remarks, with $h_{{\rm ss}1}$ having lower
energy. That is, $h_{{\rm ss}1}$ lies on the stable branch of the
bifurcation diagram~\ref{bifurcs}g and $h_{{\rm ss}2}$ lies on the
unstable branch (since $\B \Pss^{3-q} \Ass^{q-1} = 39.46$.) Note also
that the constant steady state $h := \overline{h_{{\rm ss}1}}$ is
linearly stable since $\B \overline{h}^{q-1}<1$.

\begin{figure}[h] 
\vspace{-.7cm}
 \begin{minipage}{.45\linewidth} 
    \begin{center}      
    \includegraphics[width=.8\textwidth]{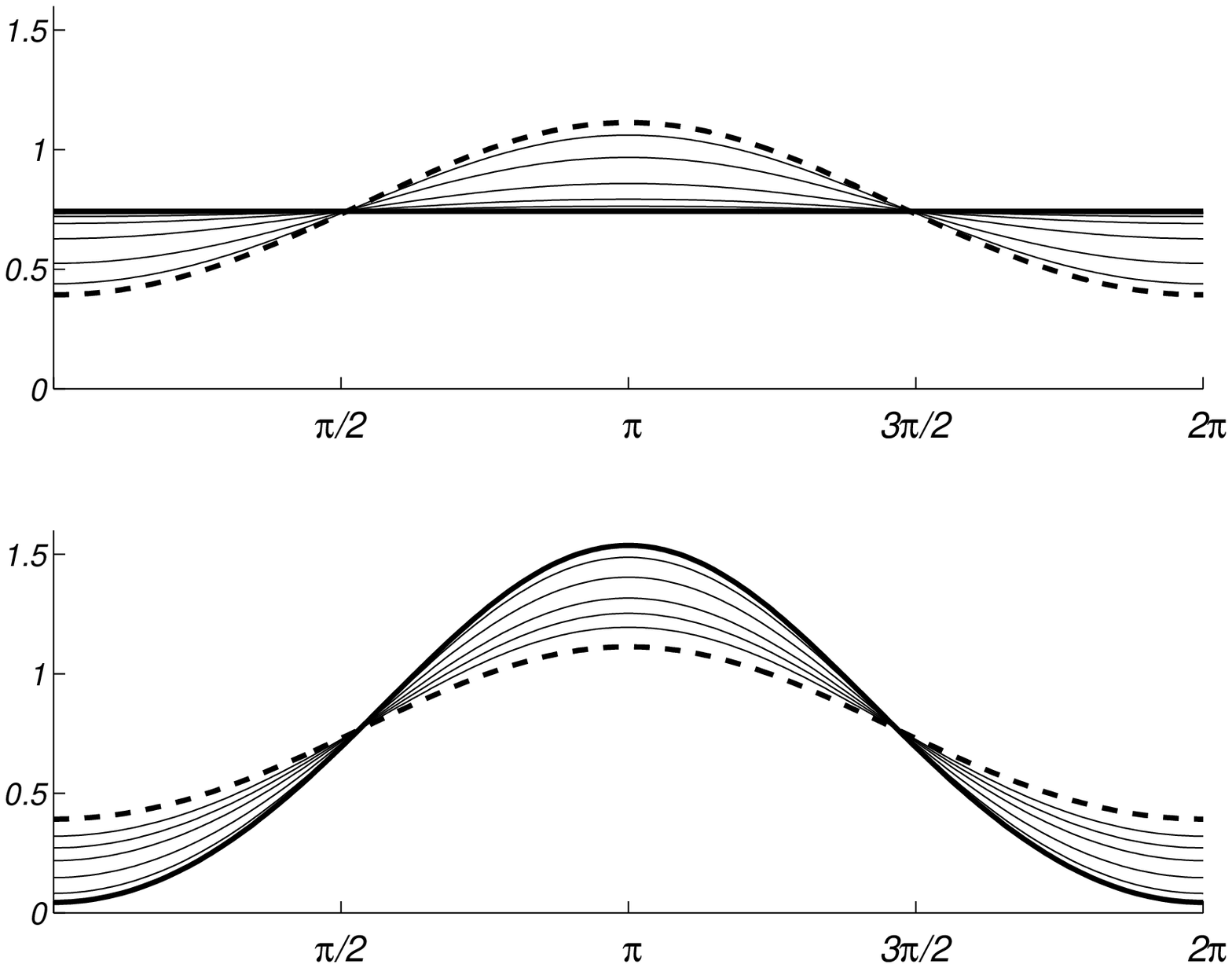}
    \vspace{-.2cm}      
\renewcommand{\figurename}{Fig.$\!$$\!$}
\setcaptionwidth{2.8in}
\caption{\label{q1.768_het}$q = 1.768$ and $n=1$. (a) Dashed:
$h_{{\rm ss}2}+10^{-4}h_{{\rm ss}2}^\dd$; solid: $\overline{h_{{\rm ss}2}}$.  Solution
relaxes to constant.  (b) Dashed: $h_{{\rm ss}2}-10^{-4}h_{{\rm ss}2}^\dd$; solid: $h_{{\rm ss}1}$.  Solution relaxes to $h_{{\rm ss}1}$.}
\end{center} \end{minipage} \begin{minipage}{.45\linewidth}
\begin{center}
    \includegraphics[width=.8\textwidth]{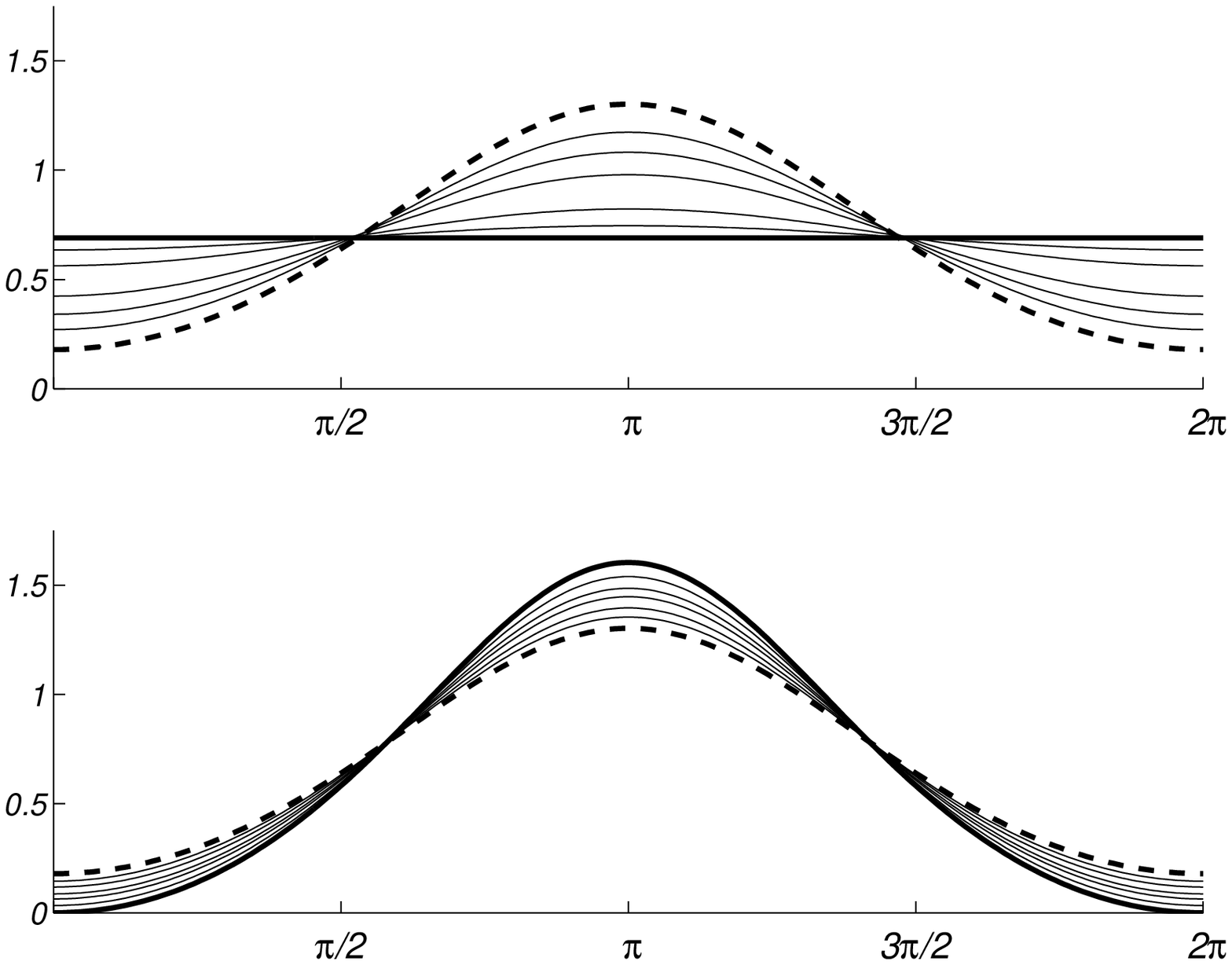}
    \vspace{-.2cm}      
\renewcommand{\figurename}{Fig.$\!$$\!$}
\setcaptionwidth{2.8in}
\caption{\label{q2.5_fig2}$q = 2.5$ and $n=1$. (a) Dashed: $\h +
10^{-4}\h^\dd$; solid: $\overline{\h}$.  Solution relaxes to
$\overline{\h}$.  (b) Dashed: $\h - 10^{-4}\h^\dd$.  Solution 
touches down in finite time.}
\end{center} \end{minipage}
  \vspace{-0.2cm}
\end{figure}

We consider even perturbations ($\e \h^\dd$), odd perturbations ($\e
\h^\d$), and random perturbations.  Our simulations show that $h_{{\rm
ss}1}$ has stability properties like the $q = 1.5$ steady state $\h$ examined in \S\ref{q1.5}. The other steady state $h_{{\rm ss}2}$ is
unstable.  For example, the initial data $h_{{\rm ss}2} + 10^{-4}
h_{{\rm ss}2}^\dd$ yields a solution that converges to the constant
steady state, as shown in Figure~\ref{q1.768_het}(a).  The initial
data $h_{{\rm ss}2} - 10^{-4} h_{{\rm ss}2}^\dd$ yields a solution
converging to $h_{{\rm ss}1}$ as $t \to \infty$, as shown in
Figure~\ref{q1.768_het}(b).  The observed behavior is very robust, and
strongly suggests existence of a heteroclinic connection from the
unstable steady state $h_{{\rm ss}2}$ to the stable one $h_{{\rm
ss}1}$.

Perturbing $h_{{\rm ss}2}$ with {\it random} zero-mean perturbations
yields evidence of heteroclinic connections from $h_{{\rm ss}2}$ to {\em translates} of
$h_{{\rm ss}1}$ and to the constant steady state.

%
%

\vskip 6pt 

We now explain how we found the two steady states $h_{{\rm ss}1}$ and
$h_{{\rm ss}2}$ having the same period and area, since this is not
completely obvious.  First we plot $E(\alpha)$ for $\alpha \in (0,1)$,
as in \cite[Figure~5]{LP3}.  We seek $\alpha_1 < \alpha_2$ such that
$E(\alpha_1) = E(\alpha_2)$ with $E'(\alpha_1) < 0$ and $E'(\alpha_2)
> 0$.  Choosing $\alpha_1 = .05$, it is graphically clear from
\cite[Figure~5]{LP3} that the desired $\alpha_2$ exists. To determine
it, we compute $k_{\alpha_1}$ and its period $P(\alpha_1)= 6.703$ and
area $A(\alpha_1) = 5.660$, and hence $E(\alpha_1)=39.46$. Taking
$\alpha=\alpha_1$ and $D=D_1=1$ and $\Pss = 2\pi$ in the rescaling
relations \cite[eq.\ (24)]{LP3} we evaluate the Bond number as $\B =
1.257$. We then find six values of $\alpha$ such that $E(\alpha) <
E(\alpha_1)$ at the first three values and $E(\alpha) > E(\alpha_1)$
at the last three values, so that $\alpha_2$ is between the third and
fourth values.  We interpolate $E$ at these six values of $\alpha$
with a quintic polynomial and use Newton--Raphson iteration to find
$\alpha_2= 0.5069$ satisfying $E(\alpha_1)=E(\alpha_2)$.  We find
$k_{\alpha_2}$ has period $P(\alpha_2) = 6.388$ and area $A(\alpha_2)
= 6.115$.  Then $P(\alpha_2)$ and the Bond number $\B$ determine the
integration constant $D_2 = 0.8010$ from \cite[eq.\ (24)]{LP3}, using
again $\Pss=2\pi$. Using $D_1$ and $\B$ and \cite[eq.\ (7)]{LP3}, we
rescale $k_{\alpha_1}$ to find $h_{{\rm ss}1}$; by similarly using
$D_2$ we rescale $k_{\alpha_2}$ to find $h_{{\rm ss}2}$. Then we use
$h_{{\rm ss}1}$ and $h_{{\rm ss}2}$ to determine finite difference
steady states, as in \S\ref{compute_id}.
%

\subsection{$\mathbf{q=2.5}$}\ 
\label{q2.5}

{\it Characteristic features of $q \in (1.8,3)$: positive periodic
steady states are linearly unstable. `Mountain pass' scenario can
occur.} (See bifurcation diagram~\ref{bifurcs}i and remarks after
\cite[Theorem~11]{LP3}.)

\vskip 3pt
We take $n=1$, $m=2.5$, $q=2.5$ and compute solutions of 
\[
h_t = - (h^1 h_{xxx})_x - \B (h^{2.5} h_x)_x .
\]

\subsubsection{$q=2.5$. Perturbing the positive periodic steady state}

For Bond number $\B = 1.561$
we consider a positive periodic steady state $\h$ with period $2\pi$
and area $\Ass = 4.335$.
This arises from rescaling $k_\alpha$ with  
minimum height $\alpha =0.2145$ and period $P_\alpha = 6.869$.
Here $\h$ is linearly unstable.

Figure~\ref{q2.5_fig2}(a) presents the solution with initial data $\h
+ 10^{-4}\h^\dd$, which converges to the constant steady state as time
passes.  Figure~\ref{q2.5_fig2}(b) shows the solution with initial
data $\h - 10^{-4}\h^\dd$, which appears to touch down in finite
time. These numerical results are qualitatively the same as for
$q=0.5$, in \S\ref{q.5}. We found this behavior was very robust --- we
considered random perturbations (see \S\ref{qm3pp})
and found for each one that 
the solutions would either converge to the constant steady state or
else touch down in finite time.

\subsubsection{$q=2.5$. Perturbing the constant steady state}

Just as in the $q=-3$ case, our numerics robustly confirm our
predictions: for $\B \overline{h}^{q-1}<1$, small zero-mean
perturbations of $\overline{h}$ yield solutions that relax to the
constant steady state, while for $\B \overline{h}^{q-1}>1$ such
perturbations yield solutions that appear to touch down in finite
time.

\subsection{$\mathbf{q=4}$}\ 
\label{q4}

{\it Characteristic features of $q \in [3,\infty)$: positive periodic
steady states are linearly unstable, and if a positive periodic steady
state and a zero-angle droplet steady state have the same area, then
the period of the former is less than the length of the latter. 
(See \cite[Theorem~7]{LP2} and the proof of \cite[Theorem~7]{LP3}.)}

We take $n=1$, $m=4$, $q=4$ and compute solutions of 
\begin{equation} \label{blow}
h_t = - (h^1 h_{xxx})_x - \B (h^4 h_x)_x .
\end{equation}
This equation is `super-critical' in the sense of Bertozzi and Pugh \cite{BPLW}, since $m>n+2$ ({\it i.e.} $q>3$). According to the conjecture in \cite{BPLW}, then, positive periodic solutions can blow up in finite time ($\|h(\cdot,t)\|_\infty
\to \infty$). Bertozzi and Pugh made the same conjecture
for compactly supported weak solutions on the line, and proved blow-up can occur in finite time when $n=1$ and $m \geq n+2 = 3$ \cite{BPFS}.  Specifically, they proved for such cases 
that if the compactly supported initial data $h_0(x)$ has negative energy
\[
\E(h_0) = \int \left[ \frac{1}{2} h_0^\d(x)^2 
                - \frac{\B}{q(q+1)} h_0(x)^{q+1} \right] dx < 0,
\]
then the compactly supported weak solution blows up in finite time, with its $L^\infty$ and $H^1$ norms both going to infinity.

Here we present computational evidence that smooth {\it periodic}
solutions of (\ref{blow}) can also blow up in finite time (there is no
proof of this).  Further, we find initial data that has positive
energy yet still appears to yield finite--time blow-up, suggesting
that negativity of the energy is not necessary for blow-up, in the
periodic case.

\subsubsection{$q=4$, perturbing the positive periodic steady state}

For Bond number $\B = 4.281$
we consider a positive periodic steady state $\h$ with period $2\pi$
and area $\Ass = 3.213$.  This arises from rescaling $k_\alpha$ with
minimum height $\alpha =0.2145$ and period $P_\alpha = 7.536$.  Here
$\h$ is linearly unstable and the energy $\E(h_0) = .01445$ is positive.
%
%
%
%

\noindent $\bullet$ {\it Even perturbations.}

We considered initial data $\h \pm .001 \h^\dd$.  The initial
data $\h + .001 \h^\dd$ yielded a solution that relaxed to the
constant solution as time evolved.  The local extrema were fixed
in space and, after a short transient, relaxed monotonically to the mean.  

\begin{figure}[h]
    \vspace{-.2cm}	
\begin{center}   
    \epsfig{file=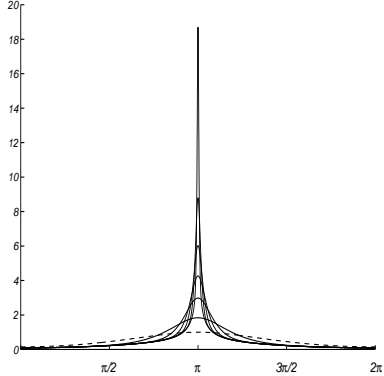, width=2in, height=2in}
\renewcommand{\figurename}{Fig.$\!$$\!$}
\setcaptionwidth{5.7in}
    \vspace{-.5cm}	
\caption{\label{q4_blow}$q=4,n=1$.  Dashed: initial data $\h -.001
\h^\dd$.  In the figure, the maximum increases with time.  The
final profile, at $t=10.0282848199749$, has $4096$ meshpoints.}
\end{center}   
\vspace{-.5cm}	
\end{figure}
The initial data $\h - .001 \h^\dd$ yielded a solution that appears to
blow up in finite time (see Figure \ref{q4_blow}).  The extrema are
fixed in space and, after a short transient, $h_{max}$ increases
monotonically towards infinity (the figure shows $h_{max}$ increasing
by a factor of $18.7$).  After a short transient, $h_{min}$ decreases
monotonically to a positive value as the singular time approaches.

A self-similarity ansatz similar to the case $q=-3$ suggests that
$h(x,t) \sim (t_c - t)^{-1/7} H((x-\pi)/(t_c - t)^{3/14})$ as blowup
approaches at $x=\pi$; here $t_c$ is the time of blowup and $H$ is a
positive function with $H(\eta) \sim C |\eta|^{-2/3}$ for large
$\eta$.  Our computations are consistent with the above ansatz. Again,
if we make the ansatz then we can estimate the blowup time $t_c$,
since taking the ratio of the computed values of $h(\pi,t)$ at two
late times $t_1$ and $t_2$ gives the value of $(t_c - t_1)^{-1/7}/(t_c
- t_2)^{-1/7}$, from which $t_c$ can be determined.  We find that
$t_c$ is slightly larger than the final resolved time.

Self--similar blow-up for super-critical exponents has also been found
for $h_t = -(h^3 h_{xxx})_x -\B (h^6 h_x)_x$ in \cite{Bpc}.

\noindent $\bullet$ {\it Odd perturbations.}

The odd perturbation $10^{-3} \h^\d$ yielded an evolution qualitatively
like that shown in Figure \ref{q4_blow}. The solution appears to blow up in
finite time, with the location of the maximum moving in time.

\noindent $\bullet$ {\it Random perturbations.}

Random perturbations led to the same type of dynamics as seen with
$\pm \h^\dd$ perturbations.


%
%
\subsubsection{$q=4$, perturbing the constant steady state}

Our numerics robustly confirm our predictions: for $\B
\overline{h}^{q-1}<1$, small zero-mean perturbations of $\overline{h}$
yield solutions that relax to the constant steady state, while for $\B
\overline{h}^{q-1}>1$ such perturbations yield solutions that appear
to blow up in finite time.

%
%

%
\section{The effect of changing the mobility exponents, $\mathbf{n}$ and $\mathbf{m}$}
\label{mobility}

In this section, we vary $n$ and $m$ in the equation $h_t
= - (h^n h_{xxx})_x - \B (h^m h_x)_x$. We think of these exponents as
{\em mobility} parameters, since they determine the diffusion
coefficients of the fourth and second order terms in the equation.
When we change $n$ and $m$, we keep $q=m-n+1$ fixed, which means the
steady states of the evolution are unchanged by \S\ref{set_up}. This allows us to ask
three natural questions about the effects of changing the mobility
exponents:
\begin{enumerate}
\item Can a heteroclinic orbit between steady states be {\em broken},
or is it merely perturbed?
\item Can the {\em type} of a singularity be altered ({\it e.g.} from
finite-time to infinite-time)?
\item Can the {\em number} of singularities be altered ({\it e.g.}
from one to two per period)?
\end{enumerate}
The next three subsections address these questions.

Note that while one often thinks of the initial and terminal points of
a heteroclinic connection as being {\em isolated} equilibrium points,
here our equilibria are not isolated, for two reasons. First, the
translates of a steady state are themselves all steady states with the
same area. (Of course this translational freedom disappears when
considering even perturbations of even initial data.) And second, when
we consider configurations of droplet (compactly supported) steady
states that consist of several droplets, the droplets can be
translated, shrunk or expanded, subject only to the requirements that
the total area (volume) be fixed and that the droplets remain
disjoint.

\subsection{Perturbing heteroclinic orbits}
\label{pert_hets}

For $q=2.5$, we computed a $2\pi$-periodic steady state $\h$ with Bond
number $\B = 1.561$ and area $4.335$, by rescaling $k_\alpha$ with
$\alpha=0.2145$ and $P_\alpha = 6.869$.
The steady state is linearly unstable (see bifurcation diagram
\ref{bifurcs}i with $\B \Pss^{3-q} \Ass^{q-1} = 35.32$.)  We take
initial data $\h + .001 \h^\dd$; we expect to find solutions that
relax to the constant steady state, $\overline{\h} \equiv 0.6899$.  We
vary the mobility exponents, by taking $n = 0$, $1$, $2$, and $3$ in
turn, and determining $m$ from $q=2.5=m-n+1$. That is, we compute
solutions for the four evolution equations, all with the same initial
data.

\begin{figure}[h] 
  \vspace{-.225cm}	
    \begin{center}
    \includegraphics[width=.35\textwidth]{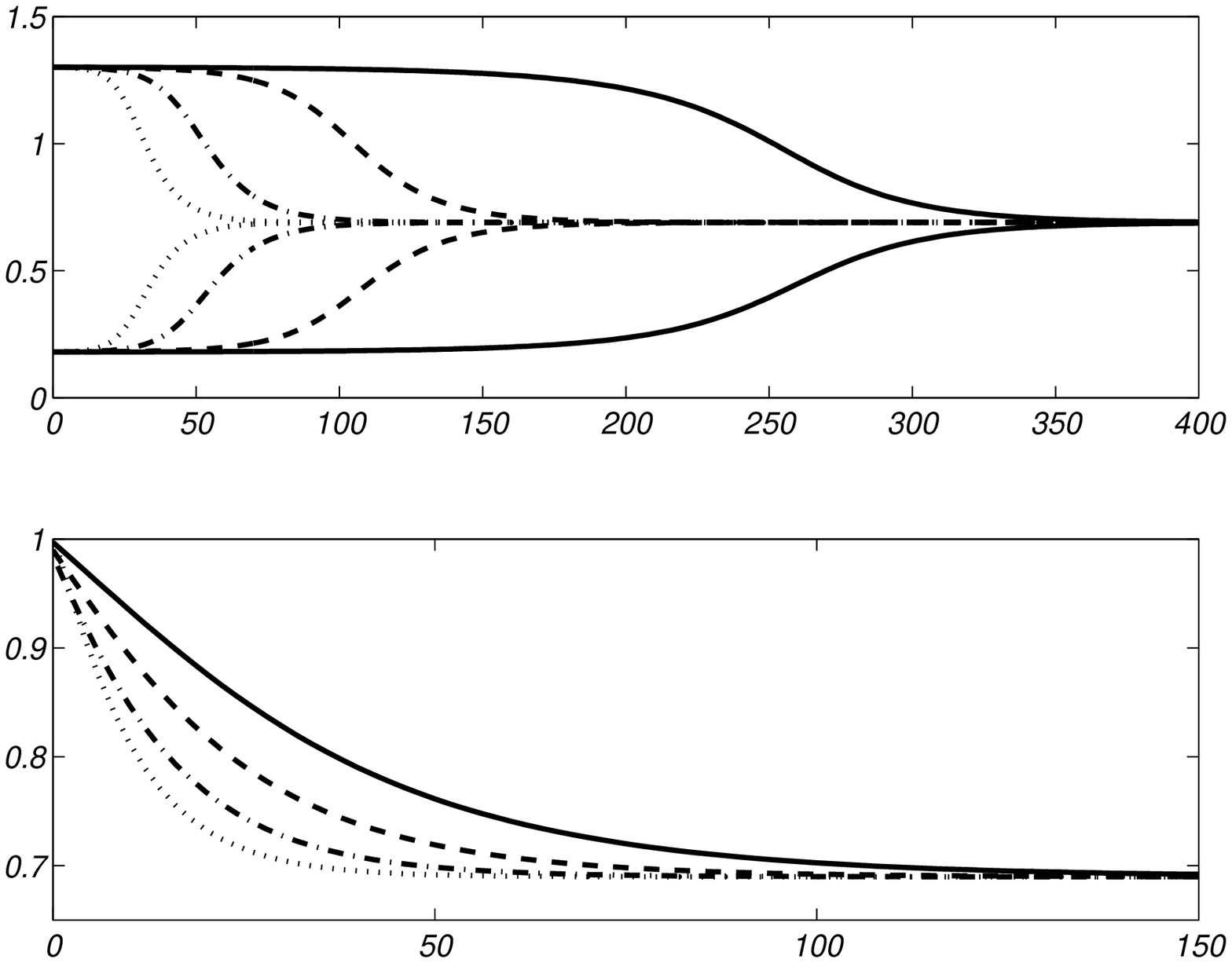}
    \hspace{1.5cm}	
    \includegraphics[width=.35\textwidth]{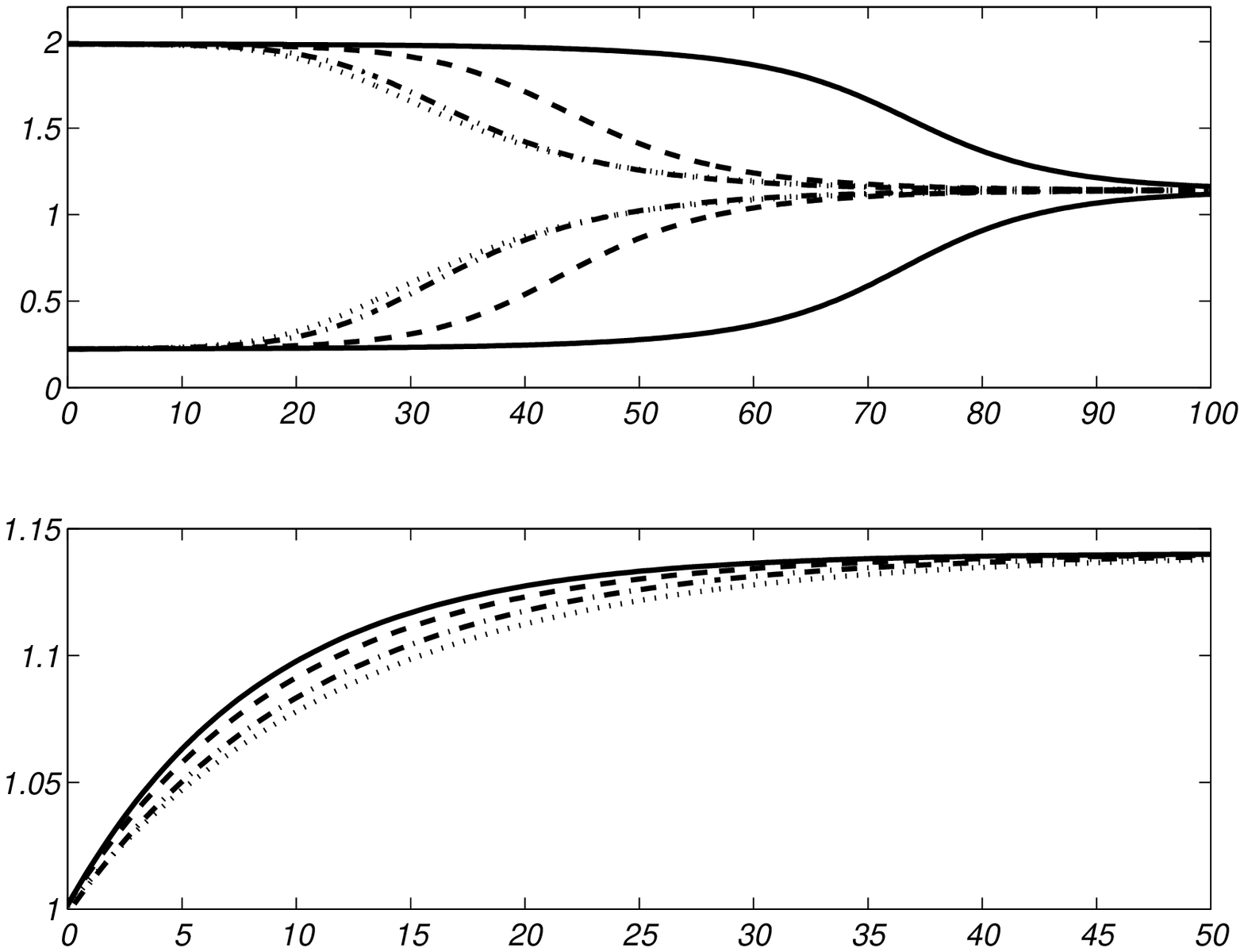}
  \end{center}
\vspace{-0.2cm}
\renewcommand{\figurename}{Fig.$\!$$\!$}
\setcaptionwidth{5.7in}
\caption{\label{time_scales} 
Dotted $n = 0$, dot-dashed $n = 1$, dashed $n = 2$, solid $n = 3$. 
Left: $\!q=2.5$. 
Top left: $h_{max}(t)$ and $h_{min}(t)$ versus $t$.  
Bottom left: $h_{max}(t-t_1)$ versus $t-t_1$, where $h_{max}(t_1) = 1$.  
Right: $\! q=0.5$.  
Top right: $h_{max}(t)$ and $h_{min}(t)$ versus $t$.
Bottom right: $h_{min}(t-\tilde{t}_1)$ versus $t-\tilde{t}_1$ where $h_{min}(\tilde{t}_1) = 1$.}
\vspace{-.8cm}
\end{figure}

We find that all four solutions relax to the constant steady state;
the apparent heteroclinic orbit is robust under this change in
mobility.  In the top left plot of Figure~\ref{time_scales}, we plot
$h_{min}$ and $h_{max}$ versus time for the four solutions.  The
larger the exponent $n$, the longer it takes for the solution to relax
to the constant steady state.

Since $\overline{\h} < 1$, $h_{min}(t) < 1$ for all time and there is
a time $t_1$, dependent on $n$, such that $h_{max}(t) > 1$ for $t <
t_1$ and $h_{max}(t) < 1$ for $t > t_1$.  Near the minimum point,
$h^0>h^1>h^2>h^3$, suggesting that the larger $n$ is, the slower the
diffusion will be (near the minimum).  Similarly, so long as the
maximum is larger than $1$ we have $h^0<h^1<h^2<h^3$ near the maximum,
suggesting that the larger $n$ is, the faster the diffusion will be
(near the maximum).  This conflict of timescales appears to be
mediated through the conservation of mass.  Since the mean of the
solution is conserved, we find that the solution moves as slowly as
its slowest part (which is around the minimum): thus the larger $n$
is, the slower the diffusion.  This is demonstrated in the upper left
plot of Figure \ref{time_scales}: the $n=3$ solution takes longer to
relax than the $n=0$ solution.  Beyond the time $t_1$ there is no
conflict; both the minimum and maximum should relax more slowly as $n$
increases.  We demonstrate this by plotting $h_{max}$ versus $t-t_1$
in the bottom left plot of Figure \ref{time_scales}.

\vskip 6pt
For $q=.5$ we compute a $2\pi$-periodic steady state $\h$ with Bond
number $\B = 0.9817$
and area $A_{ss} = 7.165$.  This arises from rescaling $k_\alpha$ with
minimum height $\alpha=0.2145$ and period $P_\alpha = 6.168$.  The
steady state is linearly unstable (see bifurcation diagram
\ref{bifurcs}b with $\B \Pss^{3-q} \Ass^{q-1} = 36.29$.)
As for $q=2.5$, we take initial data $\h + .001 \h^\dd$ and exponents
$n=0$, $1$, $2$, and $3$.

Again, we find that all four solutions relax to the constant steady
state $\overline{\h} \equiv 1.140$; the apparent heteroclinic orbit is robust
under this change in mobility.  In the top left plot of
Figure~\ref{time_scales}, we plot $h_{min}$ and $h_{max}$ versus time
for the four evolutions.  We see that the larger the exponent $n$, the
longer it takes for the solution to relax to the constant steady
state.

Since $\overline{\h} > 1$, $h_{max}(t) > 1$ for all times and there is
a time $\tilde{t}_1$, dependent on $n$, such that $h_{min}(t) < 1$ for
$t < \tilde{t}_1$ and $h_{min}(t) > 1$ for $t > \tilde{t}_1$.  By the
same logic as before, for $t < \tilde{t}_1$, the time-scales will be
dominated by the dynamics of $h_{min}$.  This is demonstrated in the
upper right plot of Figure \ref{time_scales}: the $n=3$ solution takes
longer to relax than the $n=0$ solution.  Beyond the time
$\tilde{t}_1$, however, both the minimum and maximum should relax more
{\it quickly} as $n$ increases.  We demonstrate this by plotting
$h_{min}$ versus $t-\tilde{t}_1$ in the bottom right plot of Figure
\ref{time_scales}.  There we see the speeds of relaxation reverse, as
expected.

%

\vskip 6pt

For $q=2.5$, we verify as follows that the profile of the solution is
not largely affected by the mobility function $h^n$.  First, we find
the four `half--times': the times at which $h_{min}(t_{1/2}) = .5$.
(The use of $.5$ is essentially arbitrary.)  We find that the four
solutions at their half--times differ by only 0.1\% in the $L^\infty$
norm.
We find analogous results for $q=.5$.  This suggests that using
$h_{min}(t)$ to set a time--scale is an effective way of closely
correlating two points on two different heteroclinic orbits.

\subsection{Changing the type of singularities}
\label{time_pinch}

The choice of mobility coefficients in equation (\ref{evolve}), $h_t =
-(h^n h_{xxx})_x - \B (h^m h_x)_x$, affects whether a positive
solution can become zero somewhere in finite time.  For example, if
$3.5 < n \leq m < n+2$ then it cannot: the solution stays positive for
all time \cite[\S4.2]{BPLW}.  (Note that $m<n+2$ means $q<3$.) On the
other hand, if $m \geq n+2$ then it is possible that the solution
could blow up: $\|h(\cdot,t)\|_{H^1} \to \infty$ in finite time. But
even then we know from the methods of \cite{BPFS} that if $3.5 < n
\leq m$ then the solution remains positive as long as it exists.

Here, we seek the critical exponent $n_c(q)$ such that if $n > n_c(q)$
then positive initial data yield positive solutions for all time,
while if $n < n_c(q)$ it is possible for a positive smooth solution to
touch down in finite time (becoming then a nonnegative weak solution).
From above, if $q<3$ then $n_c(q)\leq 3.5$.

\vskip 6pt 

For $q=1$, Goldstein et al.\ \cite[\S4]{GPS97} presented
simulations with $n=1$ that suggest a finite-time singularity is
possible if $\B > 1$.  Bertozzi and Pugh presented numerical
simulations for $q=1$ and $n=3$ in which the solutions remain positive
for all time and appear to converge to one droplet per period as $t
\to \infty$ \cite{BPLW}. This suggests that $1 \leq n_c(1) \leq 3$.

Here we consider two further $q$-values, $q=2.5$ and $q=.5$.  In each
case we take initial data $\h - .001 \h^\dd$, with the same steady
states $\h$ as in \S\ref{pert_hets}.

\begin{figure}[h] 
  \vspace{-.2cm}	
    \begin{center}
    \includegraphics[width=.35\textwidth]{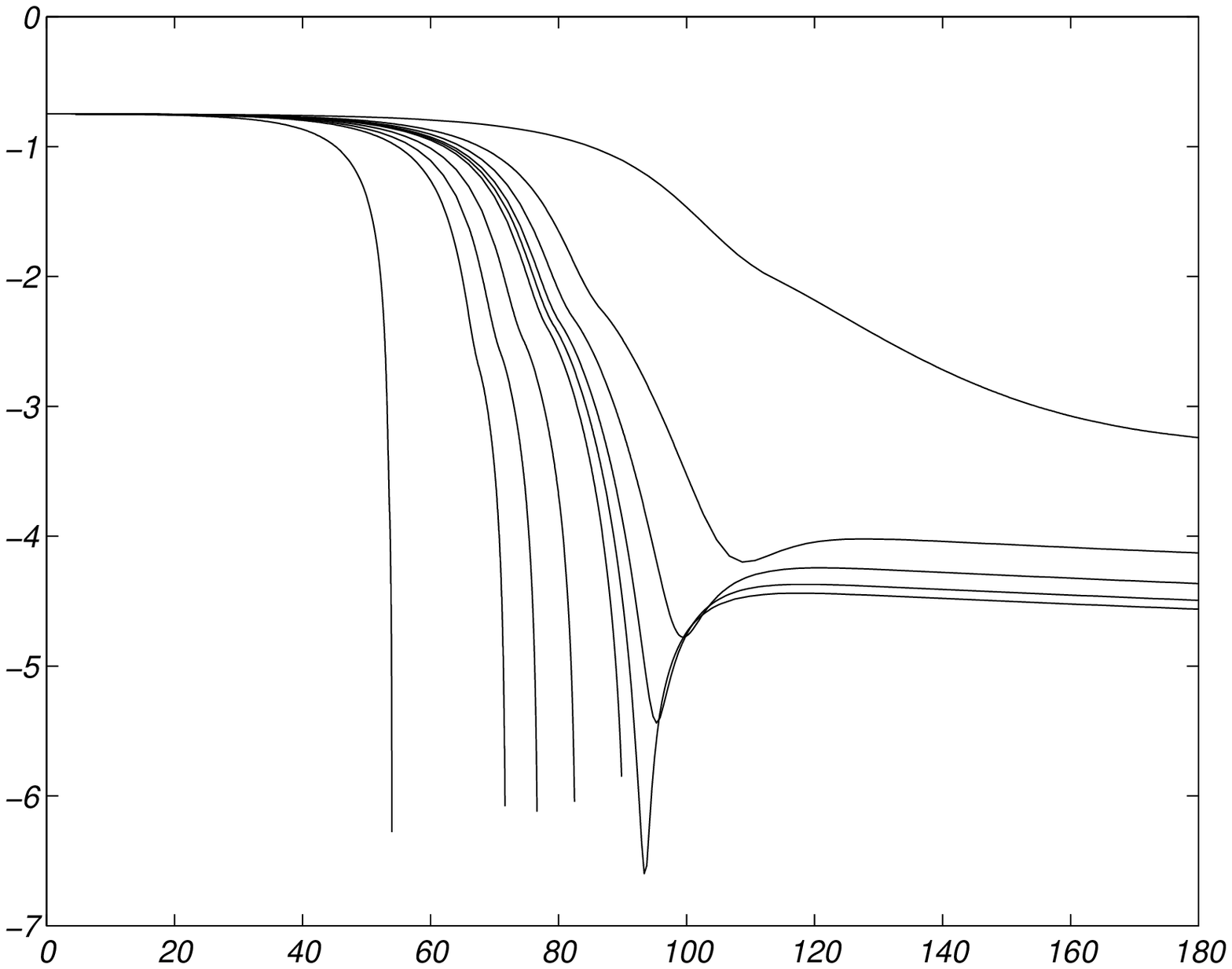}
    \hspace{1.5cm}	
    \includegraphics[width=.35\textwidth]{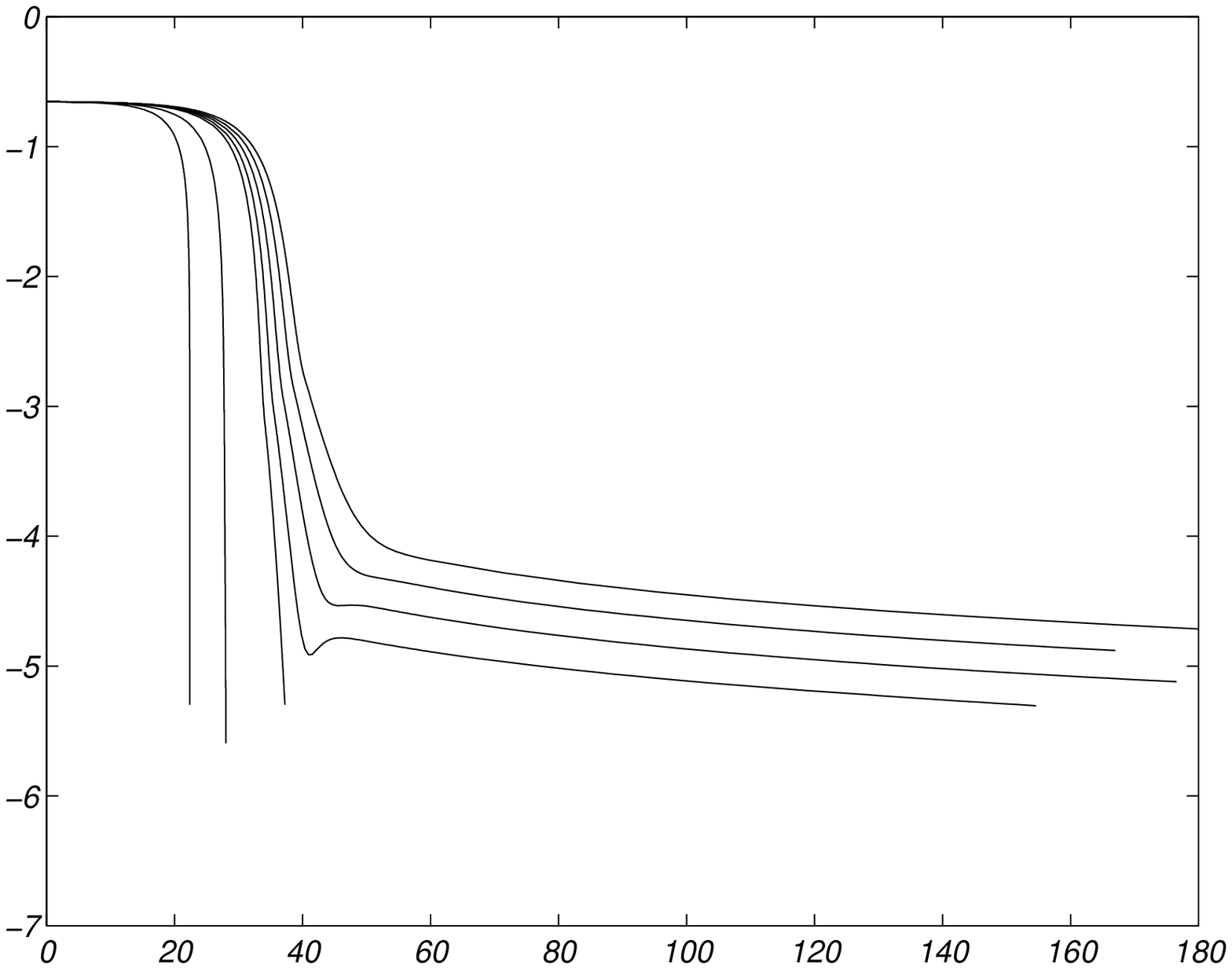}
  \end{center}
\vspace{-0.4cm}
\renewcommand{\figurename}{Fig.$\!$$\!$}
\setcaptionwidth{5.5in}
\caption{\label{find_nc}$\log_{10}(h_{min}(t))$ versus $t$.  Left:
$q=2.5$; $n=1$, $1.25$, $1.5$, $1.55$, $1.6$, $1.65$, $1.6625$,
$1.675$, $1.7$, $1.75,2$.  Right: $q=.5$; $n=1$, $1.5$, $1.8$, $1.85$,
$1.9$, $1.95$, $2$.}
\vspace{-0.6cm}
\end{figure}

We saw for $q=2.5$ and $n=1$, in \S\ref{q2.5}, that solutions appeared
to touch down in finite time, hence $1 \leq n_c(q) \leq 3.5$. To
further approximate the critical $n$-value, we performed simulations
with $n = 1.25$, $1.5$, $1.55$, $1.6$, $1.65$, $1.6625$, $1.675$,
$1.7$, $1.75$, $2$, $2.25$, $2.5$, $2.75$, $3$, and $4$.  Our findings
suggest
\[
1.65 < n_c(2.5) < 1.6625.
\]
In the left plot of Figure \ref{find_nc} we plot $\log_{10}
h_{min}(t)$ versus $t$ for most of these exponents.  If $h_{min}$ is
decreasing at an exponential rate then the graph will be linear at
large times.  If $h_{min}$ is decreasing to zero in finite time with
an algebraic rate then the graph will go to $-\infty$ at some finite
time, dropping down with a vertical slope.  From the plot, if $n=2$
(the rightmost graph) then $h_{min}$ decreases monotonically in time,
eventually decreasing with an exponential rate.  For $n=1.75$, $1.7$,
$1.675$, or $1.6625$, $h_{min}$ decreases, then increases, and then
ultimately decreases with an exponential rate.  The solutions with $n
< 1.6625$ appear to be touching down in finite time.  However the
$n=1.6625$ simulation gives a note of caution; it is possible that the
simulations with $n < 1.6625$ would run until $h_{min}$ became quite
small but would then increase and ultimately decrease exponentially.
Note: all of the simulations were run until the $8192$ meshpoint
simulation lost resolution, except for the $n=1.6625$ simulation which
required $65,\!536$ meshpoints to resolve the solution when $h_{min}$
was at its smallest.  We ran the simulations many decades beyond those
shown to verify the exponential rate of decrease.

\vspace{.2cm}

For $q = .5$ we similarly considered $n=1.25$, $1.4$, $1.5$,
$1.55$, $1.6$, $1.65$, $1.7$, $1.75$, $1.8$, $1.85$, $1.9$, $1.95$,
$2$, $2.25$, $2.5$, $2.75$, $3$, and $4$.  We see very similar
phenomena to the $q=2.5$ case.  The right plot of Figure~\ref{find_nc}
is the analogue of the left plot and suggests that
\[
1.8 \leq n_c(.5) < 1.85.
\]
For $n \geq 1.85$, the solutions stayed positive for
the length of the computation and $h_{min}$ decreased exponentially
in time.

\subsection{Splitting singularities}
\label{one_two_sec}

Our work in \S\ref{pert_hets} suggests that heteroclinic orbits can be
preserved under some changes of the mobility. On the other hand,
qualitative features of a solution can change significantly when the
mobility is changed.  For example, in \S\ref{time_pinch} we
demonstrated that the mobility can affect the regularity of the
solutions; for sufficiently large $n$, solutions are classical for all
time while for smaller $n$, solutions can become nonnegative weak
solutions in finite time.  In this section, we demonstrate another
effect of changing the mobility: a solution that touches down at one
point per period can change into one that touches down at two.

\subsubsection{$\mathbf{q=2.5}$}
\label{one_two_2.5}

We first consider $q = 2.5$ and $n > 0$. Then positive periodic smooth
solutions of (\ref{evolve}) remain bounded in $H^1$ for as long as they exist:
$\|h(\cdot,t)\|_{H^1} \leq M < \infty$ by \cite{BPLW}. We expect that
these solutions will converge to a steady state, as $t \to
\infty$. But the positive periodic steady state is linearly unstable,
and so we expect solutions to converge either to the constant steady
state or to configurations of steady droplets (zero or nonzero contact
angle).

As before, we take initial data $\h - .001 \h^\dd$ with the same $\h$
as in \S\ref{pert_hets}.  In \S\ref{time_pinch}, we found that if $n >
1.6625 \approx n_c(2.5)$ then solutions appear to stay positive for
all time, with $h_{min}$ decreasing to zero exponentially slowly in
time.

\begin{figure}[h] 
  \vspace{-.2cm}	
    \begin{center}
    \includegraphics[width=.35\textwidth]{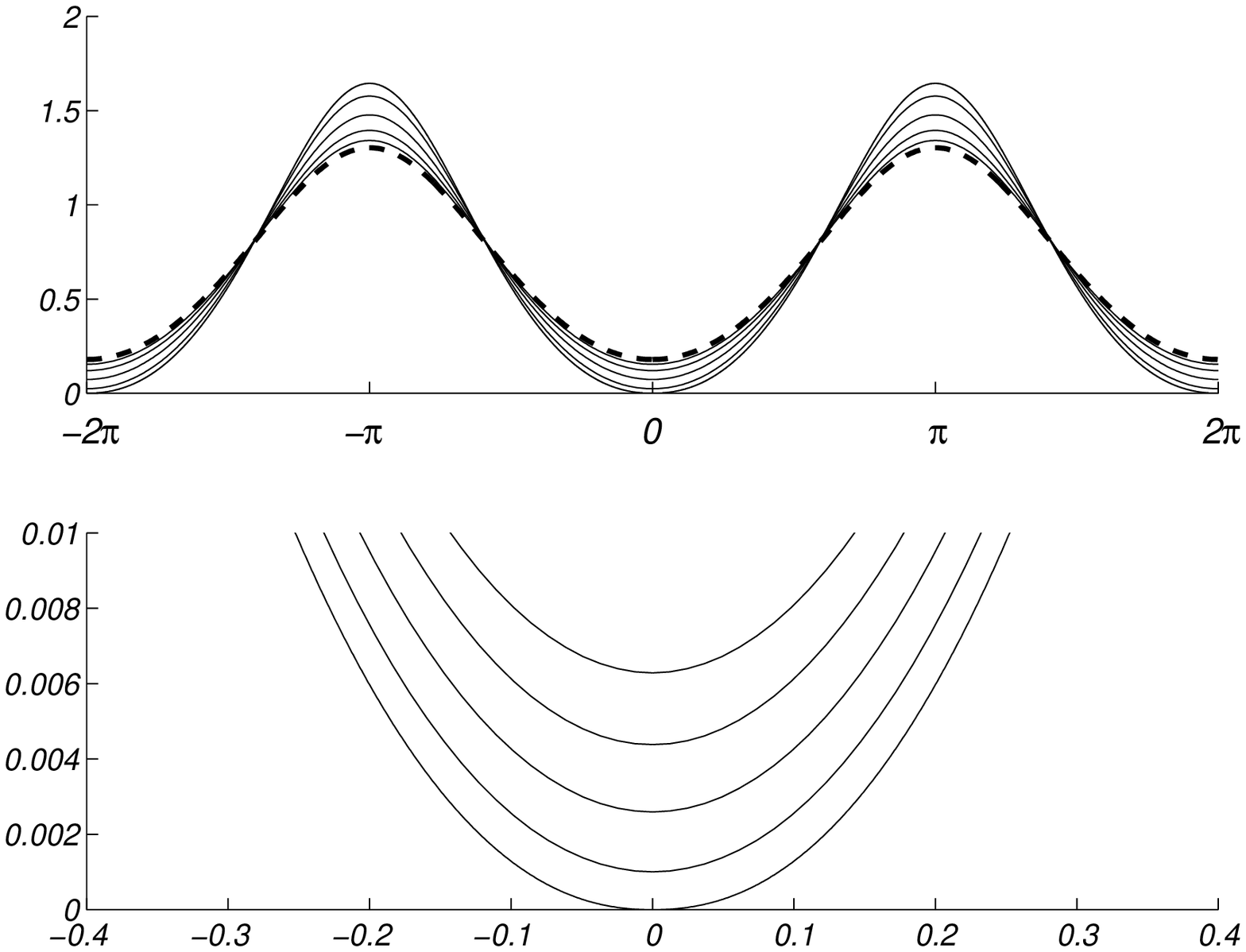}
    \hspace{1.5cm}	
    \includegraphics[width=.35\textwidth]{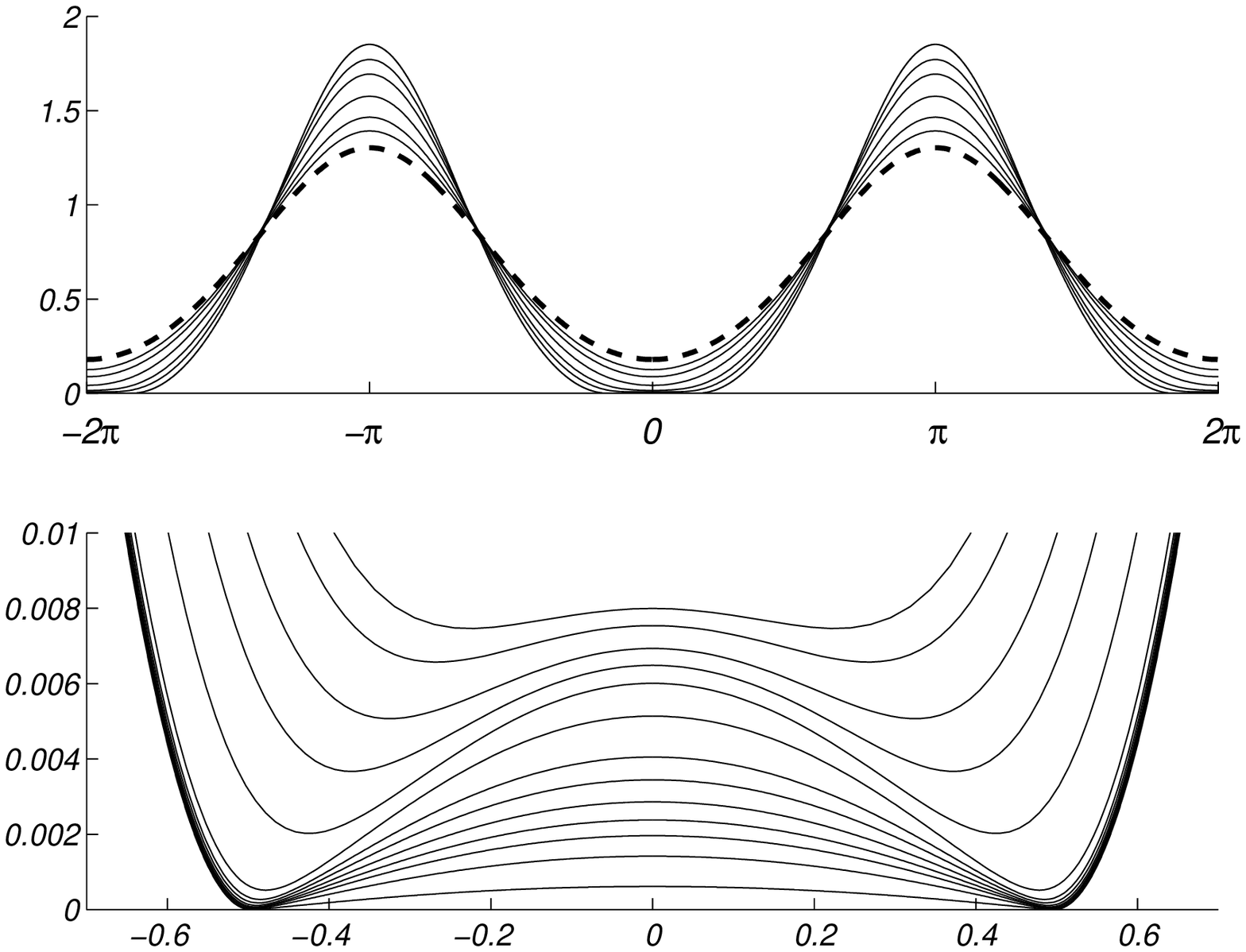}
  \end{center}
\vspace{-0.2cm}
\renewcommand{\figurename}{Fig.$\!$$\!$}
\setcaptionwidth{5.5in}
\caption{\label{mobility_q2.5_n1}$q=2.5$. Left: $n=1$; $h_{min}(t)$
occurs at $x=0$.  At all times there appears to be one droplet per
period.  Right: $n=2$; at late times, the global minima flank $x=0$, 
suggesting a long--time limit of two droplets per period. But in fact
the smaller droplet appears to be vanishing as $h_{min} \to 0$.}
  \vspace{-0.2cm}
\end{figure}

We find for $n = 1$ that the solution touches down in finite time at
one point per period, consistent with a long--time limit of one
droplet per period (see left plot of Figure~\ref{mobility_q2.5_n1}).
For $n=2$ the solution appears to be positive at all times and to
touch down at two points per period in the long--time limit (see right
plot of Figure \ref{mobility_q2.5_n1}).  This suggests a long--time
limit of two steady droplets per period. But it is {\em impossible} to
contain two zero contact angle steady droplets in an interval of
length $2 \pi$, as we argue shortly.  In fact, we find that the small
`proto-droplet' is actually draining, with its maximum decreasing to
zero like $t^{-.4024}$. A similar phenomenon was observed by
Constantin et al.\ \cite[\S{III},IV]{First} with $n=1$ and $\B=0$ (and
with different boundary conditions), although their proto-droplet
seemed to decay like $1/t$.  In our case, we find that the draining
rate depends on $n$.

Our simulations suggest that a second critical exponent,
$\tilde{n}_c(q)$, governs the number of touch--downs per period, at
least for the even perturbations we are using.  If $n <
\tilde{n}_c(q)$ then there appears to be one touch--down per period,
occurring at $x=0$.  If $n > \tilde{n}_c(q)$ then there appear to be
two touch--downs per period, with the position of the local minimum
moving in time and with the solution being non-symmetric about the
local minimum.  The singularity splits as $n$ increases through
$\tilde{n}_c(q)$.  Goldstein et al.\ \cite[\S4C]{GPS97} observed
something similar for $h_t = - (h h_{xxx})_x - \B (h h_x)_x$
($q=n=m=1$).  Specifically, they found a single symmetric singularity
that splits into a pair of asymmetric singularities as $\B$ increases
from $1$ past $\B \approx 1.35$.

\begin{figure}[h] 
  \vspace{-.2cm}	
    \begin{center}
    \includegraphics[width=.35\textwidth]{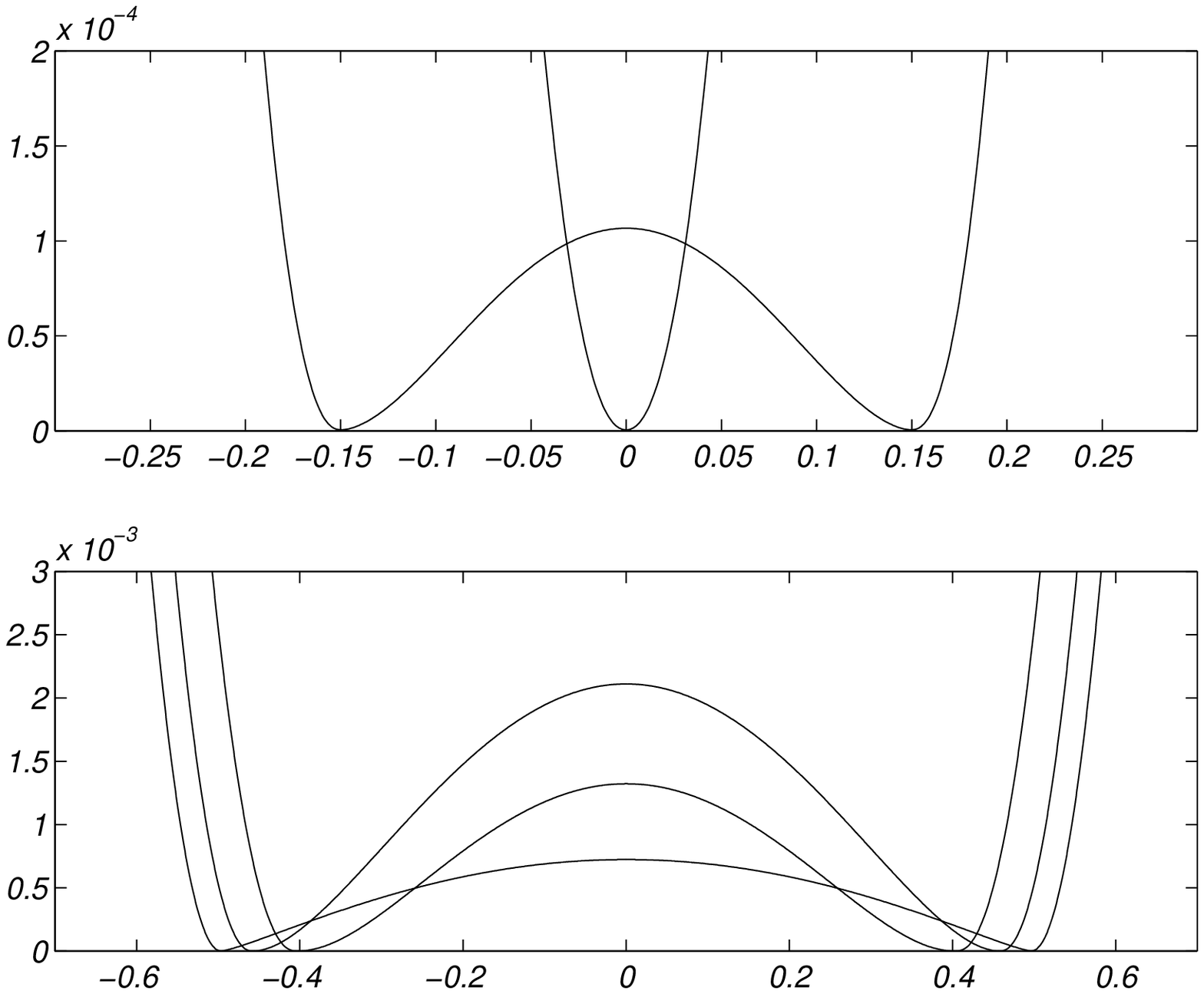}
    \hspace{1.5cm}	
    \includegraphics[width=.35\textwidth]{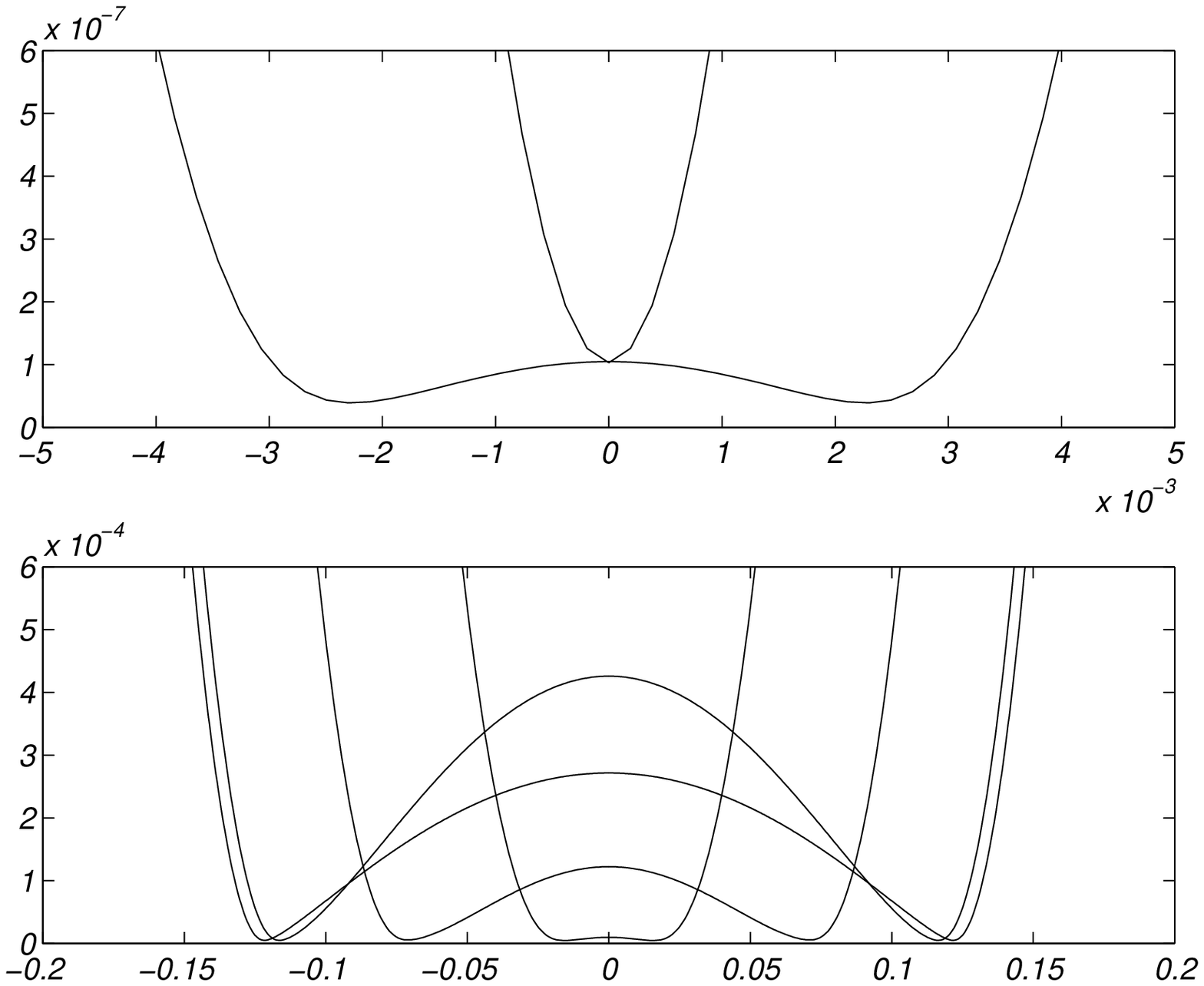}
  \end{center}
\vspace{-0.2cm}
\renewcommand{\figurename}{Fig.$\!$$\!$}
\setcaptionwidth{5.5in}
\caption{\label{one_two}Late--time profiles of $h$.  Left:
$q=2.5$. Top: $n=1$ (single minimum) and $1.25$.  Bottom: $n = 1.5$,
$1.6$, $1.7$; width of proto-droplet increases with $n$.  Right:
$q=.5$. Top: $n=1$ (single minimum) and $1.25$.  Bottom: $n = 1.4$,
$1.6$, $1.8$, $2$; width of proto-droplet increases with $n$.}
\end{figure}

We find that 
\[
1 < \tilde{n}_c(2.5) < 1.25 .
\]
In the left plot of
Figure~\ref{one_two} we plot the late--time profiles for a range of
$n$.  All five profiles shown are the final resolved solution with
$8192$ meshpoints.  In the top plot, we plot the profiles from $n=1$
and $1.25$.  The $n=1$ profile has only one local minimum, while the
$n=1.25$ profile has two.  In the bottom plot, we plot the profiles
from $n=1.5$, $1.6$ and $n=1.7$.  Each profile has two local minima,
with the distance between the minima increasing with $n$.  

The $n=3$ and $4$ evolutions appear to be very similar to the $n=2$
evolution.  We did not do any simulations beyond $n=4$ since the
larger the value of $n$, the longer the simulation had to run before
we could observe anything tangible --- the mobility coefficient $h^n$
is very small near the local minimum, where $h$ is small.  We could
not compensate by taking large time--steps, because $h^n$ could be
quite large near the local maximum, since $h_{max}(t) > 1$.  (Large
differences in time--scales are difficult to handle numerically.)

\vskip 3pt 

We now prove our earlier claim that for $q=2.5$ and $\B = 1.561$,
there cannot be two disjoint zero contact angle steady droplets in an
interval of length $2 \pi$, if the total area of the droplets is
$\Ass=4.335$ (the area of the initial data in
\S\ref{pert_hets}). There does exist a {\em single} zero angle droplet
steady state with that area and with length less than $2\pi$ (the
length is $\Pss = \Ass^{-3} \left( E_0(2.5)/\B \right)^2 = 5.287 < 2
\pi$), by \cite[Theorem~7]{LP3}.
Hence there {\em is} a zero contact angle droplet to which the solution on the right of Figure~\ref{mobility_q2.5_n1} could relax.  But if there were a pair of zero contact angle steady state
droplets, with areas $A_1 = \lambda \Ass$ and $A_2 = (1-\lambda) \Ass$, then
the combined length of the two droplets would be
\[
P_1 + P_2 = \left( \lambda^{-3} + (1-\lambda)^{-3} \right) \Ass^{-3}
                             \left( E_0(2.5)/\B \right)^2. 
\]
The righthand side is a convex function of $\lambda$, achieving its
minimum value $84.5905$ at $\lambda = 1/2$.  This minimum value is
greater than $2 \pi$, and so one cannot fit the two droplet steady
states in an interval of length $2 \pi$.  Thus the simulations
described earlier cannot be converging to a pair of steady zero
contact angle droplets.  Indeed, our computations show the
proto-droplet is slowly draining.

\subsubsection{$\mathbf{q=.5}$}
\label{one_two_.5}

For $q=.5$ we took initial data $\h - .001 \h^\dd$ (for $\h$ as in
\S\ref{pert_hets}) and observed phenomena very similar to those in the
$q=2.5$ case.  We find there is an exponent $\tilde{n}_c(.5)$ such
that if $n < \tilde{n}_c(.5)$ then there is only one touch--down per
period and if $n > \tilde{n}_c(.5)$ then there can be two. Again, our
simulations suggest
\[
1 < n_c(.5)<1.25 ;
\]
see the plots to the right of Figure~\ref{one_two}.  In the right top
plot, we present the final resolved solutions for $n=1$ and $1.25$
with $32,\!768$ meshpoints.  The local minima for $n=1.25$ are
much closer to each other than in the $q=2.5$ plot.  In fact, at
$8192$ meshpoints it appeared that there would be only one local
minimum; as more time passed it split into two.  In the right bottom
plot of Figure~\ref{one_two} we plot solutions with $n=1.4$, $1.6$,
$1.8$, and $2$.  The $n=1.4$, $1.6$, and $1.8$ solutions are with
$8192$ meshpoints, and the $n=2$ solution is with $16,\!384$
meshpoints.  The profiles are all resolved and were chosen to have
comparable $h_{min}$.  The distance between the local minima in the
plot increase monotonically with $n$.

The evolutions for $n=1$ and $2$ are very similar to those shown in
the plots of Figure~\ref{mobility_q2.5_n1}.  Specifically, for $n=2$
the long--time limit appears to be one droplet.  This is interesting
since for $q=.5$ (unlike for $q=2.5$), it {\em is} possible to have
two disjoint zero contact angle steady solutions in an interval of
length $2 \pi$, as we now show.  We have a Bond number $\B = 0.9817$
and area $\Ass = 7.165$ (from \S\ref{pert_hets}).
A single zero contact angle steady state would satisfy \cite[\S
3.1]{LP1}
\[
\Pss = \Ass^{1/5} \left( E_0(.5)/\B \right)^{2/5},
\]
where $E_0(.5) = 32.86$.
We find $\Pss = 6.038 < 2 \pi$.
Thus a single zero contact angle steady state is a potential
long--time limit.  For two zero contact angle steady states, we find
$$
P_1 + P_2 = \left(\lambda^{1/5}+(1-\lambda)^{1/5}\right)
                 A^{1/5}
                 \left( E_0(.5)/\B \right)^{2/5}.
$$
This is a concave function with its maximum at $\lambda = 1/2$.  We
find that $P_1+P_2 < 2\pi$ if $\lambda \leq 10^{-7}$ (approx.) and so
one {\em can} have two zero contact angle droplets --- but one of them
must be fairly small.  For example, if $\lambda = 10^{-7}$ then the
length of the smaller droplet is $P_1 = 0.2404$.
In the bottom right plot of Figure~\ref{one_two} the profile for $n=2$
(the longest proto-droplet shown) has length $.2431$, which is close
to $.2404$; thus a two--droplet long--time limit is at least a
possibility.  However, like for $q=2.5$, the proto-droplet appears to
shrink in time, with its maximum value decreasing to zero like
$t^{-.4056}$.  (As for $q=2.5$, the draining rate also depends on
$n$.)  We did not succeed in finding a $q$, $n$, and initial data that
yield a solution with a multi-droplet configuration as its longtime
limit.  But we believe it should be possible to find this somehow.

\section{Numerical methods} 
\label{numerics}

The numerical simulations are done using a finite-difference evolution
code.  Throughout this section, the diffusion coefficients are
represented as functions $f$ and $g$; we use power law coefficients
$f(y)=y^n$ and $g(y)=\B y^m$ in our simulations for this paper.  The
exponent $\ell$ represents the $\ell^{\rm th}$ time-step and $h^\ell$
is the numerical approximation of the solution at time $t_\ell = \ell
\dt.$

\subsection{The evolution code}
\label{the_code}
We use an adaptive time--stepping scheme based on a Crank-Nicolson
scheme:
\begin{eqnarray} \nonumber 
\frac{h^{\ell+1}-h^\ell}{\dt} 
               &=& - \half \left(f(h^{\ell+1/2})h^{\ell+1}_{xxx}\right)_x
               - \half \left(f(h^{\ell+1/2})h^{\ell}_{xxx}\right)_x \\
&& \qquad \qquad  - \half \left(g(h^{\ell+1/2})h^{\ell+1}_{x}\right)_x
\nonumber              - \half \left(g(h^{\ell+1/2})h^{\ell}_{x}\right)_x.
\end{eqnarray}
Here the diffusion coefficients are evaluated at $h^{\ell+1/2}$, which
we find by linearly extrapolating the solutions $h^{\ell-1}$ and
$h^\ell$ to time $t_\ell + \dt/2$.  The scheme is $O(\dt^2)$ on time
intervals where $\dt$ is fixed and is $O(\dt)$ at the time when the
timestep is changed.  We explain the adaptive timestepping in
\S\ref{acc}.

Finding $h^{\ell+1}$ reduces to solving a linear problem, which we
write in a residual formulation:
\begin{equation} \label{resid_crank_nic}
\LL z = - \dt \left[ f(h^{\ell+1/2}) (h^{\ell}_{xxx} +
                                      r(h^{\ell+1/2}) h^{\ell}_x) \right]_x
\end{equation}
where $z = h^{\ell+1}-h^{\ell}$ and $r(y)=g(y)/f(y)$, and the linear
operator $\LL$ is defined by
\[
\LL z := z 
+ \dt \; \half \left[ f(h^{\ell+1/2})(z_{xxx} + r(h^{\ell+1/2}) z_x ) \right]_x.
\]
The scheme uses $h^{\ell-1}$ and $h^\ell$ to compute $h^{\ell+1}$; for
the first step we take $h^{-1} = h^0$.  \\

We now perform a linear stability analysis of the scheme about the
constant steady state $h \equiv \hbar$, for power law coefficients
$f(y) = y^n$ and $g(y) = \B y^m$.  The linearized equation is $h_t = -
\hbar^{\: n} h_{xxxx} - \B \hbar^{\: m} h_{xx}$.  Initial data $h^0(x) = \e
\cos(k x)$ yields $h^{\ell +1} = \sigma(k)^{\ell+1} \; \e \cos(k x)$
where
\[
\sigma(k) = 
  \frac{1 - \frac{\dt}{2} \hbar^n k^2 (k^2 - \B \hbar^{m-n})}
 {1 + \frac{\dt}{2} \hbar^n k^2 (k^2 - \B \hbar^{m-n})} =:
  \frac{1-\mu(k)}{1+\mu(k)}.
\]
Hence there is a band of unstable modes: if $0 < k^2 < \B \hbar^{m-n}$
({\it i.e.} $\mu(k) < 0$) then $|\sigma(k)| > 1$ and the initial data
$h_0(x) = \e \cos(k x)$ yields a solution that grows exponentially in
time.  We consider a numerical scheme linearly stable if perturbations
{\it outside} this band are not amplified:
\[
k^2 \geq k_c^2 := \B \hbar^{m-n} \Longrightarrow
| \sigma(k) | \leq 1.
\]
It follows immediately from the form of the growth factor that the
Crank-Nicolson scheme is linearly stable, as expected.  Such a linear
stability analysis provides a useful guide, though it is directly
relevant only for small perturbations of flat steady states. \\

For the spatial discretization, the key issue is to implement the
scheme in a way that preserves steady states.  By (\ref{ss1}), a steady state $\h$ satisfies ${\h}_{xxx} + r(\h) {\h}_x = 0$.  Such an `analytic steady state' will
not generally be a `finite-difference steady state', although $\h$
will be $O(\dx^2)$-close to the finite-difference steady state
$\tilde{h}$ satisfying the following discretization:
\begin{equation} \label{fd_ss}
\tilde{h}_{i+2}
- 3 \tilde{h}_{i+1}
+ 3 \tilde{h}_{i}
- \tilde{h}_{i-1}
+ \dx^2 \frac{ r(\tilde{h}_{i+1})+ r(\tilde{h}_{i})}{2}
\left( \tilde{h}_{i+1} -  \tilde{h}_{i} \right) = 0 ,
\quad i=1\dots N.
\end{equation}
The meshpoints are $x_1 = \dx, x_2 = 2 \dx ,\cdots ,x_{N} = X$, where
$X$ is the length of the interval, and we denote the function values
at the meshpoints with subscripts: $\tilde{h}_1, \tilde{h}_2, \dots
\tilde{h}_{N}$.  The function is periodic: $\tilde{h}_0 =
\tilde{h}_{N}$.

To implement the residual formulation (\ref{resid_crank_nic}), we
apply the $O(\dx^2)$ approximation
\begin{eqnarray}
\left[ \left(f(h) (z_{xxx} + r(h) z_x) \right)_x \right]_i
& \simeq &
\frac{f_{i+}}{\dx}\left(\frac{z_{i+2}-3z_{i+1}+3z_i-z_{i-1}}{\dx^3}
           + r_{i+}\frac{z_{i+1}-z_i}{\dx} \right) \nonumber \\
&& \qquad \qquad 
- \frac{f_{i-}}{\dx}\left(\frac{z_{i+1}-3z_{i}+3z_{i-1}-z_{i-2}}{\dx^3}
                                 + r_{i-}\frac{z_{i}-z_{i-1}}{\dx} \right) , \label{resid}
\end{eqnarray}
where the subscripts $i+$ and $i-$ denote the right-average and
left-average, for example:
\[
r_{i+} := \frac{r(h^{\ell+1/2}_{i+1}) 
       + r(h^{\ell+1/2}_i)}{2} , \qquad \qquad
r_{i-} := \frac{r(h^{\ell+1/2}_{i}) 
       + r(h^{\ell+1/2}_{i-1})}{2} .
\]
This approximation yields a $O(\dx^2)$-accurate matrix approximation
$\tilde{\LL}$ of the operator $\LL$.  Using this, the residual
formulation is written 
\[
\tilde{\LL} \vec{z} = \Vec{RHS}(h^{\ell-1},h^{\ell}) .
\] 
The $N \times N$ matrix $\tilde{\LL}$ is pentadiagonal periodic.  The
righthand side of (\ref{resid_crank_nic}) is discretized analogously.

It follows immediately from the condition (\ref{fd_ss}) for a finite
difference steady state that the above time-stepping scheme preserves
finite-difference steady states.  To ensure this we factored $f$ out
in (\ref{resid}) before doing the finite--difference approximation of
$f(h) h_{xxx} + g(h) h_x$, because the relation $fr=g$ does not hold
in the discrete setting: $f_{i+} r_{i+} \neq g_{i+}$ for example.
Also, by factoring $f$ out we have also isolated the pressure gradient
term, $(h_{xx} + \B/q h^q)_x$, in the equation.

\subsection{Timestepping and accuracy}
\label{acc}

The adaptive timestepping controls the accuracy as follows.  An error
tolerance is set, $\e = 10^{-11}$.  At each time step, we first use
the Crank-Nicolson scheme to compute $h_1$, an approximation of the
solution at time $t + \dt$.  We then take two timesteps with $\dt/2$
to compute $h_2$, another approximation of the solution at time
$t+\dt$.  For some constant $C$, the error is
bounded \cite[\S 5.2]{Iserles} by the difference of $h_1$ and $h_2$:
\[
\|h(\cdot,t+\dt) - h_2 \|_\infty \leq C \| h_1-h_2 \|_\infty.
\]
If $\|h_1-h_2\|_\infty> \e$ then we replace $\dt$ with $\dt/2$ and try again
(without advancing in time).  If $\|h_1-h_2\|_\infty < \e/10$ then we replace
$\dt$ with $2 \dt$ and try again. If $\|h_1-h_2\|_\infty$ lies between $\e$ and $\e/10$ then we take the solution at time $t+\dt$ to be $h_2$.

Adaptive timestepping takes at least three times longer than using the
Crank-Nicolson scheme with a fixed timestep.  For this reason, we
performed all of our exploratory studies using a fixed time-step.
Once we found phenomena of interest, we re-ran using adaptive
time--stepping.  Since the first time--step has $O(\dt^2)$ local
truncation error, rather than $O(\dt^3)$, the adaptive time--stepper
initially refines $\dt$ to meet the tolerance.  Also, most runs had an
initial fast transient (see \S\ref{qm3pp}) which required early
refinement of the timestep.

Admittedly, since we do not know the constant $C$ this error control
is valid only as long as $C$ is not large.  In practice we find that,
after the initial transient, the timestep is rarely reduced, except
near times when the run has to be stopped anyway in order to
point--double.

\subsection{Computing the finite-difference steady state}
\label{compute_id}
Here, we only discuss the case of power law coefficients.  Given $N$
uniformly distributed meshpoints between $0$ and $X=2 \pi$, we seek a
finite-difference steady state $\tilde{h}$ that solves the $N$
equations (\ref{fd_ss}) to the level of round-off error.

We solve the $N$ equations (\ref{fd_ss}) simultaneously using
Newton--Raphson iteration.  To do this, we need a good first guess for
$\tilde{h}$.  In the following, we describe how we find a first guess
and then how we execute the Newton--Raphson iteration.

Given the exponent $q \neq 0$ we first compute the rescaled steady
state $k=k_\alpha$ at the $N$ points $x=P/N, 2P/N, \cdots, P$. As
described in \cite[\S6.1]{LP1}, we do this by viewing the steady state
equation
\begin{equation} \label{k_ss_eqn}
  k_{xx} + \frac{k^q - 1}{q} = 0, \qquad \qquad k(0) = \alpha, \quad
  k_x(0) = 0 ,
\end{equation}
as an initial value problem in $x$.  We verify that $k$ is spectrally
accurate by using a discrete fast Fourier transform to check that the
power-spectrum is fully resolved. Since for numerical purposes, $k$ is
an exact solution of equation (\ref{k_ss_eqn}), in the following we
refer to it as an `analytic steady state'.

Once the analytic steady state $k$ is known, we rescale it to find an
analytic steady state $\h$ of period $2 \pi$, by taking $D=1$ in
the rescaling \cite[eq.\ (7)]{LP3} and defining
\[
\B = \left(\frac{P}{2 \pi}\right)^{\! \! 2q}, \qquad \text{and} \qquad
\h = \left(\frac{1}{\B}\right)^{\! \! 1/q} k.
\]
This gives $\h$ at $N$ meshpoints.  By construction, $\h$
is an analytic steady state for $(\h)_{xx} + (\B \h^q - 1)/q =
0$ with period $2\pi$.

To find the finite-difference steady state $\tilde{h}$ close to $\h$
we need to solve the $N$ equations (\ref{fd_ss}), which we write as
\[
\vec{F}(\tilde{h}) := M \tilde{h} + \vec{V}({\tilde{h}}) = \vec{0}
\]  
where $M$ is a tetradiagonal periodic matrix and
$\vec{V}(\tilde{h})_i$ is a nonlinear function of $\tilde{h}_i$ and
$\tilde{h}_{i+1}$. The Newton--Raphson iteration is
\[
\tilde{h}^{new} = \tilde{h}^{old}
     - (D\vec{F}(\tilde{h}^{old}))^{-1} \vec{F}(\tilde{h}^{old}).
\]
To iterate, one has to solve $D\vec{F}\vec{x} =
\vec{F}(\tilde{h}^{old})$.  We find that $D\vec{F}$ is a singular
matrix of rank $N-1$.  We solve $D\vec{F} \vec{x} =
\vec{F}(\tilde{h}^{old})$ using the singular value decomposition of
$D\vec{F}$ obtained using LAPACK's `dgesvd.f' to solve for $\vec{x}$.
The iteration is then started with the initial guess $\h$ and stopped
when the largest error in the $N$ equations (\ref{fd_ss}) is less than
$10^{-14}$.
%

\vskip 6pt

In a different approach, one might attempt to compute the finite
difference steady state $\tilde{h}$ with a relaxation method, by
computing the finite-difference solution of the evolution equation
$h_t = \Delta h + (\B h^q - 1)/q$ for large $t$. This seems unlikely
to succeed because steady states of this PDE have the same linear
stability properties (with respect to zero--mean perturbations) as the
steady states we are trying to compute, making it extremely difficult
to obtain convergence to linearly unstable steady states.

\vskip 6pt

{\it Note:} For $q=1$, the finite--difference steady states satisfy a
linear problem:
$$
\tilde{h}_{i+2}
- 3 \tilde{h}_{i+1}
+ 3 \tilde{h}_{i}
- \tilde{h}_{i-1}
+ \dx^2 \tilde{\B}
\left( \tilde{h}_{i+1} -  \tilde{h}_{i} \right) = 0 ,
\quad i=1\dots N.
$$
In \S\ref{q1}, we perturb a nontrivial $2\pi$-periodic
finite-difference steady state.  This will be close to an analytic
steady state.  Even analytic steady states are $a + b \cos(\sqrt{\B}
x)$ and are $2\pi$-periodic if $\sqrt{\B}$ is an integer.  Sampling
such a steady state on a uniform mesh gives
$$
\tilde{h}_j = a + b \cos(\sqrt{\B} \; j \Delta x),
$$ 
and one can check that $\tilde{h}$ is a finite--difference steady state provided
$$
\tilde{\B} = 2 \; \frac{1-\cos(\sqrt{\B} \Delta x)}{\Delta x^2} 
           = \B - O(\Delta x^2).
$$
That is, there are nontrivial analytic steady states for a countable
collection of Bond numbers, $\B \in \{ 1^2,2^2,3^2,\ldots \}$, and nontrivial finite difference steady states for a nearby countable set of Bond numbers $\tilde{\B}$.

\subsection{Stopping criteria and issues for singularities} \label{stop}
To test whether to stop the code, we compute the minimum value of
$h^\ell$ at each time-step.  If this minimum is ever less than or
equal to zero, we stop the code.  As discussed in \S\S
\ref{simulations} and \ref{mobility}, we find that the stopping
criterion is indeed often met.  While this suggests a finite-time
singularity, we emphasize that the code has not been written in a way
to carefully resolve such singularities --- the singularity may occur
in infinite time, with the stopping criterion being met in finite time
because of instabilities causing oscillations in the profile.  In
practice, the stopping criterion was always met after the solution
became spectrally unresolved.

The code was designed to preserve the periodic steady states, rather
than being written to preserve positivity.  (There are a number of
different approaches to ensuring positivity: we refer interested
readers to \cite{barrett-blowey, Gruen99,LiyaAndrea}.)  For this
reason, it might be that the code is stopping spuriously.  Also, as
the code has no local mesh refinement, we have to over-resolve much of
the solution in order to resolve the solution where it is tending to
zero.  This over--resolution away from the singular points slows the
computation significantly.  

That being said, we use the code primarily to study nonsingular
behavior.  Some of our results are suggestive of finite--time
singularities; we present them with the above caveats.  To find fine
details of the temporal and spatial scales of the singularities, we
would implement a code that preserves positivity and has an adaptive
spatial mesh.

Finally, we are not using a fully--implicit timestepping scheme; a
basic such scheme would only be $O(\dt)$ but could be made $O(\dt^2)$
using Richardson extrapolation \cite{First}.  It would likely be very
stable, but slow.  For speed, we chose a {\em partially} implicit
scheme and then checked for numerical instabilities when we held the
timestep fixed.  We observed none. Also, it seems unlikely that our
adaptive timestepper was taking small timesteps in order to control
numerical instabilities, since after the initial transient, the
timesteps were refined only when the solution was becoming singular.

\section{Conclusions and Future Directions}
\label{conclusions}

%

We have shown numerically for the evolution equation $h_t = -(h^n
h_{xxx})_x - \B (h^m h_x)_x$ that our linear stability theorems in
\cite{LP2} accurately predict the short-time nonlinear behavior of the
solutions near positive periodic and constant steady states.  We have
found strong evidence for the existence of heteroclinic connections
between steady states, as suggested by our theorems on the energy
levels of steady states \cite{LP3}. We have further observed a
mountain pass scenario in which perturbations of a periodic steady
state relax towards either a droplet or a constant steady state.

All of this suggests that the energy landscape through which the
solutions travel is fairly simple, and that understanding the relative
energy levels of the steady states gives considerable insight into
that landscape.

It is worth recalling that the evolution equation (\ref{evolve})
describes {\em gradient flow} for the energy $\E$ defined in
(\ref{energy_is}), with respect to the following 
weighted $H^{-1}$ inner product. Let $h(x,t)$ be a positive smooth
function that is $X$-periodic in $x$, and for each $t$, define an
inner product on functions $a,b \in C(\T_X)$ having mean value zero by
\[
\langle a,b \rangle := \int_0^X A^\d(x) B^\d(x) h(x,t)^n \, dx ,
\]
where $A,B \in C^2(\T_X)$ satisfy $(h^n A_x)_x=a$ and $(h^n B_x)_x=b$. Then the equation (\ref{evolve}) is equivalent to: 
\[
\frac{\delta \E}{\delta \phi} = - \langle h_t , \phi \rangle \qquad \text{for all $\phi \in C(\T_X)$ having mean value zero.}
\]
Hence the variation of the energy is most negative in the direction
$\phi = h_t$, so that the evolution equation for $h(x,t)$ simply
describes flow by steepest descent on the energy surface of $\E$, with
respect to the inner product $\langle , \rangle$. Note that this inner
product is {\it time-dependent} since it depends on $h^n$, {\it i.e.}
on the solution itself.

The above gradient flow formulation with weighted $H^{-1}$ inner
product was observed by Taylor and Cahn \cite{TC}; their evolution
(7a) contains our equation (\ref{evolve}). For the special case of the
Cahn--Hilliard equation the inner product is unweighted, since the
fourth order term in the equation is linear. Gradient flow ideas for
related equations have been used in \cite{BBW,WB00}, and a
Wasserstein--flow idea in \cite{Otto}.

In \S\ref{mobility} we presented numerical results on the persistence
of heteroclinic connections under changes in the mobility parameters
$n$ and $m$. There we changed $n$ and $m$ in a way that preserved
$q=m-n+1$, and thus preserved the energy $\E$ and also the steady
states (which are critical points of the energy). Hence our change in
mobilities does not change the energy landscape. But it does change
the weight $h^n$ appearing in the inner product $\langle a,b \rangle$
and in the equations for $A$ and $B$, and this is how changing the
mobility affects the evolution. In \S\ref{mobility} we found the {\em
timescale} of the solution changed noticeably in response to changes
in the mobilities, even though the {\em shape} of the solution changed
little.

\vskip 6pt 

Lastly, in \S\ref{mobility} we further investigated {\em
critical} mobility exponents, such as the critical $n$ above which
solutions remain positive for all time (in other words, the critical
exponent for film rupture or pinch-off).  An interesting question for
the future is to find formulas for the critical mobility
exponents. These critical exponents determine important qualitative
features of the evolution and determining them would shed considerable
light not only on the equation (\ref{evolve}) studied here, but also
on related equations that arise from physical models.

\subsection*{{\bf Acknowledgments}}
Laugesen was partially supported by NSF grant number DMS-9970228, and
a grant from the University of Illinois Research Board.  He is
grateful for the hospitality of the Department of Mathematics at
Washington University in St. Louis.

Pugh was partially supported by NSF grant number DMS-9971392, by the
MRSEC Program of the NSF under Award Number DMR-9808595, by the ASCI
Flash Center at the University of Chicago under DOE contract B341495,
and by an Alfred P. Sloan fellowship. The computations were done using
a network of workstations paid for by an NSF SCREMS grant,
DMS-9872029.  Part of the research was conducted while enjoying the
hospitality of the Mathematics Department and the James Franck
Institute of the University of Chicago.

Pugh thanks Todd Dupont and Bastiaan Braams for illuminating
conversations regarding numerical issues.


%
{\sc email contact:}  laugesen@math.uiuc.edu,  mpugh@math.upenn.edu
\end{document}